\documentclass[secthm,seceqn,amsthm,ussrhead,12pt]{amsart}
\usepackage[utf8]{inputenc}
\usepackage[english]{babel}

\usepackage{amssymb,amsmath,amsthm,amsfonts,xcolor,enumerate,hyperref,comment,longtable,cleveref}
\usepackage{soul}
\usepackage{times}
\usepackage{cite}
\usepackage{pdflscape}
\usepackage{ulem}
\usepackage[mathcal]{euscript}
\usepackage{tikz}
\usepackage{hyperref}
\usepackage{cancel}
\usepackage{stmaryrd}

\usetikzlibrary{arrows}

\newcommand{\Co}{{\mathbb C}}

\newtheorem{theorem}{Theorem}
\newtheorem{lemma}[theorem]{Lemma}

\newtheorem{remark}[theorem]{Remark}

\newtheorem{definition}[theorem]{Definition}

\newenvironment{Proof}[1][Proof.]{\begin{trivlist}
\item[\hskip \labelsep {\bfseries #1}]}{\flushright
$\Box$\end{trivlist}}

\usepackage{stmaryrd}
\usepackage{xcolor}

\newcommand{\aut}[1]{\operatorname{\mathrm{Aut}}{#1}}

\newcommand{\orb}{\mathrm{Orb}}

\newcommand{\A}{\mathbf{A}}

\newcommand{\cd}[2]{\mathfrak{CD}^{#1}_{#2}}

\newcommand{\la}{\langle}
\newcommand{\ra}{\rangle}
\newcommand{\La}{\Big\langle}
\newcommand{\Ra}{\Big\rangle}
\newcommand{\D}[2]{\mathbf{D}^{#1}_{#2}}
\newcommand{\Dt}[2]{\Delta_{#1#2}}
\newcommand{\Dl}[2]{[\Delta_{#1#2}]}
\newcommand{\nb}[1]{\nabla_{#1}}

\newcommand{\lb}{\lambda}
\newcommand{\0}{\theta}
\newcommand{\af}{\alpha}
\newcommand{\bt}{\beta}
\newcommand{\gm}{\gamma}

\begin{document}

{\Large\noindent 
The algebraic  classification of nilpotent $\mathfrak{CD}$-algebras
\footnote{The work was partially  supported by  FAPESP 18/00726-5, 18/15712-0;
CNPq 404649/2018-1, 302980/2019-9; RFBR 20-01-00030.
The authors thank  Prof. Dr. Pasha Zusmanovich for   discussions about  $\mathfrak{CD}$-algebras.}
%$^,$\footnote{Corresponding Author: Ivan Kaygorodov -- kaygorodov.ivan@gmail.com}
}

\
  
\ 

   {\bf Ivan Kaygorodov$^{a}$ \&   Mykola Khrypchenko$^{b}$ \\

    \medskip
}

{\tiny

$^{a}$ CMCC, Universidade Federal do ABC, Santo Andr\'e, Brazil

$^{b}$ Departamento de Matem\'atica, Universidade Federal de Santa Catarina, Florian\'opolis, Brazil

\

\ 

\smallskip

   E-mail addresses:

    \smallskip

    Ivan Kaygorodov (kaygorodov.ivan@gmail.com)
    
    Mykola Khrypchenko (nskhripchenko@gmail.com)

}

\

\ 

\

\noindent{\bf Abstract}: 
{\it We give an algebraic  classification of complex
$4$-dimensional nilpotent $\mathfrak{CD}$-algebras. }

\


\ 


\section*{Introduction}

There are many results related to the algebraic and geometric 
classification
of low-dimensional algebras in the varieties of Jordan, Lie, Leibniz and 
Zinbiel algebras;
for algebraic classifications  see, for example, 
\cite{ack, cfk19, degr3, usefi1, ck13, degr2, demir, degr1, lisa, demir, jkk19, hac18, ikm19,  kkk18, kpv19, kv16}; 
for geometric classifications and descriptions of degenerations see, for example, 
\cite{pilar}.
In the present paper, we give algebraic classification of nilpotent $\mathfrak{CD}$-algebras.
This is a new class of non-associative algebras introduced in \cite{ack}.
The idea of the definition of a $\mathfrak{CD}$-algebra  comes from the following property of Jordan and Lie algebras: {\it the commutator of any pair of multiplication operators is a derivation}.
This gives three identities of degree four,    which reduce to only one identity of degree four in the commutative or anticommutative case.
Commutative and anticommutative  $\mathfrak{CD}$-algebras are related to many interesting varieties of algebras.
Thus, anticommutative  $\mathfrak{CD}$-algebras is a generalization of Lie algebras, 
containing the intersection of Malcev and Sagle algebras as a proper subvariety. 
Moreover, the following intersections of varieties coincide:
Malcev and Sagle algebras; 
Malcev and anticommutative  $\mathfrak{CD}$-algebras; 
Sagle and anticommutative  $\mathfrak{CD}$-algebras.
On the other hand, 
the variety of anticommutative  $\mathfrak{CD}$-algebras is a proper subvariety of 
the varieties of binary Lie algebras 
and almost Lie algebras \cite{kz}.
The variety of anticommutative  $\mathfrak{CD}$-algebras  coincides with the intersection of the varieties of binary Lie algebras and almost Lie algebras.
Commutative  $\mathfrak{CD}$-algebras is a generalization of Jordan algebras, 
which is a generalization of associative commutative algebras.
On the other hand, the variety of commutative  $\mathfrak{CD}$-algebras also known as the variety of almost-Jordan algebras (sometimes,  called Lie triple algebras)  \cite{peresi,osborn69,Sidorov_1981} 
and the bigger variety of generalized almost-Jordan algebras \cite{arenas,hl,labra}.
The $n$-ary  version of commutative  $\mathfrak{CD}$-algebras was introduced in a recent paper by 
Kaygorodov, Pozhidaev and Saraiva \cite{kps19}. Commutative  $\mathfrak{CD}$-algebras are also related to assosymmetric algebras  \cite{askar18}.

The algebraic classification of nilpotent algebras will be achieved by the calculation of central extensions of algebras from the same variety which have a smaller dimension.
Central extensions of algebras from various varieties were studied, for example, in \cite{ss78,zusmanovich,kkl18,omirov}.
Skjelbred and Sund \cite{ss78} used central extensions of Lie algebras to classify nilpotent Lie algebras.
Using the same method,  
all non-Lie central extensions of  all $4$-dimensional Malcev algebras \cite{hac16},
all non-associative central extensions of all $3$-dimensional Jordan algebras \cite{ha17},
all anticommutative central extensions of $3$-dimensional anticommutative algebras \cite{cfk182},
all central extensions of $2$-dimensional algebras \cite{cfk18}
and some others were described.
One can also look at the classification of
$3$-dimensional nilpotent algebras \cite{fkkv},
$4$-dimensional nilpotent associative algebras \cite{degr1},
$4$-dimensional nilpotent Novikov algebras \cite{kkk18},
$4$-dimensional nilpotent bicommutative algebras \cite{kpv19},
$4$-dimensional nilpotent commutative algebras in \cite{fkkv},
$4$-dimensional nilpotent assosymmetric algebras in \cite{ikm19},
$4$-dimensional nilpotent noncommutative Jordan algebras in \cite{jkk19},
$4$-dimensional nilpotent terminal algebras \cite{kkp19geo},
$5$-dimensional nilpotent restricted Lie algebras \cite{usefi1},
$5$-dimensional nilpotent associative commutative algebras \cite{krs19},
$5$-dimensional nilpotent Jordan algebras \cite{ha16},
$6$-dimensional nilpotent Lie algebras \cite{degr3, degr2},
$6$-dimensional nilpotent Malcev algebras \cite{hac18},
$6$-dimensional nilpotent Tortkara algebras \cite{gkk,gkks},
$6$-dimensional nilpotent binary Lie algebras \cite{ack},
$6$-dimensional nilpotent anticommutative $\mathfrak{CD}$-algebras \cite{ack},
$6$-dimensional nilpotent anticommutative algebras \cite{kkl19},
$8$-dimensional dual mock-Lie algebras \cite{lisa}.

\ 

\paragraph{\bf Motivation and contextualization} 
Given algebras ${\bf A}$ and ${\bf B}$ in the same variety, we write ${\bf A}\to {\bf B}$ and say that ${\bf A}$ {\it degenerates} to ${\bf B}$, or that ${\bf A}$ is a {\it deformation} of ${\bf B}$, if ${\bf B}$ is in the Zariski closure of the orbit of ${\bf A}$ (under the base-change action of the general linear group). The study of degenerations of algebras is very rich and closely related to deformation theory, in the sense of Gerstenhaber. It offers an insightful geometric perspective on the subject and has been the object of a lot of research.
In particular, there are many results concerning degenerations of algebras of small dimensions in a  variety defined by a set of identities.
One of the main problems of the {\it geometric classification} of a variety of algebras is a description of its irreducible components. In the case of finitely-many orbits (i.e., isomorphism classes), the irreducible components are determined by the rigid algebras --- algebras whose orbit closure is an irreducible component of the variety under consideration. 
The algebraic classification of complex  $4$-dimensional nilpotent $\mathfrak{CD}$-algebras gives  a way to obtain a geometric classification of  all
complex $4$-dimensional nilpotent $\mathfrak{CD}$-algebras, as well as an algebraic classification of all complex $4$-dimensional nilpotent algebras.

 \ 
 
\paragraph{\bf Main result}
The main result of the paper is the complete classification of complex $4$-dimensional nilpotent $\mathfrak{CD}$-algebras.
The full list of non-isomorphic algebras consists of three parts:
\begin{enumerate}
    \item Trivial $\mathfrak{CD}$-algebras were classified in 
    \cite[Theorem 2.1, Theorem 2.3  and Theorem 2.5]{demir} and the only anticommutative 
    trivial $\mathfrak{CD}$-algebra $\mathfrak{CD}_{03}^{4*}$;
    
    \item Terminal (non-Leibniz) non-trivial extensions of the algebra $\cd {3*}{04}(\lambda)$ were classified in \cite[1.4.5. $1$-dimensional central extensions of ${\bf T}^{3*}_{04}$]{kkp19geo};
    
    \item The rest of algebras were found in the present paper.

\end{enumerate}

\newpage
\section{The algebraic classification of nilpotent $\mathfrak{CD}$-algebras}

The class of  $\mathfrak{CD}$-algebras is defined by the 
property that the commutator of any pair of multiplication operators is a derivation; 
namely, an algebra $\mathfrak{A}$ is a $\mathfrak{CD}$-algebra if and only if  
\[ [T_x,T_y]   \in \mathfrak{Der} (\mathfrak{A}),\]
for all $x,y \in \mathfrak{A}$, where  $T_z \in \{ R_z,L_z\}$. Here we use the notation $R_z$ (resp. $L_z$) for the operator of right (resp. left) multiplication in $\mathfrak{A}$. We will denote the variety of $\mathfrak{CD}$-algebras by $\mathfrak{CD}$.
In terms of identities, the class of $\mathfrak{CD}$-algebras is defined by the following three:
\begin{align*}
    ((xy)a)b-((xy)b)a&=((xa)b-(xb)a)y+x((ya)b-(yb)a),\\
    (a(xy))b-a((xy)b)&=((ax)b-a(xb))y+x((ay)b-a(yb)),\\
    a(b(xy))-b(a(xy))&=(a(bx)-b(ax))y+x(a(by)-b(ay)).
\end{align*}

\subsection{Method of classification of nilpotent algebras}
Throughout this paper, we  use the notations and methods described in \cite{ha17,hac16,cfk18}
and adapted for the $\mathfrak{CD}$-algebra case with some modifications. 
From now on, we will give only some important definitions.

Let $({\bf A}, \cdot)$ be a $\mathfrak{CD}$-algebra over $\mathbb C$ and $\mathbb V$ a vector space over the same base field. Then the $\mathbb C$-linear space ${\rm Z}_{\mathfrak{CD}}^{2}\left(
\bf A,\mathbb V \right) $ is defined as the set of all bilinear maps $
\theta :{\bf A} \times {\bf A} \longrightarrow {\mathbb V},$
such that 
\begin{align*}
    \theta((xy)a,b)-\theta((xy)b,a)&=\theta((xa)b-(xb)a,y)+\theta(x,(ya)b-(yb)a),\\
    \theta(a(xy),b)-\theta(a,(xy)b)&=\theta((ax)b-a(xb),y)+\theta(x,(ay)b-a(yb)),\\
    \theta(a,b(xy))-\theta(b,a(xy))&=\theta(a(bx)-b(ax),y)+\theta(x,a(by)-b(ay)).
\end{align*}

 Its elements will be called \textit{cocycles}. For a
linear map $f$ from $\bf A$ to  $\mathbb V$, if we define $\delta f\colon {\bf A} \times
{\bf A} \longrightarrow {\mathbb V}$ by $\delta f \left( x,y \right) =f(xy )$, then $
\delta f\in {\rm Z}_{\mathfrak{CD}}^{2}\left( {\bf A},{\mathbb V} \right) $. Let ${\rm B}^{2}\left(
{\bf A},{\mathbb V}\right) =\left\{ \theta =\delta f\ :f\in {\rm Hom}\left( {\bf A},{\mathbb V}\right) \right\} $.
One can easily check that ${\rm B}^{2}(\bf A,\mathbb V)$ is a linear subspace of $
{\rm Z}_{\mathfrak{CD}}^{2}\left( {\bf A},{\mathbb V}\right) $ whose elements are called \textit{%
coboundaries}. We define the \textit{second cohomology space} $
{\rm H}_{\mathfrak{CD}}^{2}\left( {\bf A},{\mathbb V}\right) $ as the quotient space ${\rm Z}_{\mathfrak{CD}}^{2}%
\left( {\bf A},{\mathbb V}\right) \big/ {\rm  B}^{2}\left( {\bf A},{\mathbb V}\right) $.
% The equivalence class of $
%\theta \in Z^{2}\left( {\bf A},{\mathbb V}\right) $ will be denoted by $\left[
%\theta \right] \in H^{2}\left( {\bf A},{\mathbb V}\right) $.

\bigskip

Let $\operatorname{Aut}\left( {\bf A} \right) $ be the automorphism group of the $\mathfrak{CD}$-algebra $
{\bf A} $ and let $\phi \in \operatorname{Aut} \left( {\bf A}\right) $. For $\theta \in
{\rm Z}_{\mathfrak{CD}}^{2}\left( {\bf A},{\mathbb V}\right) $ define $\phi \theta \left( x,y\right)
=\theta \left( \phi \left( x\right) ,\phi \left( y\right) \right) $. Then $
\phi \theta \in {\rm Z}_{\mathfrak{CD}}^{2}\left( {\bf A},{\mathbb V}\right) $. So, $\operatorname{Aut}\left( {\bf A}\right) $
acts on ${\rm Z}_{\mathfrak{CD}}^{2}\left( {\bf A},{\mathbb V}\right) $. It is easy to verify that $
{\rm B}^{2}\left( {\bf A},{\mathbb V}\right) $ is invariant under the action of $\operatorname{Aut}\left(
{\bf A}\right) $, and thus $\operatorname{Aut}\left( {\bf A}\right) $ acts on $
{\rm H}_{{\mathfrak{CD}}}^{2}\left( {\bf A},{\mathbb V}\right) $.

\bigskip

Let $\bf A$ be a $\mathfrak{CD}$-algebra of dimension $m<n$ over $\mathbb C$, and ${\mathbb V}$ be a $\mathbb C$-vector
space of dimension $n-m$. For any $\theta \in {\rm Z}_{\mathfrak{CD}}^{2}\left(
{\bf A},{\mathbb V}\right) $ define on the linear space ${\bf A}_{\theta }:={\bf A}\oplus {\mathbb V}$ the
bilinear product ``$\left[ -,-\right] _{{\bf A}_{\theta }}$'' by $
\left[ x+x^{\prime },y+y^{\prime }\right] _{{\bf A}_{\theta }}=
 xy +\theta \left( x,y\right) $ for all $x,y\in {\bf A},x^{\prime },y^{\prime }\in {\mathbb V}$.
Then ${\bf A}_{\theta }$ is a $\mathfrak{CD}$-algebra called an {\it $(n-m)$-dimensional central extension} of ${\bf A}$ by ${\mathbb V}$. In fact, ${\bf A_{\theta}}$ is a $\mathfrak{CD}$-algebra \textit{if and only if} $\theta \in {\rm Z}_{\mathfrak{CD}}^2({\bf A}, {\mathbb V})$.

We also call the
set $\operatorname{Ann}(\theta)=\left\{ x\in {\bf A}:\theta \left( x, {\bf A} \right) +\theta( {\bf A},x)=0\right\} $
the {\it{annihilator}} of $\theta $. We recall that the {\it{annihilator}} of an  algebra ${\bf A}$ is defined as
the ideal $\operatorname{Ann}\left( {\bf A} \right) =\left\{ x\in {\bf A}:  x{\bf A}+{\bf A}x=0\right\}$ and observe
 that
$\operatorname{Ann}\left( {\bf A}_{\theta }\right) = \operatorname{Ann}(\theta) \cap \operatorname{Ann}\left( {\bf A}\right)
 \oplus {\mathbb V}.$

\medskip

We have the next  key result:

\begin{lemma}
Let ${\bf A}$ be an $n$-dimensional $\mathfrak{CD}$-algebra such that $dim(\operatorname{Ann}({\bf A}))=m\neq0$. Then there exists, up to isomorphism, a unique $(n-m)$-dimensional $\mathfrak{CD}$-algebra ${\bf A}^{\prime }$ and a bilinear map $\theta \in {\rm Z}_{\mathfrak{CD}}^2({\bf A}, {\mathbb V})$ with $\operatorname{Ann}({\bf A})\cap \operatorname{Ann}(\theta)=0$, where $\mathbb V$ is a vector space of dimension m, such that ${\bf A}\cong {\bf A}^{\prime }_{\theta}$ and $
{\bf A}/\operatorname{Ann}({\bf A})\cong {\bf A}^{\prime }$.
\end{lemma}

\begin{Proof}

Let ${\bf A}^{\prime }$ be a linear complement of $\operatorname{Ann}({\bf A})$ in ${\bf A}$. Define a linear map $P: {\bf A} \longrightarrow {\bf A}^{\prime }$ by $P(x+v)=x$ for $x\in {\bf A}^{\prime }$ and $v\in \operatorname{Ann}({\bf A})$ and define a multiplication on ${\bf A}^{\prime }$ by $[x, y]_{{\bf A}^{\prime }}=P(x y)$ for $x, y \in {\bf A}^{\prime }$.
For $x, y \in {\bf A}$ we have
$$P(xy)=P((x-P(x)+P(x))(y- P(y)+P(y)))=P(P(x) P(y))=[P(x), P(y)]_{{\bf A}^{\prime }}.$$

Since $P$ is a homomorphism, then $P({\bf A})={\bf A}^{\prime }$ is a $\mathfrak{CD}$-algebra and $
{\bf A}/\operatorname{Ann}({\bf A})\cong {\bf A}^{\prime },$ which gives us the uniqueness of ${\bf A}^{\prime }$. Now, define the map $\theta: {\bf A}^{\prime }\times{\bf A}^{\prime }\longrightarrow \operatorname{Ann}({\bf A})$ by $\theta(x,y)=xy-[x,y]_{{\bf A}^{\prime }}$. Then ${\bf A}^{\prime }_{\theta}$ is ${\bf A}$ and therefore $\theta \in {\rm Z}_{\mathfrak{CD}}^2({\bf A}, {\mathbb V})$ and $\operatorname{Ann}({\bf A})\cap \operatorname{Ann}(\theta)=0$.
\end{Proof}

\bigskip

However, in order to solve the isomorphism problem, we need to study the
action of $\operatorname{Aut}\left( {\bf A}\right) $ on ${\rm H}_{\mathfrak{CD}}^{2}\left( {\bf A},{\mathbb V}%
\right) $. To this end, let us fix $e_{1},\ldots ,e_{s}$ a basis of ${\mathbb V}$, and $
\theta \in {\rm Z}_{\mathfrak{CD}}^{2}\left( {\bf A},{\mathbb V}\right) $. Then $\theta $ can be uniquely
written as $\theta \left( x,y\right) =\underset{i=1}{\overset{s}{\sum }}%
\theta _{i}\left( x,y\right) e_{i}$, where $\theta _{i}\in
{\rm Z}_{\mathfrak{CD}}^{2}\left( {\bf A},\mathbb C\right) $. Moreover, $\operatorname{Ann}(\theta)=\operatorname{Ann}(\theta _{1})\cap \operatorname{Ann}(\theta _{2})\cap\ldots \cap \operatorname{Ann}(\theta _{s})$. Further, $\theta \in
{\rm B}^{2}\left( {\bf A},{\mathbb V}\right) $\ if and only if all $\theta _{i}\in {\rm B}^{2}\left( {\bf A},%
\mathbb C\right) $.

\bigskip

\begin{definition}
Let ${\bf A}$ be an algebra and $I$ be a subspace of ${\rm Ann}({\bf A})$. If ${\bf A}={\bf A}_0 \oplus I$ for some subalgebra ${\bf A}_0$ of ${\bf A}$
then $I$ is called an {\it annihilator component} of ${\bf A}$. We say that an algebra is  {\it split} if it has a nontrivial annihilator component.
\end{definition}

It is not difficult to prove (see \cite[%
Lemma 13]{hac16}) that, given a $\mathfrak{CD}$-algebra ${\bf A}_{\theta}$, if we write as
above $\theta \left( x,y\right) =\underset{i=1}{\overset{s}{\sum }}$
 $\theta_{i}\left( x,y\right) e_{i}\in {\rm Z}_{\mathfrak{CD}}^{2}\left( {\bf A},{\mathbb V}\right) $ and we have
$\operatorname{Ann}(\theta)\cap \operatorname{Ann}\left( {\bf A}\right) =0$, then ${\bf A}_{\theta }$ has an
annihilator component if and only if $\left[ \theta _{1}\right] ,\left[
\theta _{2}\right] ,\ldots ,\left[ \theta _{s}\right] $ are linearly
dependent in ${\rm H}_{\mathfrak{CD}}^{2}\left( {\bf A},\mathbb C\right) $.

\bigskip

Let ${\mathbb V}$ be a finite-dimensional vector space over $\mathbb C$. The {\it{%
Grassmannian}} $G_{k}\left( {\mathbb V}\right) $ is the set of all $k$-dimensional
linear subspaces of $ {\mathbb V}$. Let $G_{s}\left( {\rm H}_{\mathfrak{CD}}^{2}\left( {\bf A},\mathbb C%
\right) \right) $ be the Grassmannian of subspaces of dimension $s$ in $
{\rm H}_{\mathfrak{CD}}^{2}\left( {\bf A},\mathbb C\right) $. There is a natural action of $
\operatorname{Aut} \left( {\bf A}\right) $ on $G_{s}\left( {\rm H}_{\mathfrak{CD}}^{2}\left( {\bf A},\mathbb C%
\right) \right) $. Let $\phi \in \operatorname{Aut} \left( {\bf A}\right) $. For $W=\left\langle %
\left[ \theta _{1}\right] ,\left[ \theta _{2}\right] ,\dots,\left[ \theta _{s}%
\right] \right\rangle \in G_{s}\left( {\rm H}_{\mathfrak{CD}}^{2}\left( {\bf A},\mathbb C%
\right) \right) $ define $\phi W=\left\langle \left[ \phi \theta _{1}\right]
,\left[ \phi \theta _{2}\right] ,\dots,\left[ \phi \theta _{s}\right]
\right\rangle $. Then $\phi W\in G_{s}\left( {\rm H}_{\mathfrak{CD}}^{2}\left( {\bf A},\mathbb C \right) \right) $. We denote the orbit of $W\in G_{s}\left(
{\rm H}_{\mathfrak{CD}}^{2}\left( {\bf A},\mathbb C\right) \right) $ under the action of $
\operatorname{Aut} \left( {\bf A}\right) $ by $\mathrm{Orb}\left( W\right) $. Since, given
\begin{equation*}
W_{1}=\left\langle \left[ \theta _{1}\right] ,\left[ \theta _{2}\right] ,\dots,%
\left[ \theta _{s}\right] \right\rangle ,W_{2}=\left\langle \left[ \vartheta
_{1}\right] ,\left[ \vartheta _{2}\right] ,\dots,\left[ \vartheta _{s}\right]
\right\rangle \in G_{s}\left( {\rm H}_{\mathfrak{CD}}^{2}\left( {\bf A},\mathbb C\right)\right),
\end{equation*}%
we easily have that if $W_{1}=W_{2}$, then $\underset{i=1}{\overset{s}{%
\cap }}\operatorname{Ann}(\theta _{i})\cap \operatorname{Ann}\left( {\bf A}\right) =\underset{i=1}{\overset{s}%
{\cap }}\operatorname{Ann}(\vartheta _{i})\cap \operatorname{Ann}\left( {\bf A}\right) $, we can introduce
the set

\begin{equation*}
{\rm T}_{s}\left( {\bf A}\right) =\left\{ W=\left\langle \left[ \theta _{1}\right] ,%
\left[ \theta _{2}\right] ,\dots,\left[ \theta _{s}\right] \right\rangle \in
G_{s}\left( {\rm H}_{\mathfrak{CD}}^{2}\left( {\bf A},\mathbb C\right) \right) :\underset{i=1}{%
\overset{s}{\cap }}\operatorname{Ann}(\theta _{i})\cap \operatorname{Ann}\left( {\bf A}\right) =0\right\},
\end{equation*}
which is stable under the action of $\operatorname{Aut}\left( {\bf A}\right) $.

\medskip

Now, let ${\mathbb V}$ be an $s$-dimensional linear space and let us denote by $
{\rm E}\left( {\bf A},{\mathbb V}\right) $ the set of all {\it non-split $s$-dimensional} central extensions of ${\bf A}$ by
${\mathbb V}.$ We can write
\begin{equation*}
{\rm E}\left( {\bf A},{\mathbb V}\right) =\left\{ {\bf A}_{\theta }:\theta \left( x,y\right) =\underset{%
i=1}{\overset{s}{\sum }}\theta _{i}\left( x,y\right) e_{i}\mbox{
and }\left\langle \left[ \theta _{1}\right] ,\left[ \theta _{2}\right] ,\dots,%
\left[ \theta _{s}\right] \right\rangle \in T_{s}\left( {\bf A}\right) \right\} .
\end{equation*}%
Also, we have the next result, which can be proved the same way as \cite[Lemma 17]{hac16}.

\begin{lemma}
 Let ${\bf A}_{\theta },{\bf A}_{\vartheta }\in {\rm E}\left( {\bf A},{\mathbb V}\right).$ 
 Suppose that $\theta \left( x,y\right) =\underset{i=1}{\overset{s}{\sum }}%
\theta _{i}\left( x,y\right) e_{i}$ and $\vartheta \left( x,y\right) =%
\underset{i=1}{\overset{s}{\sum }}\vartheta _{i}\left( x,y\right) e_{i}$.
Then the $\mathfrak{CD}$-algebras ${\bf A}_{\theta }$ and ${\bf A}_{\vartheta } $ are isomorphic
if and only if 
\[\mathrm{Orb}\left\langle \left[ \theta _{1}\right] ,%
\left[ \theta _{2}\right] ,\dots,\left[ \theta _{s}\right] \right\rangle =%
\mathrm{Orb}\left\langle \left[ \vartheta _{1}\right] ,\left[ \vartheta
_{2}\right] ,\dots,\left[ \vartheta _{s}\right] \right\rangle. \]
\end{lemma}

Thus, there exists a one-to-one correspondence between the set of $\operatorname{Aut}
\left( {\bf A}\right) $-orbits on ${\rm T}_{s}\left( {\bf A}\right) $ and the set of
isomorphism classes of ${\rm E}\left( {\bf A},{\mathbb V}\right) $. Consequently, we have a
procedure that allows us, given a $\mathfrak{CD}$-algebra ${\bf A}^{\prime }$ of
dimension $n$, to construct all non-split central extensions of ${\bf A}^{\prime }.$ This procedure is as follows:

\medskip

{\centerline{\it Procedure}}

\begin{enumerate}
\item For a given [nilpotent] $\mathfrak{CD}$-algebra ${\bf A}^{\prime }$
of dimension $n-s$, determine ${\rm T}_{s}({\bf A}^{\prime })$
and $\operatorname{Aut}({\bf A}^{\prime })$.

\item Determine the set of $\operatorname{Aut}({\bf A}^{\prime })$-orbits on $
{\rm T}_{s}({\bf A}^{\prime })$.

\item For each orbit, construct the $\mathfrak{CD}$-algebra corresponding to one of its
representatives.
\end{enumerate}

\medskip

%%\newpage \subsection{Notations}
Let us introduce the following notations. Let ${\bf A}$ be a $\mathfrak{CD}$-algebra with
a basis $e_{1},e_{2}, \ldots, e_{n}$. Then by $\Delta _{ij}$\ we will denote the
bilinear form
$\Delta _{ij}:{\bf A}\times {\bf A}\longrightarrow \mathbb C$
with $\Delta _{ij}\left( e_{l},e_{m}\right) = \delta_{il}\delta_{jm}.$ 
The set $\left\{ \Delta_{ij}:1\leq i, j\leq n\right\} $ is a basis of the linear space of 
bilinear forms on ${\bf A}$, so every $\theta \in
{\rm Z}_{\mathfrak{CD}}^{2}\left( {\bf A},\bf \mathbb V \right) $ can be uniquely written as $
\theta =\underset{1\leq i,j\leq n}{\sum }c_{ij}\Delta _{ij}$, where $
c_{ij}\in \mathbb C$. We also denote by

$$\begin{array}{lll}
\mathfrak{CD}^{i*}_j& \mbox{the }j\mbox{th }i\mbox{-dimensional nilpotent trivial   $\mathfrak{CD}$-algebra}, \\
\mathfrak{CD}^i_j& \mbox{the }j\mbox{th }i\mbox{-dimensional nilpotent  non-trivial  $\mathfrak{CD}$-algebra}.
%\\
%{\mathfrak{N}}_i& \mbox{the }i\mbox{-dimensional algebra with zero product}. \\
%({\bf A})_{i,j} & \mbox{---}& j\mbox{th }i\mbox{-dimensional central extension of }\bf A. \\
\end{array}$$

%%\newpage \subsection{The algebraic classification of low dimensional nilpotent $\mathfrak{CD}$-algebras}

\subsection{Trivial $\mathfrak{CD}$-algebras}
Recall that the class of $n$-dimensional algebras defined by the identities $(xy)z=0$ and $x(yz)=0$ 
lies in the intersection of all well-known varieties of algebras defined by some family of polynomial identities of degree $3,$ such as Leibniz algebras, Zinbiel algebras, associative, Novikov and many other algebras.
On the other hand, 
all algebras defined by the identities $(xy)z=0$ and $x(yz)=0$ are central extensions of some suitable algebra with zero product.
The list of all  non-anticommutative $4$-dimensional algebras defined by the identities $(xy)z=0$ and $x(yz)=0$  can be found in \cite{demir}.
Note that there is only one $4$-dimensional nilpotent anticommutative algebra with identity  $(xy)z=0.$
Obviously all algebras from this list are  $4$-dimensional nilpotent $\mathfrak{CD}$-algebras.

\subsection{$2$-dimensional nilpotent $\mathfrak{CD}$-algebras}
There is only one non-zero $2$-dimensional nilpotent $\mathfrak{CD}$-algebra:
$$\begin{array}{ll llll}
\cd {2*}{01} &:& e_1 e_1 = e_2.
\end{array}$$

\subsection{$3$-dimensional nilpotent $\mathfrak{CD}$-algebras}
Thanks to \cite{cfk18}, we have the classification of all $3$-dimensional nilpotent algebras.
It is easy to see, that each $3$-dimensional nilpotent algebra is a $\mathfrak{CD}$-algebra.

\begin{longtable}{lllll llll}
$\cd 3{01}$&$:$& $e_1 e_1 = e_2$ & $e_2 e_2=e_3$ \\
\hline
$\cd 3{02}$&$:$& $e_1 e_1 = e_2$ & $e_2 e_1= e_3$ & $e_2 e_2=e_3$ \\
\hline$\cd 3{03}$&$:$& $e_1 e_1 = e_2$ & $e_2 e_1=e_3$ \\
\hline$\cd 3{04}(\lambda)$&$:$& $ e_1 e_1 = e_2$ & $e_1 e_2=e_3$ & $e_2 e_1=\lambda e_3$ \\
\hline$\cd {3*}{01}$&$:$& $e_1 e_1 = e_2$\\
\hline$\cd {3*}{02}$&$:$& $e_1 e_1 = e_3$ &$ e_2 e_2=e_3$ \\
\hline$\cd {3*}{03}$&$:$& $e_1 e_2=e_3$ & $e_2 e_1=-e_3$ \\
\hline$\cd {3*}{04}(\lambda)$&$:$& $e_1 e_1 = \lambda e_3$ & $e_2 e_1=e_3$  & $e_2 e_2=e_3$
\end{longtable}

\subsection{$4$-dimensional nilpotent non-trivial $\mathfrak{CD}$-algebras with $2$-dimensional annihilator }
Thanks to \cite{cfk18}, we have the classification of all $4$-dimensional non-split nilpotent non-trivial algebras with $2$-dimensional annihilator:
\begin{longtable}{lllll llll}
$\cd 4{05}$&:  & $e_1 e_1 = e_2$ & $e_2 e_1=e_4$ & $e_2 e_2=e_3$ \\
\hline$\cd 4{06}$ &: & $ e_1 e_1 = e_2$ & $e_1 e_2=e_4$ & $e_2 e_1=e_3$  \\
\hline$\cd 4{07}(\lambda)$&: & $e_1 e_1 = e_2$ & $e_1 e_2=e_4$ & $e_2 e_1=\lambda e_4$ & $e_2 e_2=e_3$ \\
%\cd 4{04}&:&(\mathfrak{N}_2)_{4,4} &:&  e_1 e_2=e_3, & e_2 e_1=e_4;  \\
%\cd 4{05}&:&(\mathfrak{N}_2)_{4,5} &:& e_1 e_1 = e_3,  & e_2 e_1=e_4; \\
%\cd 4{06}&:&(\mathfrak{N}_2)_{4,6} &:& e_1 e_1 = e_3,  & e_2 e_2=e_4; \\
%\cd 4{07}&:&(\mathfrak{N}_2)_{4,7} &:& e_1 e_1 = e_3, & e_1 e_2=e_3, & e_2 e_1=e_4,  & e_2 e_2=e_3; \\
%\cd 4{08}(\lambda)&:&(\mathfrak{N}_2)_{4,8}(\lambda \in {\mathbb C}_{\geq 0}) &:& e_1 e_1 = e_3,  & e_2 e_1=e_4,  & e_2 e_2=e_3+\lambda e_4; \\
%\cd 4{09}(\lambda)&:&(\mathfrak{N}_2)_{4,9}(\lambda \in {\mathbb C}) &:& e_1 e_1 = e_3, & e_1 e_2=e_4, & e_2 e_1=\lambda e_4.
\end{longtable}

%%\newpage 

\section{The algebraic classification of $4$-dimensional nilpotent $\mathfrak{CD}$-algebras}

\subsection{Automorphism and cohomology groups of $3$-dimensional nilpotent $\mathfrak{CD}$-algebras}

\

{ \tiny
\begin{longtable}{|c|c|c|c|}
\hline 
$\A$ & $\aut\A$ & ${\rm Z}^2_{\mathfrak{CD}}(\A)$ & ${\rm H}^2_{\mathfrak{CD}}(\A)$\\
\hline

$\cd 3{01}$ 
& 
$\begin{pmatrix}
x &    0  &  0\\
0 &  x^2  &  0\\
y &   0  &  x^4
\end{pmatrix}$
&
$\la \Dt 11, \Dt 12, \Dt 21, \Dt 22\ra$
&
$\la \Dl 12, \Dl 21\ra$
\\
\hline

$\cd 3{02}$ 
& 
$\begin{pmatrix}
1 &   0  &  0\\
0 &   1  &  0\\
x &   0  &  1
\end{pmatrix}$
&
$\la \Dt 11, \Dt 12, \Dt 21, \Dt 22\ra$
&
$\la \Dl 12, \Dl 21\ra$
\\
\hline

$\cd 3{03}$ 
& 
$\begin{pmatrix}
x &    0  &  0\\
y &  x^2  &  0\\
z &   xy  &  x^3
\end{pmatrix}$ 
&
$\La 
\begin{array}{l}
     \Dt 11, \Dt 12, \Dt 13-2\Dt 31,\\
     \Dt 21, \Dt 22 
\end{array}
\Ra$
&
$\la \Dl 12, \Dl 22, \Dl 13-2\Dl 31\ra$
\\
\hline

$\cd 3{04}$ 
& 
$\begin{pmatrix}
x &               0  &  0\\
y &             x^2  &  0\\
z &   (\lambda+1)xy  &  x^3
\end{pmatrix}$
&
$\La 
\begin{array}{l}
     \Dt 11, \Dt 12,\Dt 21, \Dt 22,\\
      (\lb-2)\Dt 13-(2\lb-1)\Dt 31

\end{array}
\Ra$
&
$\La 
\begin{array}{l}
(\lb-2)\Dl 13-(2\lb-1)\Dl 31,\\
\Dl 21, \Dl 22
\end{array}
\Ra$
\\
\hline
$\cd {3*}{01}$ 
& 
$\begin{pmatrix}
x &    0  &  0\\
y &  x^2  &  u\\
z &   0  &  v
\end{pmatrix}$ 
&
$\La 
\begin{array}{l}
     \Dt 11, \Dt 12, \Dt 13, \Dt 21, \Dt 22\\
     \Dt 23, \Dt 31, \Dt 32, \Dt 33
\end{array}
\Ra$
&
$\La 
\begin{array}{l}
     \Dl 12, \Dl 13, \Dl 21, \Dl 22\\
     \Dl 23, \Dl 31, \Dl 32, \Dl 33
\end{array}
\Ra$\\
\hline

$\cd {3*}{02}$ 
& 
$\begin{pmatrix}
x &    y  &  0\\
(-1)^{n+1} y &  (-1)^n x  &  0\\
z &   u  &  x^2+y^2
\end{pmatrix}$ 
&
$\La 
\begin{array}{l}
     \Dt 11, \Dt 12, \Dt 13, \Dt 21, \Dt 22\\
     \Dt 23, \Dt 31, \Dt 32, \Dt 33 
\end{array}
\Ra$
&
$\La 
\begin{array}{l}
    \Dl 11, \Dl 12, \Dl 13, \Dl 21,\\
    \Dl 23, \Dl 31, \Dl 32, \Dl 33
\end{array}
\Ra$
\\
\hline

$\cd {3*}{03}$ 
& 
$\begin{pmatrix}
x &    y  &  0\\
z &    u  &  0\\
v &    w  &  xu-yz
\end{pmatrix}$ 
&
$\La 
\begin{array}{l}
     \Dt 11, \Dt 12, \Dt 13, \Dt 21, \Dt 22\\
     \Dt 23, \Dt 31, \Dt 32, \Dt 33 
\end{array}
\Ra$
&
$\La 
\begin{array}{l}
    \Dl 11, \Dl 12, \Dl 13, \Dl 22,\\
    \Dl 23, \Dl 31, \Dl 32, \Dl 33
\end{array}
\Ra$
\\
\hline

$\cd {3*}{04}$ 
& 
$\begin{pmatrix}
x      &      y    &  0\\
-\lb y &    x-y    &  0\\
z      &      u    &  x^2-xy+\lb y^2
\end{pmatrix}$
&
$\La 
\begin{array}{l}
     \Dt 11, \Dt 12, \Dt 13, \Dt 21, \Dt 22\\
     \Dt 23, \Dt 31, \Dt 32, \Dt 33 
\end{array}
\Ra$
&
$\La 
\begin{array}{l}
    \Dl 11, \Dl 12, \Dl 13, \Dl 21,\\
    \Dl 23, \Dl 31, \Dl 32, \Dl 33
\end{array}
\Ra$
\\

\hline
\end{longtable}
}

\vskip0.5cm

Let us denote 
second cohomology of 
Jordan-${\mathfrak{CD}}$-algebras,
Lie-${\mathfrak{CD}}$-algebras and 
terminal-${\mathfrak{CD}}$-algebras
as
${\rm H}^2_{J\mathfrak{CD}},$
${\rm H}^2_{L\mathfrak{CD}}$
and 
${\rm H}^2_{T\mathfrak{CD}}.$
Now we have the following table of cohomology spaces:

\begin{longtable}{|c|c|c|}
\hline 
$\A$ & ${\rm H}^2_{J\mathfrak{CD}}(\A)$ & ${\rm H}^2_{\mathfrak{CD}}(\A)$\\
\hline

$\cd {3*}{01}$ 

&
$\La 
\begin{array}{l}
\Dl 12+\Dl 21, \Dl 13+\Dl 31,\\
\Dl 23+\Dl 32,\Dl 33
\end{array}
\Ra$
&
${\rm H}^2_{J\mathfrak{CD}}(\A)\oplus\la\Dl 21,\Dl 22,\Dl 31,\Dl 32\ra$
\\
\hline

$\cd {3*}{02}$ 
& 
$\La 
\begin{array}{l}
     \Dl 11, \Dl 12+\Dl 21,\\
     \Dl 13+\Dl 31, \Dl 23+\Dl 32
\end{array}
\Ra$
&
${\rm H}^2_{J\mathfrak{CD}}(\A)\oplus\la\Dl 21,\Dl 31,\Dl 32,\Dl 33\ra$
\\
\hline

\multicolumn{3}{|c|}{}
\\

\hline 
$\A$ & ${\rm H}^2_{L\mathfrak{CD}}(\A)$ & ${\rm H}^2_{\mathfrak{CD}}(\A)$\\
\hline

$\cd {3*}{03}$ 

&
$\La 
\begin{array}{l}
\Dl 13-\Dl 31,\Dl 23-\Dl 32
\end{array}
\Ra$
&
${\rm H}^2_{L\mathfrak{CD}}(\A)\oplus\la \Dl 11, \Dl 12, \Dl 13, \Dl 22,
     \Dl 23, \Dl 33\ra$
\\
\hline
\multicolumn{3}{|c|}{}
\\
\hline 
  $\A$  &  ${\rm H}^2_{T\mathfrak{CD}}(\A)$  &  ${\rm H}^2_{\mathfrak{CD}}(\A)$ \\
\hline
  $\cd {3*}{04}$  
&
 $\La 
\begin{array}{l}
    \Dl 11, \Dl 12, \Dl 13, \Dl 21,\\
    \Dl 23, \Dl 31, \Dl 32
\end{array}
\Ra$ 
&
  ${\rm H}^2_{T\mathfrak{CD}}(\A)\oplus\la\Dl 33\ra$ \\
\hline
\end{longtable}

\begin{remark}
It is easy to see that each $4$-dimensional central extension of $\cd 3{01}$ or $\cd 3{02}$ is an algebra with $2$-dimensional annihilator and hence it was found before.
\end{remark}

	\subsection{$1$-dimensional central extensions of $\cd 3{03}$}
	
	Let us use the following notations 
	\begin{align*}
	\nb 1 = \Dl 12, \nb 2 = \Dl 22, \nb 3 = \Dl 13-2\Dl 31.    
	\end{align*}
	Take $\0=\sum_{i=1}^3\af_i\nb i\in {\rm H}^2_{\mathfrak{CD}}(\cd 3{03})$.
	If 
	$$
	\phi=
	\begin{pmatrix}
	x &    0  &  0\\
	y &  x^2  &  0\\
	z &   xy  &  x^3
	\end{pmatrix}\in\aut{\cd 3{03}},
	$$
	then
	$$
	\phi^T\begin{pmatrix}
	0        &  \af_1  & \af_3\\
	0        &  \af_2  & 0\\
	-2\af_3  &     0   & 0
	\end{pmatrix} \phi=
	\begin{pmatrix}
	\af^*     &  \af_1^* & \af_3^*\\
	\af^{**}  &  \af_2^* & 0\\
	-2\af_3^* &       0  & 0
	\end{pmatrix},
	$$
	where
	\[ \begin{array}{lll}
	\af^*_1 = x^2(\af_1x + (\af_2 + \af_3)y), & 	\af^*_2 = \af_2x^4, &	\af^*_3 = \af_3x^4.
	\end{array} \]
 	Hence, $\phi\langle\0\rangle=\langle\0^*\rangle$, where $\0^*=\sum\limits_{i=1}^3 \af_i^*  \nb i.$
We are only interested in elements with $\alpha_3\neq 0.$ 

\begin{enumerate}
    \item $\alpha_2+\alpha_3 \neq 0$. Choosing $y= -\frac{ \alpha_1x}{\alpha_2 + \alpha_3}$ and $x=1,$
    we have the family of representatives of distinct orbits $\la \af \nb 2+\nb 3\ra_{\af\ne-1}.$
    
    \item $\alpha_2=-\alpha_3$ and $\alpha_1 \neq 0$. Choosing $x=\frac{\alpha_1}{\alpha_3},$
    we have the representative $\la \nb 1- \nb 2+\nb 3\ra_.$
    
    \item $\alpha_2=-\alpha_3$ and $\alpha_1 = 0$. Then 
    we have the representative $\la - \nb 2+\nb 3\ra_.$
    \end{enumerate}
Thus, we have the following distinct orbits $\la \af \nb 2+\nb 3\ra$ and $\la\nb 1-\nb 2+\nb 3\ra$. The corresponding algebras are
\begin{longtable}{lllllllllllllllllll}
$\cd 4{08}(\alpha)$ &$:$&
$e_1 e_1 = e_2$ & $e_1e_3=e_4$ & $e_2 e_1=e_3$ & $e_2e_2=\alpha e_4$ & $e_3e_1=-2e_4$ \\

\hline$\cd 4{09}$ &$:$& 
$e_1 e_1 = e_2$ & $e_1e_2=e_4$ & $e_1e_3=e_4$ & $e_2 e_1=e_3$ & $e_2e_2=- e_4$ & $e_3e_1=-2e_4.$

\end{longtable}

  \subsection{$1$-dimensional central extensions of $\cd 3{04}$}

Let us use the following notations 
\begin{align*}
    \nb 1 = (\lb-2)\Dl 13-(2\lb-1)\Dl 31, \nb 2 = \Dl 21, \nb 3 = \Dl 22.    
\end{align*}
Take $\0=\sum_{i=1}^3\af_i\nb i\in {\rm H}^2_{\mathfrak{CD}}(\cd 3{04})$.
If 
$$
\phi=
\begin{pmatrix}
x &               0  &  0\\
y &             x^2  &  0\\
z &   (\lambda+1)xy  &  x^3
\end{pmatrix}\in\aut{\cd 3{04}},
$$
then
$$
\phi^T\begin{pmatrix}
0               &      0  & (\lb-2)\af_1\\
\af_2           &  \af_3  & 0\\
-(2\lb-1)\af_1  &     0   & 0
\end{pmatrix} \phi=
	\begin{pmatrix}
	\af^*                &  \af^{**} & (\lb-2)\af_1^*\\
	\af_2^*+\lb\af^{**}  &  \af_3^*  & 0\\
	-(2\lb-1)\af_1^*     &       0  & 0
\end{pmatrix},
$$
where
\[ \begin{array}{lll}
    \af^*_1 = \af_1x^4, \ & 
    \af^*_2 = x^2((1-\lb)(\af_1(\lb+1)^2+\af_3)y + \af_2x),\ & 
    \af^*_3 = \af_3x^4.
\end{array} \]
Hence, $\phi\langle\0\rangle=\langle\0^*\rangle$, where $\0^*=\sum\limits_{i=1}^3 \af_i^*  \nb i.$ We are only interested in elements with $\alpha_1\neq 0.$ 

\begin{enumerate}
    \item $\lambda=1.$ 
    
    \begin{enumerate} 
    \item 
    $\alpha_2\neq 0$. Then choosing $x=\frac{\alpha_2}{\alpha_1}$ we have the family of representatives $\la \nb 1+\nb 2+\alpha\nb 3\ra$ of distinct orbits.

    \item 
    $\alpha_2= 0$. Then we have the family of representatives $\la \nb 1+\alpha\nb 3\ra$ of distinct orbits.
\end{enumerate}

    \item $\lambda=-1.$ 
    
\begin{enumerate}
    \item 
     $\alpha_3= 0$ and $\alpha_2 \neq  0$. Then choosing $x=\frac{\alpha_2}{\alpha_1}$ we have the  representative $\la \nb 1+\nb 2 \ra$.

    \item 
     $\alpha_3= 0$ and $ \alpha_2 = 0$. Then we have the representative $\la \nb 1  \ra$.

    \item 
    $\alpha_3 \neq 0$. Then choosing $y=-\frac{\alpha_2}{2\alpha_3}$ and $x=1$ we have the family of representatives $\la \nb 1+ \af \nb 3 \ra _{\af \neq 0}$ of distinct orbits.
    
\end{enumerate}
        
    \item $\lambda\neq \pm 1.$ 
    
    \begin{enumerate}
        \item 
    $(\lambda+1)^2 \neq -\frac{\alpha_3}{\alpha_1}$. Then we can suppose that $\alpha_2\neq 0$ and choosing $x=-\frac{(1-\lb)(\af_1(\lb+1)^2+\af_3)}{\alpha_2}$ and $y=1$  we have the family of representatives $\la \nb 1+\alpha\nb 3\ra_{\af \neq -(\lambda+1)^2}$ of distinct orbits.
    
    \item $(\lambda+1)^2 = -\frac{\alpha_3}{\alpha_1}$ and $\alpha_2=0$. Then we have the representative $\la \nb 1-(\lb+1)^2\nb 3\ra.$
    
    \item $(\lambda+1)^2 = -\frac{\alpha_3}{\alpha_1}$ and $\alpha_2 \neq 0$. Then choosing  $x=\frac{\alpha_2}{\alpha_1}$ we have the representative $\la \nb 1+\nb 2-(\lb+1)^2\nb 3\ra.$
          \end{enumerate}  

\end{enumerate}

Summarizing, we have the following representatives of distinct orbits:
$$
\la\nb 1+\nb 2+\af\nb 3\ra_{\lb=1},\ 
\la\nb 1+\nb 2-(\lb+1)^2\nb 3\ra_{\lb\ne 1},\ 
\la\nb 1+\af\nb 3\ra.
$$
The corresponding algebras are:
\begin{longtable}{lllllllllllllllllll}

$\cd 4{10}(\alpha)$ &$:$&
$e_1 e_1 = e_2$ & $e_1 e_2=e_3$ & $e_1e_3=-e_4$ \\
&& $e_2 e_1=  e_3 +e_4$ & $e_2e_2=\alpha e_4$ & $e_3e_1=- e_4$ \\

\hline$\cd 4{11}(\lambda \neq1)$ &$:$& 
$e_1 e_1 = e_2$ & $e_1 e_2=e_3$ & $e_1e_3=(\lambda-2)e_4$ \\&& $e_2 e_1=  \lb e_3+e_4$  &
$e_2e_2=-(\lambda+1)^2 e_4$ &$e_3e_1= (1-2\lb) e_4$ \\

\hline$\cd 4{12}(\alpha, \lambda)$ &$:$& 
$e_1 e_1 = e_2$ & $e_1 e_2=e_3$ & $e_1e_3=(\lambda-2)e_4$ \\&& $e_2 e_1=\lambda e_3$ & $e_2e_2=\alpha e_4$ & $e_3e_1=(1-2\lambda) e_4.$ 
\end{longtable}

\subsection{$1$-dimensional central extensions of $\cd {3*}{01}$}
	Let us use the following notations
	$$ 
	\begin{array}{rclrclrclrcl}
	\nb 1& = &\Dl 12+\Dl 21, &\nb 2& = &\Dl 13+\Dl 31, &\nb 3& = &\Dl 23+\Dl 32, &\nb 4& = &\Dl 33,\\
	\nb 5& = &\Dl 21,        &\nb 6& = &\Dl 22,        &\nb 7& = &\Dl 31,        &\nb 8& = &\Dl 32.    
	\end{array}
	$$
	Take $\0=\sum_{i=1}^8\af_i\nb i\in {\rm H}^2_{\mathfrak{CD}}(\cd {3*}{01})$.
	If 
	$$
	\phi=
	\begin{pmatrix}
	x &    0  &  0\\
	y &  x^2  &  u\\
	z &   0  &  v
	\end{pmatrix}\in\aut{\cd {3*}{01}},
	$$
	then
	$$
	\phi^T
	\begin{pmatrix}
	0          &       \af_1& \af_2\\
	\af_1+\af_5&       \af_6& \af_3\\
	\af_2+\af_7& \af_3+\af_8& \af_4
	\end{pmatrix} 
	\phi=
	\begin{pmatrix}
	\af^*          &         \af^*_1& \af^*_2\\
	\af^*_1+\af^*_5&         \af^*_6& \af^*_3\\
	\af^*_2+\af^*_7& \af^*_3+\af^*_8& \af^*_4	
	\end{pmatrix},
	$$
	where
	\begin{longtable}{ll}
	$\af^*_1 = x^2(\af_1x + (\af_3+\af_8)z + \af_6y)$ &
	$\af^*_2= \af_1ux + \af_2vx + \af_3vy + \af_4vz + \af_6uy  + (\af_3+\af_8)uz$\\
	
	$\af^*_3 = x^2(\af_3v+\af_6u)$ &
	$\af^*_4 = \af_4v^2 + \af_6u^2 +(2\af_3+\af_8)uv$\\

	$\af^*_5 = x^2(\af_5x - \af_8z)$ &
	$\af^*_6 = \af_6x^4$\\
	
	$\af^*_7 = u(\af_5x - \af_8z) + v(\af_7x + \af_8y)$ &
	$\af^*_8 = \af_8vx^2.$
	\end{longtable}
	Hence, $\phi\langle\0\rangle=\langle\0^*\rangle$, where $\0^*=\sum\limits_{i=1}^8 \af_i^*  \nb i$. 
	We are interested in $\0$ such that 
	$(\af_2,\af_3,\af_4,\af_7,\af_8)\ne(0,0,0,0,0)$
	and 
	$(\af_1,\af_3,\af_5,\af_6,\af_8)\ne(0,0,0,0,0).$
	If $(\af_5,\af_6,\af_7,\af_8)=(0,0,0,0)$ we have a Jordan algebra, 
	on the other case it is a non-Jordan $\mathfrak{CD}$-algebra.
	
\begin{enumerate}
		\item $\af_8\ne 0$ and $\alpha_6 \neq 0$. 
		Then choosing   $v=\frac{ \alpha_6x^2}{\alpha_8},$ 
		                $y=-\frac{\af_7x}{\af_8},$ 
		                $z=\frac{\af_5x}{\af_8}$, 
		                $u= -\frac{ \alpha_3 x^2}{\alpha_8},$ 
		                we have $\af^*_3=\af_5^5=\af^*_7=0$, $\af^*_6=\af^*_8$  and
		\begin{longtable}{ll}
		$\af^*_1 = \frac{(\af_1\af_8 + (\af_3+\af_8)\af_5 - \af_6\af_7)x^3}{\af_8} $&
		$\af^*_2 =- \frac{ (\alpha_3^2 \alpha_5+\alpha_3 (\alpha_1+\alpha_5) \alpha_8-\alpha_6 (\alpha_4 \alpha_5+\alpha_2 \alpha_8))x^3}{\af^2_8}$\\
		$\af^*_3 = 0 $& 
		$\af^*_4 = \frac{ \alpha_6 (-\alpha_3^2+\alpha_4 \alpha_6-\alpha_3 \alpha_8)x^4}{\alpha_8^2}$\\
		$\af^*_5 = 0$ & 
		$\af^*_6 = \af_6x^4,$\\
		$\af^*_7 = 0$ & 
		$\af^*_8 = \af_6x^4.$
		\end{longtable}
		
		\begin{enumerate}
				\item $\af_1\af_8 + (\af_3+\af_8)\af_5 - \af_6\af_7\ne 0$. Then choosing $x=\frac{\af_1\af_8 + (\af_3+\af_8)\af_5 - \af_6\af_7}{\af_6\af_8}\ne 0$, we get the family of representatives of distinct orbits $\la\nb 1+\af\nb 2+\bt\nb 4+\nb 6+\nb 8\ra$.
				
				\item $\af_1\af_8 + (\af_3+\af_8)\af_5 - \af_6\af_7=0$. Then $\af_1^*=0$ and 
				$\af^*_2 = \frac{\alpha_6 (\alpha_4 \alpha_5-\alpha_3 \alpha_7+\alpha_2 \alpha_8)x^3}{\af_8^2}$.
				\begin{enumerate}

\item $\alpha_4 \alpha_5-\alpha_3 \alpha_7+\alpha_2 \alpha_8\ne 0$. Then choosing $x=\frac{\alpha_4 \alpha_5 - \alpha_3 \alpha_7 + \alpha_2 \alpha_8}{\alpha_8^2},$ we get the family of representatives of distinct orbits $\la\nb 2+\af\nb 4+\nb 6+\nb 8\ra$.
					
\item $\af_4\af_5-\af_3\af_7+ \af_2\af_8=0$. Then $\af_2^*=0$ and  we get the family of representatives of distinct orbits $\la \af \nb 4+\nb 6+\nb 8\ra$. 

				\end{enumerate}
			\end{enumerate}

			\item $\af_8 \neq 0$ and $\af_6=0$. Then choosing 
			$y = -\frac{ \alpha_7x}{\alpha_8}$
			and $z = \frac{\alpha_5x}{\alpha_8},$ we have
			\begin{longtable}{ll}
$\af^*_1 = \frac{(\af_1\af_8 + (\af_3+\af_8)\af_5)x^3}{\af_8}$ &
$\af^*_2 = \frac{(u \alpha_3 \alpha_5+u (\alpha_1+\alpha_5) \alpha_8+v (\alpha_4 \alpha_5-\alpha_3 \alpha_7+\alpha_2 \alpha_8))x}{\af_8}$ \\
$\af^*_3 = \alpha_3 v x^2 $ & 
$\af^*_4 = v (v \alpha_4+u (2 \alpha_3+\alpha_8))$\\
$\af^*_5 = 0$ & 
$\af^*_6 = 0$ \\
$\af^*_7 = 0$ & 
$\af^*_8 = \alpha_8 v x^2.$ 
\end{longtable}
			\begin{enumerate}

\item $2\af_3+\af_8\ne 0$. Then choosing $u=  -\frac{\alpha_4 v}{2 \alpha_3 + \alpha_8}$ we have $\af^*_4=0$.

				\begin{enumerate}
\item $\af_1\af_8 + (\af_3+\af_8)\af_5\ne 0,$
$2 \alpha_3^2 \alpha_7 \neq \alpha_8 (\alpha_2 \alpha_8-\alpha_1 \alpha_4)+\alpha_3 (\alpha_4 \alpha_5+2 \alpha_2 \alpha_8-\alpha_7 \alpha_8).$

Then choosing 
$x=\frac{\alpha_3 \alpha_4 \alpha_5-2 \alpha_3^2 \alpha_7+2 \alpha_2 \alpha_3 \alpha_8-\alpha_1 \alpha_4 \alpha_8-\alpha_3 \alpha_7 \alpha_8+\alpha_2 \alpha_8^2}{\alpha_8^2(2 \alpha_3+\alpha_8)}$
and 
$v=\frac{(\alpha_3 \alpha_5+(\alpha_1+\alpha_5) \alpha_8) (-2 \alpha_3^2 \alpha_7+\alpha_8 (-\alpha_1 \alpha_4+\alpha_2 \alpha_8)+\alpha_3 (\alpha_4 \alpha_5+2 \alpha_2 \alpha_8-\alpha_7 \alpha_8))}{\alpha_8^4(2 \alpha_3+\alpha_8)}$ 
we get the family of representatives of distinct orbits
    $\la\nb 1+ \nb 2+\af\nb 3+\nb 8\ra_{\af \neq -\frac{1}{2}}.$

\item $\af_1\af_8 + (\af_3+\af_8)\af_5\neq 0,$ 
$2 \alpha_3^2 \alpha_7= \alpha_8 (\alpha_2 \alpha_8-\alpha_1 \alpha_4)+\alpha_3 (\alpha_4 \alpha_5+2 \alpha_2 \alpha_8-\alpha_7 \alpha_8).$
Then choosing 
$x=1$ and $v= \frac{\alpha_3 \alpha_5+(\alpha_1+\alpha_5) \alpha_8}{\alpha_8^2}$
we get the family of representatives of distinct orbits
    $\la\nb 1+  \af\nb 3+\nb 8\ra_{\af \neq -\frac{1}{2}}.$

\item  $\af_1\af_8 + (\af_3+\af_8)\af_5= 0$
and 
$\alpha_4 \alpha_5-\alpha_3 \alpha_7+\alpha_2 \alpha_8 \neq 0.$
Then choosing $x=\frac{\alpha_4 \alpha_5-\alpha_3 \alpha_7+\alpha_2 \alpha_8}{\alpha_8^2}$ and $v=1$
we get the family of representatives of distinct orbits
    $\la\nb 2+  \af\nb 3+\nb 8\ra_{\af \neq -\frac{1}{2}}.$

\item  $\af_1\af_8 + (\af_3+\af_8)\af_5= 0$
and 
$\alpha_4 \alpha_5-\alpha_3 \alpha_7+\alpha_2 \alpha_8 = 0.$
Then 
we get the family of representatives of distinct orbits
    $\la  \af\nb 3+\nb 8\ra_{\af \neq -\frac{1}{2}}.$

\end{enumerate}

				\item $2\af_3+\af_8=0$. Then

\begin{longtable}{lll}
$\af^*_1 = \frac{(2\af_1 + \af_5)x^3}{2}$ & 
\multicolumn{2}{l}{
$\af^*_2 = \frac{(2 v \alpha_4 \alpha_5+u (2 \alpha_1+\alpha_5) \alpha_8+v (2 \alpha_2+\alpha_7) \alpha_8)x}{2\af_8}$} \\
$\af^*_3 = -\frac 12  v x^2 \alpha_8$ &
$\af^*_4 =  \af_4 v^2$ &
$\af^*_8 =\alpha_8 v x^2$ 
				\end{longtable}

				\begin{enumerate}
					\item $2\af_1 + \af_5\neq 0$ and $\alpha_4\neq 0$. 
Choosing 
    $x=\frac{\alpha_4 (2 \alpha_1+\alpha_5)}{2 \alpha_8^2},$
    $u=-\frac{\alpha_4 (2 \alpha_1+\alpha_5) (2 \alpha_4 \alpha_5+(2 \alpha_2+\alpha_7) \alpha_8)}{4 \alpha_8^4}$
    and 
    $v=\frac{\alpha_4 (2 \alpha_1+\alpha_5)^2}{4 \alpha_8^3}$ we get the representative $\la\nb 1-\frac 12\nb 3+\nb 4+\nb 8\ra$. 
    
    					\item $2\af_1 + \af_5\neq 0$ and $\alpha_4= 0$. 
Choosing 
    $x=1$
    $u=- \frac{2 \alpha_2+\alpha_7}{2 \alpha_8}$
    and 
    $v=\frac{2 \alpha_1+\alpha_5}{2 \alpha_8}$ we get the representative 
    $\la\nb 1-\frac 12\nb 3+\nb 8\ra$.

	\item $2\af_1 + \af_5= 0,$ $\alpha_4\neq 0$ and $4 \alpha_1 \alpha_4 \neq (2 \alpha_2 + \alpha_7) \alpha_8.$ Then choosing
$x=\frac{(2 \alpha_2+\alpha_7) \alpha_8-4 \alpha_1 \alpha_4}{2 \alpha_8^2}$
and
$v=\frac{ ((2 \alpha_2+\alpha_7) \alpha_8-4 \alpha_1 \alpha_4)^2}{4 \alpha_4 \alpha_8^3}$
we get the representative 
    $\la\nb 2-\frac 12\nb 3+\nb 4+\nb 8\ra$.

	\item $2\af_1 + \af_5= 0,$ $\alpha_4\neq 0$ and $4 \alpha_1 \alpha_4 = (2 \alpha_2 + \alpha_7) \alpha_8.$ Then choosing
$x=1$
and
$v=\frac{   \alpha_8}{ \alpha_4}$
we get the representative 
    $\la  -\frac 12\nb 3+\nb 4+\nb 8\ra$.

 	\item $2\af_1 + \af_5= 0,$ $\alpha_4= 0$ and $2 \alpha_2 \neq -\alpha_7.$ Then choosing
$x=\frac{2 \alpha_2+\alpha_7}{2 \alpha_8}$
and
$v=1$
we get the representative 
    $\la \nb 2  -\frac 12\nb 3+\nb 8\ra$.

  	\item $2\af_1 + \af_5= 0,$ $\alpha_4= 0$ and $2 \alpha_2 = -\alpha_7.$ Then
we get the representative 
    $\la   -\frac 12\nb 3+\nb 8\ra$.

\end{enumerate}

\end{enumerate}

		\item $\af_8=0$. Then $\af^*_8=0$ and

		\begin{longtable}{llr}
$\af^*_1 = x^2(\af_1x + \af_3z + \af_6y)$ &
\multicolumn{2}{l}{$\af^*_2 = u(\af_1x + \af_3z + \af_6y) + v(\af_2x + \af_3y + \af_4z)$} \\

$\af^*_3 = x^2(\af_3v+\af_6u)$ & 
\multicolumn{2}{l}{$\af^*_4 = \af_4v^2 + \af_6u^2 +2\af_3uv$} \\

$\af^*_5 = \af_5x^3$& 
$\af^*_6 = \af_6x^4$&
$\af^*_7 = x(\af_5u + \af_7v).$\\

		\end{longtable}

		\begin{enumerate}
			\item $\af_5\ne 0$ and $\af_6\ne 0$. Then choosing $u=-\frac{\af_7v}{\af_5}$ and $x=\frac {\af_5}{\af_6}$, $y=-\frac{\alpha_1 \alpha_5+z \alpha_3 \alpha_6}{\alpha_6^2}$, we have $\af_1^*=0$, $\af^*_5=\af^*_6=\frac {\af_5^4}{\af_6^3}\ne 0$, $\af^*_7=0$ and

			\begin{longtable}{lll}
$\af^*_2 = \frac{(\alpha_6 (-z \alpha_3^2+\alpha_2 \alpha_5+z \alpha_4 \alpha_6)-\alpha_1 \alpha_3 \alpha_5)v}{\alpha_6^2}$ &

$\af^*_3 = \frac{ \alpha_5 (\alpha_3 \alpha_5-\alpha_6 \alpha_7)v}{\alpha_6^2}$& 

$\af^*_4 = \left(\alpha_4+\frac{\alpha_7 (\alpha_6 \alpha_7-2 \alpha_3 \alpha_5)}{ \alpha_5^2}\right)v^2.$\\
			\end{longtable}
			
			\begin{enumerate}
				\item $\af_3^2 \ne \af_4\af_6   $. Then choosing $z=\frac{\af_5(\af_2\af_6 - \af_1\af_3)}{\af_6(\af_3^2 - \af_4\af_6)}$, we have $\af^*_2=0$. 
	
				\begin{enumerate}
				
\item $\af_3\af_5\ne \af_6\af_7$. Then choosing 
$v=\frac{\alpha_5^3}{\alpha_6(\alpha_3 \alpha_5 -\alpha_6 \alpha_7)}$
we have the family of representatives of distinct orbits
$\la \nb 3+\af\nb 4+\nb 5+\nb 6\ra_{\af\ne 1}$.

\item $\af_3\af_5=\af_6\af_7.$  
    Then choosing $v=\frac{\alpha_5^2}{\alpha_6\sqrt{(\alpha_4 \alpha_6-\alpha_3^2)}}$ we have the representative $\la\nb 4+\nb 5+\nb 6\ra$.

\end{enumerate}

				\item $\af_3^2 =\af_4\af_6$ Then $\af_4=\frac{\af_3^2}{\af_6}$, so
				\begin{longtable}{lll}
$\af^*_2 = \frac{\af_5}{\af_6^2}(\af_2\af_6 - \af_1\af_3)v$ &
$\af^*_3 = \frac {\af_5}{\af_6^2}(\af_3\af_5-\af_6\af_7)v$ &
$\af^*_4 = \frac{v^2}{\af_5^2\af_6}(\af_3\af_5 - \af_6\af_7)^2.$
				\end{longtable}
				
				\begin{enumerate}

\item $\af_3\af_5 \ne \af_6\af_7.$ Then by choosing 
$v=\frac{\alpha_5^3}{\alpha_6 (\alpha_3 \alpha_5 - \alpha_6 \alpha_7)}$ 
we get the family of representatives of distinct orbits $\la \af  \nb 2+  \nb 3+  \nb 4+\nb 5+\nb 6\ra.$

\item  $\af_3\af_5 = \af_6\af_7$ and $\af_2\af_6 \neq \af_1\af_3.$  
Then by choosing 
$v=\frac{\alpha_5^3}{\alpha_6 (\alpha_2 \alpha_6 - \alpha_1 \alpha_3)}$ 
we get the  representative $\la  \nb 2+ \nb 5+\nb 6\ra.$

\item $\af_3\af_5 =\af_6\af_7$ and $\af_2\af_6 = \af_1\af_3 $. 
Then we have a split extension.
				\end{enumerate}
			\end{enumerate}

				\item $\af_5\ne 0$ and $\af_6=0$. Then choosing $u= -\frac{ \alpha_7 v}{\alpha_5},$ we have

				\begin{longtable}{lll}
$\af^*_1 =  (\af_1x + \af_3z)x^2$ & 
$\af^*_2 = \frac{(x \alpha_2 \alpha_5+y \alpha_3 \alpha_5+z \alpha_4 \alpha_5-x \alpha_1 \alpha_7-z \alpha_3 \alpha_7)v}{\af_5}$ & 
$\af^*_3 = \af_3vx^2$\\

$\af^*_4 = \frac{ (\alpha_4 \alpha_5-2 \alpha_3 \alpha_7)v^2}{\alpha_5}$ & 
$\af^*_5 = \af_5x^3$ & 
$\af^*_7 = 0.$
				\end{longtable}

				\begin{enumerate}

					\item $\af_3\ne 0$. Then choosing 
					$z=-\frac{\af_1x}{\af_3},$  
					$y=\frac{ (\alpha_1 \alpha_4-\alpha_2 \alpha_3)x}{\alpha_3^2}$  and
					$v=\frac{ \alpha_5x}{\alpha_3}$, we get 
					$\af^*_1=\af^*_2=\af^*_7=0$, $\af^*_3=\af^*_5=\af_5x^3$ and $\af^*_4=\frac{\af_5(\af_4\af_5-2\af_3\af_7)x^2}{\af_3^2}$.
					
					\begin{enumerate}
						\item $\af_4\af_5\ne 2\af_3\af_7.$ Then by choosing 
						$x=\frac{\alpha_4 \alpha_5 - 2 \alpha_3 \alpha_7}{\alpha_3^2}$ we get the representative $\la\nb 3+\nb 4+\nb 5\ra$.
						\item $\af_4\af_5=2\af_3\af_7.$ Then we get the representative $\la\nb 3+\nb 5\ra$.
					\end{enumerate}
					
					\item $\af_3=0$. Then $\af^*_3=0$ and
					\begin{longtable}{ll}
					$\af^*_1 = \af_1x^3$ & 
					$\af^*_2 = v \left(z \alpha_4+\left(\alpha_2-\frac{\alpha_1 \alpha_7}{\alpha_5}\right)x\right)$\\
					
					$\af^*_4 = \af_4v^2$ & 
					$\af^*_5 = \af_5x^3$\\

					\end{longtable}
					\begin{enumerate}
						\item $\af_4\ne 0$. Then choosing 
						$x=1,$ $y=0,$ $z=\frac{\af_1\af_7-\af_2\af_5}{\af_4\af_5}$
						and $v=\sqrt{\frac{\af_5}{\af_4}}$, 
						we get the family of representatives of distinct orbits $\la\af\nb 1+\nb 4+\nb 5\ra$. 
						
						\item $\af_4=0$ and $\af_2\af_5 \ne \af_1\af_7$. Then we have the family of representatives of distinct orbits $\la\af\nb 1+\nb 2+\nb 5\ra$.
						
 						\item $\af_4=0$ and $\af_2\af_5=\af_1\af_7$. Then we have a split extension.
					\end{enumerate}
				\end{enumerate}

			\item $\af_5=0$ and $\af_6\ne 0$. Then choosing 
			$u = -\frac{ \alpha_3 v}{\alpha_6}$  and  $y = -\frac{x \alpha_1 + z \alpha_3}{\alpha_6},$ we have
			$\af^*_1=\af^*_3=\af^*_5=\af^*_8=0$ and
			\begin{longtable}{ll}
$\af^*_2=\frac{ (x (\alpha_2 \alpha_6-\alpha_1 \alpha_3)+z (\alpha_4 \alpha_6-\alpha_3^2))x}{\alpha_6}$ &
$\af^*_4=\frac{ (\alpha_4 \alpha_6-\alpha_3^2)v^2}{\alpha_6}$\\
$\af^*_6= \alpha_6x^4$ &
$\af^*_7= \alpha_7 v x$
			\end{longtable}

	\begin{enumerate}
					\item $\af_4\af_6-\af_3^2\ne 0$. Then choosing  $z=\frac{(\af_1\af_3-\af_2\af_6)x}{\af_4\af_6-\af_3^2}$ we have $\af^*_2=0$.
					\begin{enumerate}
						\item $\af_7\ne 0$. Then choosing $x=\frac{\af_7}{\sqrt{\af_4\af_6 - \af_3^2}}$ and $v=\frac{\af_6\af_7^2}{\sqrt{(\af_4\af_6 - \af_3^2)^3}}$ we get the representative $\la\nb 4+\nb 6+\nb 7\ra$.
						\item $\af_7=0$. Then we get the representative $\la\nb 4+\nb 6\ra$.
					\end{enumerate} 

					\item $\af_4\af_6=\af_3^2.$ Then $\af^*_4=0$ and 
					\begin{longtable}{lcr}
$\af^*_2 = \frac{vx}{\af_6}(\af_2\af_6-\af_1\af_3)$ &
$\af^*_6 = \af_6x^4$ &
$\af^*_7 = \af_7vx.$
					\end{longtable}
					
					\begin{enumerate}
						\item $\af_2\af_6-\af_1\af_3\ne 0.$ 
						Then we get the family of representatives $\la \af\nb 2+\nb 6+\nb 7\ra_{\af \ne0}.$
						\item $\af_2\af_6-\af_1\af_3=0$ and $\af_7\ne 0$. Then we get the representative $\la\nb 6+\nb 7\ra$.
						\item $\af_2\af_6-\af_1\af_3=\af_7=0$. Then we get a split extension.
					\end{enumerate}

				\end{enumerate}
				\item $\af_5=\af_6=0$ and $\af_7\neq 0.$ Then $\af^*_5=\af^*_6=0$ and

				\begin{longtable}{lcr}
$\af^*_1 = x^2(\af_1x + \af_3z)$ &
\multicolumn{2}{l}{$\af^*_2 = u(\af_1x + \af_3z) + v(\af_2x + \af_3y + \af_4z)$}\\
$\af^*_3 = \af_3vx^2$ &
$\af^*_4 = v(\af_4v + 2\af_3u)$ &
$\af^*_7 = \af_7vx.$\\
				\end{longtable}
				\begin{enumerate}
					\item $\af_3\ne 0$. Then choosing 
					$x=\frac{\af_7}{\af_3}$
					$y=-\frac{\af_1\af_4- \af_2\af_3}{\af_3^3}$,
					$z=-\frac{\af_1\af_7}{\af_3^2}$,  $u=-\frac{\af_4}{2\af_3}$ and  $v=1$, 
					we get the representative $\la\nb 3+\nb 7\ra$.
					
					\item $\af_3=0$. Then $\af^*_3=0$ and
					\begin{longtable}{ll}
					$\af^*_1 = \af_1x^3$ &
					$\af^*_2 = (\af_1u+\af_2v)x + \af_4vz$\\
					$\af^*_4 = \af_4v^2$ &
					$\af^*_7 = \af_7vx.$
					\end{longtable}
					
					\begin{enumerate}
						\item $\af_4\ne 0$ and $\af_1\ne 0$. Then choosing 
						$x=\frac{\af_7^2}{\af_1\af_4}$,
						$u=0,$
						$z=-\frac{\af_2 \af_7^2}{\af_1\af_4^2}$, 
						$v=\frac {\af_7^3}{\af_1 \af_4^2 }$ and  we get the representative $\la\nb 1+\nb 4+\nb 7\ra$. 
						
						\item $\af_4\ne 0$ and $\af_1=0$. Then choosing 
						$x=1,$
						$y=0,$
						$z=-\frac{\af_2}{\af_4}$,
						$u=0$ and 
						$v=\frac {\af_7}{\af_4},$ 
						we get the representative $\la\nb 4+\nb 7\ra$, which gives a split extension. 
				
						\item $\af_4=0$ and $\af_1\ne 0$. Then choosing 
						$x=1,$
						$y=0,$
						$z=0,$
						$u=-\frac{\af_2}{\af_7}$ and 
						$v=\frac{\af_1}{\af_7},$ 
						we get the representative $\la\nb 1+\nb 7\ra$.
				
						\item $\af_4=\af_1=0$. Then we get the family of representatives of distinct orbits $\la\af\nb 2+\nb 7\ra$, which gives  split extensions.
					\end{enumerate}
				\end{enumerate}
				
				\item $\af_5=\af_6=\af_7=0.$

\begin{enumerate}
    \item $\alpha_3\ne 0.$ Then choosing 
    $x=\af_3^2,$
    $y=\af_1\af_4-\af_2\af_3,$
    $z=-\af_1\af_3,$
    $u=-\frac{\af_4}{2}$ and 
    $v=\af_3,$ we get the representative $\la \nb 3\ra$.

    \item $\alpha_3= 0, \af_1\ne 0$ and $\af_4\ne 0.$ Then choosing 
    $x=\af_4,$
    $y=0,$
    $z=-\af_2,$
    $u=0$ and 
    $v=\sqrt{\af_1}\af_4,$ we get the representative $\la \nb 1 +\nb 4\ra$.

    \item $\alpha_3= 0$ and $\af_1= 0,$ then we have a split extension.

\end{enumerate}
			\end{enumerate}

\end{enumerate} 
 Summarizing, 
 we have the following distinct orbits giving non-Jordan algebras:

\begin{longtable}{lll} 
$\la\nb 1+\nb 2+\af\nb 3+\nb 8\ra_{\af\neq-\frac{1}{2}}$ &
$\la\nb 1+\af\nb 2+\bt\nb 4+\nb 6+\nb 8\ra$ &
$\la\af\nb 1+\nb 2+\nb 5\ra$ \\

$\la\nb 1-\frac12\nb 3+\nb 4+\nb 8\ra$ &
$\la\nb 1+\af\nb 3+\nb 8\ra$ &
$\la\af\nb 1+\nb 4+\nb 5\ra$ \\

$\la\nb 1+\nb 4+\nb 7\ra$&
$\la\nb 1+\nb 7\ra$&
$\la\af \nb 2+\nb 3+\nb 4+\nb 5+\nb 6\ra$ \\

$\la\nb 2-\frac12\nb 3+\nb 4+\nb 8\ra$ &
$\la\nb 2+\af\nb 3+\nb 8\ra$ &
$\la\nb 2+\af\nb 4+\nb 6+\nb 8\ra$ \\

$\la\nb 2+\nb 5+\nb 6\ra$ &
$\la\af\nb 2+\nb 6+\nb 7\ra$ &
$\la\nb 3+\nb 4+\nb 5\ra$ \\

$\la\nb 3+\af\nb 4+\nb 5+\nb 6\ra_{\af\ne 1}$&
$\la-\frac12\nb 3+\nb 4+\nb 8\ra$&
$\la\nb 3+\nb 5\ra$ \\

$\la\nb 3+\nb 7\ra$ &
$\la\af\nb 3+\nb 8\ra$ &
$\la\nb 4+\nb 5+\nb 6\ra$ \\

$\la\nb 4+\nb 6\ra$ &
$\la\nb 4+\nb 6+\nb 7\ra$ &
$\la\af\nb 4+\nb 6+\nb 8\ra$

\end{longtable}
and those giving Jordan algebras: 
$\la \nb 1 +\nb 4\ra$ and $\la \nb 3\ra$. 

A routine verification shows that all these orbits are indeed different.
	\begin{longtable}{lllllllllllllllllll}	

$\cd {4}{13}(\af\neq \frac{1}{2})$&$:$& 
$e_1 e_1 = e_2$ & $e_1e_2=e_4$ & $e_1e_3=e_4$& $e_2e_1=e_4$ \\
&&  $e_2e_3=\alpha e_4$& $e_3e_1=e_4$& \multicolumn{2}{l}{$e_3e_2=(\alpha +1)e_4$}\\

\hline
$\cd {4}{14}(\af, \beta)$&$:$& 
$e_1 e_1 = e_2$ & $e_1 e_2 = e_4$&  $e_1 e_3 = \af e_4$& $e_2 e_1 = e_4$\\ && 
$e_2 e_2 = e_4$& $e_3 e_1 = \af e_4$& $e_3 e_2 = e_4$& $e_3 e_3 =\beta  e_4$& \\

\hline
$\cd {4}{15}(\af)$&$:$& 
$e_1 e_1 = e_2$ &  $e_1 e_2 = \af e_4$&   $e_1 e_3 = e_4$&
\multicolumn{3}{l}{ $e_2 e_1 =(\af+1)  e_4$  $e_3 e_1 = e_4$}\\

\hline
$\cd {4}{16}$&$:$& 
$e_1 e_1 = e_2$ & $e_1e_2=e_4$ & $e_2 e_1 = e_4$\\ && $e_2 e_3 = -\frac{1}{2} e_4$& $e_3 e_2 =\frac{1}{2} e_4$& $e_3 e_3 = e_4$&\\

\hline
$\cd {4}{17}(\af)$&$:$& 
$e_1 e_1 = e_2$ &  $e_1 e_2 = e_4$&  $e_2 e_1 = e_4$&  $e_2 e_3 =\af e_4$&  $e_3 e_2 = (\af+1) e_4$&\\

\hline
$\cd {4}{18}(\af)$&$:$& 
$e_1 e_1 = e_2$ &  $e_1 e_2 =\af  e_4$&   \multicolumn{2}{l}{$e_2 e_1 =(\af+1) e_4$}&   $e_3 e_3 = e_4$&\\

\hline
$\cd {4}{19}$&$:$& 
$e_1 e_1 = e_2$ &   $e_1 e_2 = e_4$&   $e_2 e_1 = e_4$&    $e_3 e_1 = e_4$&
 $e_3 e_3 = e_4$&\\
 
\hline
$\cd {4}{20}$&$:$& 
$e_1 e_1 = e_2$ &    $e_1 e_2 = e_4$&   $e_2 e_1 = e_4$&   $e_3 e_1 = e_4$&\\

\hline
$\cd {4}{21}(\af)$&$:$& 
$e_1 e_1 = e_2$ & $e_1 e_3 =\af e_4$&   $e_2 e_1 = e_4$&  $e_2 e_2 = e_4$&\\
&&  $e_2 e_3 = e_4$&  $e_3 e_1 =\af e_4$&  $e_3 e_2 = e_4$&  $e_3 e_3 = e_4$\\

\hline
$\cd {4}{22}$&$:$& 
$e_1 e_1 = e_2$ &   $e_1 e_3 = e_4$&    $e_2 e_3 = -\frac{1}{2} e_4$\\&&
 $e_3 e_1= e_4$&   $e_3 e_2 =\frac{1}{2} e_4$&   $e_3 e_3 = e_4$&\\

\hline
$\cd {4}{23}(\af)$&$:$& 
  $e_1 e_1 = e_2$&   $e_1 e_3 = e_4$&   $e_2 e_3 = \af e_4$&
  $e_3 e_1 = e_4$&   $e_3 e_2 =(\af +1)e_4$&\\

\hline
$\cd {4}{24}(\af)$&$:$& 
  $e_1 e_1 = e_2$&   $e_1 e_3 = e_4$&    $e_2 e_2 = e_4$\\
  && $e_3 e_1 = e_4$&   $e_3 e_2 = e_4$& $e_3 e_3 =\af  e_4$&\\

\hline
$\cd {4}{25} $&$:$& 
  $e_1 e_1 = e_2$&  $e_1 e_3 = e_4$& $e_2 e_1 = e_4$& $e_3 e_1 = e_4$&    $e_2 e_2 = e_4$&\\

\hline
$\cd {4}{26}(\af)$&$:$& 
  $e_1 e_1 = e_2$&   $e_1 e_3 = \af e_4$&    $e_2 e_2 = e_4$& \multicolumn{2}{l}{$e_3 e_1 =(\af+1) e_4$}&\\

  \hline
$\cd {4}{27}$&$:$& 
  $e_1 e_1 = e_2$&   $e_2 e_1 = e_4$& $e_2 e_3 = e_4$&  $e_3 e_2 = e_4$&   $e_3 e_3 = e_4$&\\

  \hline
$\cd {4}{28}(\af\ne 1)$&$:$& 
  $e_1 e_1 = e_2$&   $e_2 e_1 = e_4$&    $e_2 e_2 = e_4$\\
  && $e_2 e_3 =   e_4$&   $e_3 e_2 =  e_4$&   $e_3 e_3 =\af  e_4$&\\

  \hline
$\cd {4}{29}$&$:$& 
  $e_1 e_1 = e_2$&    $e_2 e_3 =-\frac{1}{2} e_4$&   $e_3 e_2 =\frac{1}{2} e_4$&   $e_3 e_3 = e_4$&\\

  \hline
$\cd {4}{30}$&$:$& 
  $e_1 e_1 = e_2$&   $e_2 e_1 = e_4$& $e_2 e_3 = e_4$&   $e_3 e_2 = e_4$&\\

  \hline
$\cd {4}{31}$&$:$& 
  $e_1 e_1 = e_2$&  $e_2 e_3 = e_4$&     $e_3 e_1 = e_4$&
$e_3 e_2 = e_4$&\\

  \hline
$\cd {4}{32}(\af)$&$:$& 
  $e_1 e_1 = e_2$&   $e_2 e_3 = \af e_4$&   \multicolumn{2}{l}{$e_3 e_2 = (\af+1) e_4$}& \\

  \hline
$\cd {4}{33}$&$:$& 
  $e_1 e_1 = e_2$&   $e_2 e_1 = e_4$&   $e_2 e_2 = e_4$& $e_3 e_3 = e_4$& \\

 \hline
$\cd {4}{34}$&$:$& 
  $e_1 e_1 = e_2$&    $e_2 e_2 = e_4$& $e_3 e_3 = e_4$& \\

 \hline
$\cd {4}{35}$&$:$& 
  $e_1 e_1 = e_2$&    $e_2 e_2 = e_4$&  $e_3 e_1 = e_4$& $e_3 e_3 = e_4$& \\

 \hline
$\cd {4}{36}(\af)$&$:$& 
  $e_1 e_1 = e_2$&     $e_2 e_2 = e_4$&   $e_3 e_2= e_4$& $e_3 e_3 =\af  e_4$&\\

 \hline
$\cd {4}{37}$&$:$& 
  $e_1 e_1= e_2$&    $e_1 e_2= e_4$&   $e_2 e_1= e_4$&    $e_3 e_3= e_4$\\

  \hline
$\cd {4}{38}$&$:$& 
  $e_1 e_1= e_2$&   $e_2 e_3= e_4$&   $e_3 e_2= e_4$&\\

\end{longtable}

	  \subsection{$1$-dimensional central extensions of $\cd {3*}{02}$}
	
	Let us use the following notations 
	$$
	\begin{array}{rclrclrclrcl}
	\nb 1& = &\Dl 11, &\nb 2& = &\Dl 12+\Dl 21, &\nb 3& = &\Dl 13+\Dl 31, &\nb 4& = &\Dl 23+\Dl 32,\\
	\nb 5& = &\Dl 21, &\nb 6& = &\Dl 31,        &\nb 7& = &\Dl 32,        &\nb 8& = &\Dl 33.
	\end{array}
	$$
	Take $\0=\sum_{i=1}^8\af_i\nb i\in {\rm H}^2_{\mathfrak{CD}}(\cd {3*}{02})$.
	If 
	$$
	\phi=
	\begin{pmatrix}
	x &    y  &  0\\
	(-1)^{n+1} y &  (-1)^n x  &  0\\
	z &   u  &  x^2+y^2
	\end{pmatrix}\in\aut{\cd {3*}{02}},
	$$
	then
	$$
	\phi^T\begin{pmatrix}
	\af_1        &  \af_2        & \af_3\\
	\af_2+\af_5  &  0            & \af_4\\
	\af_3+\af_6  &  \af_4+\af_7  & \af_8
	\end{pmatrix} \phi=
	\begin{pmatrix}
	\af_1^*+\af^*    &  \af_2^*        & \af_3^*\\
	\af_2^*+\af_5^*  &  \af^*            & \af_4^*\\
	\af_3^*+\af_6^*  &  \af_4^*+\af_7^*  & \af_8^*
	\end{pmatrix},
	$$
	where
	\begin{longtable}{l}
	$\af^*_1 = \af_1(x^2 - y^2) - 2(-1)^n(2\af_2 + \af_5)xy -(-1)^n(2\af_4 + \af_7)(ux+yz)+ $\\
	\multicolumn{1}{r}{$ (2\af_3 + \af_6)(xz-uy) + \af_8(z^2 - u^2)$}\\
	$\af^*_2 = (-1)^n\af_2x^2  + \af_1xy - (-1)^n(\af_2 + \af_5)y^2 + (\af_3x - (-1)^n\af_4y + \af_8z)u+$ \\
	\multicolumn{1}{r}{$((-1)^n\af_4x + \af_3y + (-1)^n\af_7x + \af_6y)z$}\\
	$\af^*_3 = (\af_3x -(-1)^n\af_4y + \af_8z)(x^2 + y^2)$\\
	$\af^*_4= ((-1)^n\af_4x + \af_3y + \af_8u)(x^2 + y^2)$\\
	$\af^*_5 = (-1)^n\af_5(x^2 + y^2) +u(\af_6x - (-1)^n\af_7y)  -z((-1)^n\af_7x + \af_6y)$\\
	$\af^*_6 = (\af_6x - (-1)^n\af_7y)(x^2 + y^2)$\\
	$\af^*_7 = ((-1)^n\af_7x + \af_6y)(x^2 + y^2)$\\
	$\af^*_8 = \af_8(x^2 + y^2)^2.$
	\end{longtable}
	Hence, $\phi\langle\0\rangle=\langle\0^*\rangle$, where $\0^*=\sum\limits_{i=1}^8 \af_i^*  \nb i$. We are only interested in $\0$ with at least one of the coefficients $\af_5,\af_6,\af_7,\af_8$ different from zero.
	
		\begin{enumerate}
		\item $\af_8\ne 0$. Choosing $u=-\frac{(-1)^n\af_4x + \af_3y}{\af_8}$ and $z=-\frac{\af_3x -(-1)^n\af_4y}{\af_8}$, we have $\af^*_ 3=\af^*_ 4=0$, $\af^*_6,\af^*_7,\af^*_8$ are the same as above and

		\begin{longtable}{l}
		$\af^*_1 = \frac 1{\af_8}((\af_1\af_8 - \af_3(\af_3 + \af_6) + \af_4(\af_4 + \af_7))(x^2-y^2)+$  \\
	\multicolumn{1}{r}{$ 2(-1)^n(\af_3\af_7+\af_4(2\af_3 + \af_6) -\af_8(2\af_2 + \af_5))xy)$}\\
		$\af^*_2 = \frac 1{\af_8}((-1)^n(\af_2\af_8 - \af_3(\af_4+\af_7))x^2 +(\af_1\af_8 - \af_3(\af_3 + \af_6) + \af_4(\af_4 + \af_7))xy+$  \\
	\multicolumn{1}{r}{$ (-1)^n(\af_4(\af_3 + \af_6) - \af_8(\af_2 + \af_5))y^2)$}\\
		$\af^*_5 = \frac{(-1)^n}{\af_8}(\af_3\af_7- \af_4\af_6 + \af_5\af_8)(x^2 + y^2).$
		\end{longtable}
		\begin{enumerate}
			\item $\af_6^2+\af_7^2\ne 0$.
			\begin{enumerate}
				\item $\af_6\ne 0$. Then choosing $y=-\frac{(-1)^n\af_7x}{\af_6}$, we have $x^2+y^2=\frac{\af_6^2+\af_7^2}{\af_6^2}x^2\ne 0$, if $x\ne 0$. This substitution leads to $\af^*_7=0$ and
				\begin{longtable}{ll}
				$\af^*_1 = \af^{**}_1 x^2$ &  $\af^*_2 = \af^{**}_2 x^2$\\
				\multicolumn{2}{l}{$\af^*_5 = \frac{(-1)^n(\af_6^2+\af_7^2)}{\af_6^2\af_8}(\af_3\af_7- \af_4\af_6 + \af_5\af_8)x^2$}\\
				$\af^*_6 = \frac{(\af_6^2+\af_7^2)^2}{\af_6^3}x^3$ & $\af^*_8 = \frac{\af_8(\af_6^2+\af_7^2)^2}{\af_6^4}x^4$
				\end{longtable}
				for some coefficients $\af^{**}_1,\af^{**}_2$ $\in\mathbb{C}$. Choosing $x=\frac {\af_6}{\af_8}$, we get the family of representatives $\la\af\nb 1+\bt\nb 2+\gm\nb 5+\nb 6+\nb 8\ra$. Acting by an automorphism which corresponds to $x=1$, $y=z=u=0$, we can change the signs of $\bt$ and $\gm$.
				\item $\af_6=0$. Then $\af_7\ne 0$. Choosing $x=0$ and $y= (-1)^{n+1}\frac{\af_7}{\af_8}$, we get the same family of representatives found in the previous item.
			\end{enumerate}
			\item $\af_6^2+\af_7^2=0$. Then $(\af^*_6)^2+(\af^*_7)^2=0$, so $\af^*_6=\pm i\af^*_7$.
			\begin{enumerate}
				\item $\af_6\ne 0$. Then $\af_7\ne 0$. Choosing $x=0$ and $y=-\frac {(-1)^n\af_7}{\af_8}$, we have $x^2+y^2=\frac {\af_7^2}{\af_8^2}\ne 0$, and this choice leads to $\af^*_6=\af_8^*\ne 0$, and hence $\af^*_7=\pm i\af_8^*$ (we may make $\af^*_7=i\af_8^*$ or $\af^*_7=-i\af_8^*$ by the appropriate choice of $n$). We thus obtain the family of representatives 
				\[ \la\af\nb 1+\bt\nb 2+\gm\nb 5+\nb 6+i\nb 7+\nb 8\ra.\]
				
				Now, if we apply to $\la\af\nb 1+\bt\nb 2+\gm\nb 5+\nb 6+i\nb 7+\nb 8\ra$ the automorphism which corresponds to $z=u=n=0$ and $x=1-iy$ (with $y\ne -\frac i2$), we will obtain $\la\af^*\nb 1+\bt^*\nb 2+\gm^*\nb 5+\nb 6+i\nb 7+\nb 8\ra$, where
				\begin{longtable}{l}
					$\af^* = \frac 1{(x-iy)^2}\left(\af(x^2 - y^2) - 2(2\bt + \gm)xy\right)$\\
					$\bt^* = \frac 1{(x-iy)^2} \left(\bt x^2  + \af xy - (\bt + \gm)y^2\right)$\\
					$\gm^* = \frac \gm{x-iy}.$
				\end{longtable}
				If $\af\ne 0$ or $\bt\ne 0$, then we may suppose that $\af\ne 0$ and $\bt\ne 0$ by applying a suitable automorphism. So, we have the following cases.
				\begin{enumerate}
					\item $\af\ne 0$, $\bt\ne 0$ and $\gm\ne 0$. The discriminants of the quadratic trinomials in $\af^*$ and $\bt^*$ are
					\begin{longtable}{l}
					$D_{\af^*}=4(\af^2+4\bt^2+4\bt\gm+\gm^2)$\\
					$D_{\bt^*}=\af^2+4\bt^2+4\bt\gm.$
					\end{longtable} 
					Since $\gm\ne 0$, then at least one of the discriminants is nonzero. 
					
$\bullet$					If $D_{\bt^*}\ne 0$, then the trinomial in $\bt^*$ has two distinct roots $x_1=\mu_1y$ and $x_2=\mu_2 y$. Observe that at least one of $\mu_k$ is different from $\pm i$, since the coefficient of $xy$ is different from zero. So, the equation $1-iy=\mu_ky$ has the solution $y=\frac 1{\mu_k+i}\ne -\frac i2$. Choosing this $y$, we will have $\bt^*=0$. So, we obtain the family of representatives \[\la\af\nb 1+\gamma\nb 5+\nb 6+i\nb 7+\nb 8\ra_{\gm\ne 0}.\] 
					If $\af^2+\gm^2=0$, then we have two disjoint subfamilies of representatives of distinct orbits \[\la\af\nb 1+i\af\nb 5+\nb 6+i\nb 7+\nb 8\ra_{\af\ne 0}\] and \[\la\af\nb 1-i\af\nb 5+\nb 6+i\nb 7+\nb 8\ra_{\af\ne 0}.\] 
					
					%{\red If $\af^2+\gm^2\ne 0$, then $\orb\la\af\nb 1+\gamma\nb 5+\nb 6+i\nb 7+\nb 8\ra=\orb\la\af'\nb 1+\gamma'\nb 5+\nb 6+i\nb 7+\nb 8\ra$ if and only if either $(\af',\gm')=(\af,\gm)$, or $(\af',\gm')=\left(-\frac{\gm-i\af}{\gm+i\af}\af,\frac{\gm-i\af}{\gm+i\af}\gm\right)$.}
					
					$\bullet$ $D_{\af^*}\ne 0$, then the trinomial in $\af^*$ has two distinct roots with at least one different from $\pm iy$, since the coefficients of $x^2$ and $y^2$ are not equal. Thus, we obtain the family of representatives $\la\bt\nb 2+\gm\nb 5+\nb 6+i\nb 7+\nb 8\ra_{\gm\ne 0}$. It splits into the family 
					$\la\bt\nb 2+\gm\nb 5+\nb 6+i\nb 7+\nb 8\ra_{\gm\not\in\{0,-2\bt\}}$ and the representative 
					$\la\nb 2-2\nb 5+\nb 6+i\nb 7+\nb 8\ra$.  Observe that for $\gm\ne-2\bt$ we have 
					\[\orb\la\bt\nb 2+\gamma\nb 5+\nb 6+i\nb 7+\nb 8\ra=\orb\la\bt'\nb 2+\gamma'\nb 5+\nb 6+i\nb 7+\nb 8\ra \] if and only if either $(\bt',\gm')=(\bt,\gm)$, or $(\bt',\gm')=\left(\bt+\gm,-\gm\right)$.  Moreover, there exists an automorphism $\phi$, such that
\[\phi\la\af\nb 1+\gm\nb 5+\nb 6+i\nb 7+\nb 8\ra_{\gm\not\in \{0,\pm i\af\} }=\la\bt'\nb2+\gm'\nb 5+\nb 6+i\nb 7+\nb 8\ra_{\gm\not\in\{0,-2\bt\}},\] 
 where  
$\bt'= \frac{\left(\gm +i\af +\sqrt{\af^2 + \gm^2}\right)^2}{4\sqrt{\af^2+\gm^2}},$ 
$\gm'=-\frac{(\gm+i\af) \gm}{\sqrt{\af^2 + \gm^2}}.$
					Hence, the family
					\[\la\af\nb 1+\gm\nb 5+\nb 6+i\nb 7+\nb 8\ra_{\gm\not\in \{0,\pm i\af\} }\]  is a subfamily of 
					\[\la\bt\nb 2+\gm\nb 5+\nb 6+i\nb 7+\nb 8\ra_{\gm\not\in\{0,-2\bt\}},\] 
					so it suffices to choose only the second one.
					
					\item $\af\ne 0$, $\bt\ne 0$ and $\gm=0$. 
					
					$\bullet$ If $\af^2+4\bt^2\ne 0$, then using the same argument as above, we obtain the families of representatives
\[\la\af\nb 1+\nb 6+i\nb 7+\nb 8\ra \mbox{  and  }\la\bt\nb 2+\nb 6+i\nb 7+\nb 8\ra.\] However, there exists an automorphism $\phi$, such that 
\[\phi\la\af\nb 1+\nb 6+i\nb 7+\nb 8\ra=\la\frac{\af i}2\nb 2+\nb 6+i\nb 7+\nb 8\ra,\] so it suffices to choose only the first family. All the representatives of this family belong to different orbits. 

$\bullet$ If $\af^2+4\bt^2=0$, then
					$
					\af^* = \pm 2i\bt \frac {(x\pm iy)^2}{(x-iy)^2},
					\bt^* = \bt\frac {(x\pm iy)^2}{(x-iy)^2},
					$
					so we have two disjoint subfamilies of representatives 
\[\la 2i\af\nb 1+\af\nb 2+\nb 6+i\nb 7+\nb 8\ra_{\af\ne 0}\mbox{ and }\la -2i\af\nb 1+\af\nb 2+\nb 6+i\nb 7+\nb 8\ra_{\af\ne 0}.\] However, all the representatives of the first family belong to the same orbit, from which we will take the representative $\la 2i\nb 1+\nb 2+\nb 6+i\nb 7+\nb 8\ra$. The representatives of the second family belong to different orbits.
					
					\item $\af=\bt=0$ and $\gm\ne 0$. Then
					$
					\af^* = -\frac {2\gm xy}{(x-iy)^2},
					\bt^* = -\frac {\gm y^2}{(x-iy)^2}.
					$
					So, a choice of $x$ and $y$ with $x,y\ne 0$ leads us to the subcase (A).
					
					\item $\af=\bt=\gm=0$. Then we have the representative $\la\nb 6+i\nb 7+\nb 8\ra$ which belongs to the family found above.
				\end{enumerate}

				\item $\af_6=0$. Then $\af_7=0$. We have $\af^*_6=\af^*_7=0$ and
\begin{longtable}{lc}
\multicolumn{2}{l}{$\af^*_1 = \frac 1{\af_8}((\af_1\af_8 - \af_3^2 + \af_4^2)(x^2-y^2)+ 2(-1)^n(2\af_3\af_4 - \af_8(2\af_2 + \af_5))xy)$}\\
\multicolumn{2}{l}{$\af^*_2 = \frac 1{\af_8}((-1)^n(\af_2\af_8 - \af_3\af_4)x^2 +(\af_1\af_8 - \af_3^2 + \af_4^2)xy+ (-1)^n(\af_3\af_4 - \af_8(\af_2 + \af_5))y^2)$}\\
$\af^*_5 = (-1)^n\af_5(x^2 + y^2)$ & $\af^*_8 = \af_8(x^2 + y^2)^2.$
				\end{longtable}
				\begin{enumerate}
					\item $(\af_1\af_8 - \af_3^2 + \af_4^2)^2-4(\af_2\af_8 - \af_3\af_4)(\af_3\af_4 - \af_8(\af_2 + \af_5))\ne 0$. Then the equation $\af_2^*=0$ has two different roots $y_1=\mu_1 x$ and $y_2=\mu_2 x$ for some $\mu_1,\mu_2\in\mathbb{C}$. If $\mu_1^2+1=0$, then $\mu_2^2+1\ne 0$, so we may always choose a root such that $x^2+y^2\ne 0$, whenever $x\ne 0$. This choice leads to the family of representatives \[\la\af\nb 1+\bt\nb 5+\nb 8\ra.\] 
					
					$\bullet$ If $\bt\ne 0$, the applying to $\la\af\nb 1+\bt\nb 5+\nb 8\ra$ an automorphism $\phi$ with $x= \sqrt{\bt}$ and $y=z=u=n=0$, we come to the family $\la\af\nb 1+\nb 5+\nb 8\ra$. Observe that there exists an automorphism $\phi$, such that $\phi\la\af\nb 1+\nb 5+\nb 8\ra=\la\af'\nb 1+\nb 5+\nb 8\ra$ if and only if $\af=(-1)^n\af'$. 
					
					$\bullet$ If $\bt=0$ and $\af\ne 0$, then choosing $x=  \sqrt{\af}$, $y=z=u=0$, we get the representative $\la \nb 1+\nb 8\ra$. 
					
					$\bullet$ If $\bt=\af=0$, then we have the representative $\la \nb 8\ra$.
					
					\item $(\af_1\af_8 - \af_3^2 + \af_4^2)^2-4(\af_2\af_8 - \af_3\af_4)(\af_3\af_4 - \af_8(\af_2 + \af_5))=0$. 
					
					$\bullet$ If $\af_5\ne 0$,  then either $\af_3\af_4 - \af_8(\af_2 + \af_5)\ne 0$, or $\af_2\af_8 - \af_3\af_4\ne 0$. Both of the subcases are similar, so consider, for example, $\af_3\af_4 - \af_8(\af_2 + \af_5)\ne 0$. Then choosing $x=1$ and $y=-\sqrt{\frac{\af_2\af_8 - \af_3\af_4}{\af_3\af_4 - \af_8(\af_2 + \af_5)}}$, we have $x^2+y^2=\frac{-\af_5\af_8}{\af_3\af_4 - \af_8(\af_2 + \af_5)}\ne 0$. This substitution gives $\af^*_2=0$, so we are in the previous subcase. 
					
					%Moreover,
				%	\begin{longtable}{l}
%$\af^*_1 = \frac {\left(2\af_3\af_4 - \af_8(2\af_2 + \af_5) \right)
%\left(\af_1\af_8 - \af_3^2 + \af_4^2-2(-1)^n\sqrt{(\af_2\af_8 - \af_3\af_4)(\af_3\af_4 - \af_8(\af_2 + \af_5))} \right)}{\af_8(\af_3\af_4 - \af_8(\af_2 + \af_5))} $
					%\end{longtable}
					%which is also zero. So, we have the family of representatives $\la\af\nb 5+\nb 8\ra$, with $\af\ne 0$. But this is a subfamily of the family found above.
					
					$\bullet$ If $\af_5=0$ and $\af_3\af_4 - \af_2\af_8\ne 0$, then $\af_1\af_8 - \af_3^2 + \af_4^2\ne 0$. We have $\af_5^*=0$ and
					\begin{longtable}{ll}
$\af^*_1 = \frac {\af_1\af_8 - \af_3^2 + \af_4^2}{\af_8}(x\pm(-1)^niy)^2$ &
$\af^*_2 = \frac {\af_2\af_8 - \af_3\af_4}{\af_8}(x\pm(-1)^ny)^2.$
					\end{longtable}
					Choosing $x=1$ and $y=\mp(-1)^n$, we have $x^2+y^2=2\ne 0$, so we obtain the family of representatives $\la\af\nb 1+\nb 8\ra$, with $\af\ne 0$. This is again a subfamily of the family found above.
					
					$\bullet$ If $\af_5=0$ and $\af_3\af_4 - \af_2\af_8=0$, then $\af_1\af_8 - \af_3^2 + \af_4^2=0$. This gives us the representative $\la\nb 8\ra$, which was found above.
				\end{enumerate}
			\end{enumerate}
		\end{enumerate}

%		\item $\af_8=0$. Then $\af^*_8=0$ and
%		\begin{align*}
%		\af^*_1 &= \af_1(x^2 - y^2) - 2(-1)^n(2\af_2 + \af_5)xy - (-1)^n(2\af_4 + \af_7)(ux+yz)\\
%		&\quad+ (2\af_3 + \af_6)(xz-uy),\\
%		\af^*_2 &= (-1)^n\af_2x^2  + \af_1xy - (-1)^n(\af_2 + \af_5)y^2 + (\af_3x - (-1)^n\af_4y)u\\
%		&\quad + ((-1)^n\af_4x + \af_3y+ (-1)^n \af_7x + \af_6y)z,\\
%		\af^*_3 &= (\af_3x -(-1)^n\af_4y)(x^2 + y^2),\\
%		\af^*_4 &= ((-1)^n\af_4x + \af_3y)(x^2 + y^2),\\
%		\af^*_5 &= (-1)^n\af_5(x^2 + y^2) +u(\af_6x - (-1)^n\af_7y)  -z((-1)^n\af_7x + \af_6y),\\
%		\af^*_6 &= (\af_6x - (-1)^n\af_7y)(x^2 + y^2),\\
%		\af^*_7 &= ((-1)^n\af_7x + \af_6y)(x^2 + y^2).
%		\end{align*}
%		\begin{enumerate}

\item $\af_8=0$ and $\af_6^2+\af_7^2\ne 0$. We may assume that $\af_6 \neq 0$, since otherwise we could apply an automorphism with $y\ne 0$ to make $\af_6\ne 0$. Then choosing $y=-\frac{(-1)^n\af_7x}{\af_6}$ and $u=-\frac{(-1)^n\af_5x}{\af_6}$, we have $x^2+y^2=\frac{\af_6^2+\af_7^2}{\af_6^2}x^2\ne 0$, if $x\ne 0$. This substitution gives $\af^*_5=\af^*_7=0$ and

\begin{longtable}{l}
$\af^*_1 = \frac {((2\af_6\af_7(2\af_2 + \af_5)+\af_1(\af_6^2-\af_7^2)+2\af_5(\af_4\af_6-\af_3\af_7))x+\af_6(\af_6^2+ \af_7^2+2(\af_3\af_6+\af_4\af_7))z)x}{\af_6^2}$\\

$\af^*_2 = \frac {(((-1)^n(-\af_7^2(\af_2 + \af_5) + \af_6(\af_2\af_6 -
				\af_1\af_7)-\af_5(\af_3\af_6 + \af_4\af_7))x + (-1)^n\af_6(\af_4\af_6 - \af_3\af_7)z)x}{\af_6^2}$\\
				
$\af^*_3 = \frac 1{\af_6^3}(\af_6^2+\af_7^2)(\af_3\af_6 + \af_4\af_7)x^3$\\

$\af^*_4 = \frac{(-1)^n}{\af_6^3}(\af_6^2+\af_7^2)(\af_4\af_6 - \af_3\af_7)x^3$\\

$\af^*_6 = \frac 1{\af_6^3}(\af_6^2+\af_7^2)^2x^3.$
				\end{longtable}

				\begin{enumerate}
					\item $\af_4\af_6 - \af_3\af_7\ne 0$. Then choosing the appropriate value of $z$, we have $\af^*_2=0$, so we obtain the family of representatives $\la\af^*\nb 1+\bt^*\nb 3+\gm^*\nb 4+\nb 6\ra$, where $\gm^*\ne 0$. 
					
					\begin{enumerate}
					    \item If $\af^*\ne 0$, then choosing $x= \af^*$, $y=z=u=0$, we get the family $\la \nb 1+\af\nb 3+\bt\nb 4+\nb 6\ra_{\bt\ne 0}$. There exists an automorphism $\phi$, such that \[\phi\la \nb 1+\af\nb 3+\bt\nb 4+\nb 6\ra=\la \nb 1+\af'\nb 3+\bt'\nb 4+\nb 6\ra\] if and only if $(\af,\bt)=(\af',(-1)^n\bt')$. 
					    
					    \item If $\af^*=0$, then we get the family $\la \af\nb 3+\bt\nb 4+\nb 6\ra_{\bt\ne 0}$. There exists an automorphism $\phi$, such that \[\phi\la \af\nb 3+\bt\nb 4+\nb 6\ra=\la \af'\nb 3+\bt'\nb 4+\nb 6\ra\] if and only if $(\af,\bt)=(\af',(-1)^n\bt')$.
					
					\end{enumerate}
				
%					Then $\af^*_4=0$ and
%					\begin{align*}
%					\af^*_1 &= \frac 1{\af_6^2}((2\af_6\af_7(2\af_2 + \af_5) +\af_1(\af_6^2-\af_7^2))x^2+(\af_6^2+ \af_7^2)(2\af_3+\af_6)xz),\\
%					\af^*_2 &= \frac 1{\af_6^4}((-1)^n(-\af_6^2\af_7^2(\af_2 + \af_5) + \af_6^3(\af_2\af_6 -
%					\af_1\af_7)+\af_3\af_5(\af_6^2+\af_7^2))x^2,\\
%					\af^*_3 &= \frac {\af_3}{\af_6^4}(\af_6^2+\af_7^2)^2x^3,\\
%					\af^*_6 &= \frac 1{\af_6^3}(\af_6^2+\af_7^2)^2x^3.
%					\end{align*}
	\item $\af_4\af_6 - \af_3\af_7=0$ and $2\af_3+\af_6\ne 0$. Then the appropriate choice of $z$ gives the family of representatives $\la\af^*\nb 2+\bt^*\nb 3+\nb 6\ra$, where $\bt^*\ne-\frac 12$. 
	\begin{enumerate}
	    \item If $\af^*\ne 0$, then choosing
	    $x=\af^*$, $y=u=n=0$ and $z=0$, we get the family $\la \nb 2+\af\nb 3+\nb 6\ra_{\af\ne-\frac 12}$. 
	    \item If $\af^*=0$, then we get the family $\la \af\nb 3+\nb 6\ra_{\af\ne-\frac 12}$.
	    \end{enumerate}
	All the representatives of these families belong to different orbits.
					
	\item $\af_4\af_6 - \af_3\af_7=0$ and $2\af_3+\af_6=0$. Then $\af_3=-\frac 12\af_6\ne 0$. So, we obtain the family of representatives $\la\af^*\nb 1+\bt^*\nb 2-\frac  12\nb 3+\nb 6\ra$. 
	\begin{enumerate}
	    \item If $\af^*\ne 0$, then choosing 
	    $x=\af^*$, $n=y=u=0$, we get the family $\la\nb 1+\af\nb 2-\frac  12\nb 3+\nb 6\ra$. There exists an automorphism $\phi$, such that 
	    \[\phi\la\nb 1+\af\nb 2-\frac  12\nb 3+\nb 6\ra=\la\nb 1+\af'\nb 2-\frac  12\nb 3+\nb 6\ra\] if and only if $\af=(-1)^n\af'$. 
	    
	    \item If $\af^*=0$ and $\bt^*\ne 0$, then choosing $y=u=n=0$ and  $x=\bt^*,$ we get the representative $\la\nb 2-\frac  12\nb 3+\nb 6\ra$. 
	    
	    \item If $\af^*=\bt^*=0$, then we get the representative $\la-\frac  12\nb 3+\nb 6\ra$.

	\end{enumerate}			 

\end{enumerate}
\item $\af_8=0$ and $\af_6^2+\af_7^2=0$, where $\af_6,\af_7\ne 0$. Choosing $x=0$ and $y=-\frac {(-1)^n}{\sqrt[3]{\af_7}}$, we have $x^2+y^2=\frac 1{\sqrt[3]{\af_7^2}}\ne 0$, and this choice leads to $\af^*_6=1$, and hence $\af^*_7=\pm i$ (we may make $\af^*_7=i$ or $\af^*_7=-i$ by the appropriate choice of $n$). Moreover, choosing $u=-\frac{\af_5 +\af_6\sqrt[3]{\af_7}z}{(-1)^n\sqrt[3]{\af_7^4}}$, we have  $\af^*_3 = \frac {\af_4}{\af_7},$ $\af^*_4 = -\frac {(-1)^n\af_3}{\af_7},$ $\af^*_5=0$ and
				\begin{longtable}{l}
$\af^*_1 = \frac 1{\sqrt[3]{\af_7^5}}(-\af_1\af_7 -\af_5\af_6- 2\af_3\af_5 + 2(-\af_3\af_6 + \af_4\af_7 -\af_6^2)\sqrt[3]{\af_7}z)$\\
$\af^*_2 = -\frac {(-1)^n}{\sqrt[3]{\af_7^5}}(\af_2\af_7 + \af_5\af_7 + \af_4\af_5 + (\af_3\af_7+\af_4\af_6+\af_6\af_7)\sqrt[3]{\af_7}z).$
				\end{longtable}
				\begin{enumerate}
					\item $\af_3\af_7+\af_4\af_6+\af_6\af_7\ne 0$. Then choosing $z=-\frac{\af_2\af_7 + \af_5\af_7 + \af_4\af_5}{(\af_3\af_7+\af_4\af_6+\af_6\af_7)\sqrt[3]{\af_7}}$, we obtain the family of representatives $\la\af\nb 1+\bt\nb 3+\gm\nb 4+\nb 6+i\nb 7\ra$. 
				  If $\af\ne 0$, then applying the automorphism with $z=u=0$ and $y=x\ne 0$, we obtain the family $\la\bt\nb 3+\gm\nb 4+\nb 6+i\nb 7\ra$. Now choosing $z=u=n=0$, we obtain $\la\bt^*\nb 3+\gm^*\nb 4+\nb 6+i\nb 7\ra$, where
					\[	\bt^*=\frac{\bt x-\gm y}{x-iy} \mbox{ and }
						\gm^*=\frac{\gm x+\bt y}{x-iy}.\]
					If $\bt\ne 0$ or $\gm\ne 0$, then we may suppose that $\bt\ne 0$ and $\gm\ne 0$ by applying a suitable automorphism. So, we have the following cases.
					\begin{enumerate}
					
					\item If $\bt\ne 0$, $\gm\ne 0$ and $\bt^2+\gm^2\ne 0$, then choosing $y=-\frac {\gm x}{\bt}$ and $x\ne 0$, we obtain the family $\la\af\nb 3+\nb 6+i\nb 7\ra_{\af\ne 0}$, whose representatives belong to different orbits.
					
					\item If $\bt\ne 0$, $\gm\ne 0$ and $\bt^2+\gm^2=0$, then $\bt=\pm i\gm$. If $\bt=i\gm$, then
					$
					\bt^*=\frac{i\gm(x+iy)}{x-iy},
					\gm^*=\frac{\gm(x+iy)}{x-iy},
					$
					so $\bt^*=i\gm^*$. A suitable choice of $x$ and $y$ leads to the representative $\la i\nb 3+\nb 4+\nb 6+i\nb 7\ra$. If $\bt=-i\gm$, then $\bt^*=\bt$ and $\gm^*=\gm$, so we get the family of representatives $\la -i\af\nb 3+\af\nb 4+\nb 6+i\nb 7\ra_{\af\ne 0}$ of different orbits. 
					\end{enumerate}
					If $\bt=\gm=0$, then we have the representative $\la \nb 6+i\nb 7\ra$.
					\item $\af_3\af_7+\af_4\af_6+\af_6\af_7=0$. If $\af_6=\pm i\af_7$, then $\af_3=\mp i(\af_4+\af_7)$ and   we obtain the representative $\la \af_1^\star\nb 1+\af_2^\star\nb 2+\frac{\af_4}{\af_7}\nb 3\pm(-1)^n i(1+\frac{\af_4}{\af_7})\nb 4+\nb 6\mp(-1)^ni\nb 7\ra$, where
					\begin{align*}
					    \af_1^\star&=-\frac{\pm 2i\af_4\af_5 + \af_1\af_7   \pm i\af_5\af_7}{\sqrt[3]{\af_7^5}},\ \
					    \af_2^\star=-\frac{(-1)^n(\af_4\af_5 + \af_2\af_7 + \af_5\af_7)}{\sqrt[3]{\af_7^5}}.
					\end{align*}
					The appropriate choice of $n$ leads to the family \[\la \af\nb 1+\bt\nb 2+\gm\nb 3-i(1+\gm)\nb 4+\nb 6+i\nb 7\ra.\] 
					\begin{enumerate}
					    \item $(\af,\bt)\ne(0,0)$.
					    We may assume that $\af\ne 0$, since otherwise we could apply an automorphism corresponding to $n=0$, $u=iz$ and $x,y\ne 0$.  Then choosing $n=0$ and $u=iz$ we obtain the representative $\la \af^\star\nb 1+\bt^\star\nb 2+\gm^\star\nb 3-i(1+\gm^\star)\nb 4+\nb 6+i\nb 7\ra$, where
					    \begin{align*}
					        \af^\star=\frac{\af x^2 - 4\bt xy - \af y^2}{(x^2 + y^2)(x - iy)},\ \ \bt^\star=\frac{\bt x^2 + \af xy - \bt y^2}{(x^2 + y^2)(x - iy)},\ \ \gm^\star=\frac{\gm x + i(\gm+1)y}{x - iy}.
					    \end{align*}
					   
					    \begin{enumerate}
 					        \item $\af^2+4\bt^2\ne 0$. Then the equation $\af^\star=0$ has two distinct solutions $y_1=\mu_1x$ and $y_2=\mu_2x$. At least one of them is different from $\pm i$, since the coefficients of $x^2$ and $y^2$ in the nominator of $\af^\star$ are not equal. Choosing this solution we have $x^2+y^2\ne 0$ and $\af^\star=0$. Since the condition $\af^2+4\bt^2\ne 0$ is invariant under the automorphisms preserving the family, we may assume $\af=0$ and $\bt\ne 0$ from the very beginning of this subcase. Then $n=y=0$, $u=iz$ and $x=\bt$ gives the family of representatives 
\[\la \nb 2+\af\nb 3-i(1+\af)\nb 4+\nb 6+i\nb 7\ra.\] Observe that
					        \[\orb\la \nb 2+\af\nb 3-i(1+\af)\nb 4+\nb 6+i\nb 7\ra= 
					        \orb\la \nb 2+\af'\nb 3-i(1+\af')\nb 4+\nb 6+i\nb 7\ra\] 
					        if and only if either $\af'=\af$, or $\af'=-\af-1$.
					    
					        \item $\af^2+4\bt^2=0$. Then we have 4 distinct orbits whose representatives are 
\[\La \nb 1\pm\frac i2\nb 2+\nb 3-2i\nb 4+\nb 6+i\nb 7\Ra \mbox{ and }\La \nb 1\pm\frac i2\nb 2-\frac 12\nb 3-\frac i2\nb 4+\nb 6+i\nb 7\Ra.\]
					    \end{enumerate}
					    \item $(\af,\bt)=(0,0)$. Then we have 2 distinct orbits whose representatives are \[\la \nb 3-2i\nb 4+\nb 6+i\nb 7\ra \mbox{ and }\La -\frac 12\nb 3-\frac i2\nb 4+\nb 6+i\nb 7\Ra.\]
					\end{enumerate}

					%depending on whether the coefficient of $\nb 1$ is equal to zero or not. The representatives of the second family belong to pairwise distinct orbits.  Moreover, if $\af\not\in\{0,-\frac 12,-1\}$, then choosing $\phi$ with $n=0$, $z=-iu$, $x=-\frac{4i\af(\af + 1)^2}{(2\af + 1)^2}$ and $y=\frac{4\af^2(\af + 1)}{(2\af + 1)^2}$, we have $\phi\la \nb 2+\af\nb 3-i(1+\af)\nb 4+\nb 6+i\nb 7\ra=\la \nb 1+\frac{(2\af+1+i)(2\af+1-i)}{8\af(\af + 1)}\nb 2-i\nb 4+\nb 6+i\nb 7\ra$. If $\af\in\{0,-\frac 12,-1\}$, then choosing $n=0$, $z=-iu$ and $(x,y)$ a nonzero solution of $-4xy=(x^2 + y^2)(x - iy)$, we again obtain $\phi$ mapping $\la \nb 2+\af\nb 3-i(1+\af)\nb 4+\nb 6+i\nb 7\ra$ to a representative of the family $\la \nb 1+\af\nb 2+\bt\nb 3-i(1+\bt)\nb 4+\nb 6+i\nb 7\ra$.
					
					%So, choosing $z=-\frac{\af_1\af_7 \mp 2i\af_5(\af_4+\af_7)}{\af_7^2\sqrt[3]{\af_7}}$, we obtain the family of representatives 
					%\[\la\af\nb 2+\bt\nb 3+\gm\nb 4+\nb 6+i\nb 7\ra.\] Applying one more automorphism with $z=u=0$ and $x=y\ne 0$, we have $\af=0$, so we are in the situation of item (A) .
				\end{enumerate}

				\item $\af_8=\af_7=\af_6=0$ and $\af_5\ne 0$. 
				Then  $\af^*_6=\af^*_7=0$ and
\begin{longtable}{l}
$\af^*_1 = \af_1(x^2 - y^2) - 2(-1)^n(2\af_2 + \af_5)xy - 2((-1)^n\af_4x+\af_3y)u+ 2(\af_3x-(-1)^n\af_4y)z$\\

$\af^*_2 = (-1)^n\af_2x^2  + \af_1xy - (-1)^n(\af_2 + \af_5)y^2 + (\af_3x - (-1)^n\af_4y)u + ((-1)^n\af_4x + \af_3y)z$\\

$\af^*_3 = (\af_3x -(-1)^n\af_4y)(x^2 + y^2)$\\
$\af^*_4 = ((-1)^n\af_4x + \af_3y)(x^2 + y^2)$\\
$\af^*_5 = (-1)^n\af_5(x^2 + y^2).$
				\end{longtable}
				\begin{enumerate}
					\item $\af_3^2+\af_4^2\ne 0$ and $\af_4\ne 0$. Then choosing $y=\frac{(-1)^n\af_3x}{\af_4}$, we have $x^2+y^2=\frac{\af_3^2+\af_4^2}{\af_4^2}x^2\ne 0$ if $x\ne 0$. Moreover, $\af^*_3=0$ and
					\begin{longtable}{l}
$\af^*_1 = \frac 1{\af_4^2}((\af_1(\af_4^2 - \af_3^2) - 2\af_3\af_4(2\af_2 + \af_5))x - 2(-1)^n\af_4(\af_3^2+\af_4^2)u)x$\\
$\af^*_2 = \frac{(-1)^n}{\af_4^2}((\af_2\af_4^2  + \af_1\af_3\af_4 - (\af_2 + \af_5)\af_3^2)x + \af_4(\af_3^2+\af_4^2)z)x$\\
$\af^*_4 = \frac{(-1)^n}{\af_4^3}(\af_3^2+\af_4^2)^2x^3$\\
$\af^*_5 = \frac{(-1)^n\af_5}{\af_4^2}(\af_3^2+\af_4^2)x^2.$
\end{longtable}
					Choosing $z=-\frac{(\af_2\af_4^2  + \af_1\af_3\af_4 - (\af_2 + \af_5)\af_3^2)x}{\af_4(\af_3^2+\af_4^2)}$ and $u=\frac{(\af_1(\af_4^2 - \af_3^2) - 2\af_3\af_4(2\af_2 + \af_5))x}{2(-1)^n\af_4(\af_3^2+\af_4^2)}$ and  $x=\frac{\af_4\af_5}{\af_3^2+\af_4^2}\ne 0$, we obtain the representative $\la\nb 4+\nb 5\ra$.
					
					\item $\af_3^2+\af_4^2\ne 0$ and $\af_4=0$. Then $\af_3\ne 0$, so choosing $x=0$, $y=\frac{(-1)^n\af_5}{\af_3}$ and the suitable values of $z$ and $u$, we get the same representative $\la\nb 4+\nb 5\ra$.
					
					\item $\af_3^2+\af_4^2=0$ and $\af_4\ne 0$. e may assume that $\af_4=i\af_3$ using the suitable value of $n$. Then choosing $n=0$ we have $\af^*_4=i\af^*_3$. Moreover, choosing $x=0$ and $y= \frac{i\af_5}{\af_3}$, we obtain $\af_3^*=\af_5^*$.  Then $z=i u - \frac{\alpha_1 \alpha_5}{2 \alpha_3^2}$ gives $\af^*_1=0$ and $\af^*_2=\frac{\alpha_5^2 (-i \alpha_1+2 (\alpha_2+\alpha_5))}{2 \alpha_3^2},$ so we obtain the family of representatives $\la\af\nb 2+\nb 3+i\nb 4+\nb 5\ra$. In fact, only two orbits are distinct, whose representatives are $\la \nb 3+i\nb 4+\nb 5\ra$ and $\la-\frac 12\nb 2+\nb 3+i\nb 4+\nb 5\ra$.
					
					\item $\af_3^2+\af_4^2=0$ and $\af_4=0$. Then $\af_3=0$, so we have $\af_3^*=\af_4^*=0$ and this leads to a split extension.
				\end{enumerate}

\item $\af_8=\af_7=\af_6=\af_5=0.$ We may assume that $\af_3\neq 0.$

\begin{enumerate}
    \item $\af_3^2 + \af_4^2\neq0.$ Choosing 
    $z=-\frac{x (4  \alpha_2 \alpha_3 \alpha_4 + \alpha_1 (\alpha_3^2 - \alpha_4^2))}{2\alpha_3  (\alpha_3^2 +  \alpha_4^2))},$
  $u = \frac{(-1)^n x ( \alpha_1 \alpha_3 \alpha_4-\alpha_2 (\alpha_3^2 - \alpha_4^2))}{ \alpha_3(\alpha_3^2 + \alpha_4^2)}$ and     $y = \frac{(-1)^{1 + n} x \alpha_4}{\alpha_3}$ we have  the representative  $\la\nb 3 \ra$.
  
    \item $\af_3 = \pm i \af_4$. Choosing the appropriate value of $n$, we may assume that $\af_3 = i \af_4$.  
    \begin{enumerate}
        \item $\alpha_1-2 i \alpha_2 \neq 0.$ Choosing 
   $y=z=n=0,$
   $u =\frac{ix \alpha_2}{\alpha_4}$ and 
   $x=\frac{\alpha_1 - 2i\alpha_2}{\alpha_4}$
    we have  the representative  $\la\nb 1+i\nb 3 +\nb 4\ra$.
         \item $\alpha_1-2i\alpha_2=0.$ Choosing 
   $y=z=n=0$ and 
   $u =\frac{ix \alpha_2}{\alpha_4}$ 
    we have  the representative  $\la i\nb 3 +\nb 4 \ra$.
    \end{enumerate}  
    \end{enumerate}
 
\end{enumerate}
 
Summarizing, we obtain the following representatives of different orbits giving non-Jordan $\mathfrak{CD}$-algebras
\begin{longtable}{ll}
%{\red $\la \nb 1+\af\nb 2+\bt\nb 3-i(1+\bt)\nb 4+\nb 6+i\nb 7\ra_{\af\ne\pm\frac i2}$} & 
$\la \nb 1+\frac i2\nb 2+\nb 3-2i\nb 4+\nb 6+i\nb 7\ra$&
$\la \nb 1+\frac i2\nb 2-\frac 12\nb 3-\frac i2\nb 4+\nb 6+i\nb 7\ra$\\
$\la \nb 1-\frac i2\nb 2+\nb 3-2i\nb 4+\nb 6+i\nb 7\ra$&
$\la \nb 1-\frac i2\nb 2-\frac 12\nb 3-\frac i2\nb 4+\nb 6+i\nb 7\ra$\\
$\la\nb 1+\af\nb 2-\frac  12\nb 3+\nb 6\ra$&
$\la\af\nb 1+\bt\nb 2+\gm\nb 5+\nb 6+\nb 8\ra$\\
$\la 2i\nb 1+\nb 2+\nb 6+i\nb 7+\nb 8\ra$&
$\la -2i\af\nb 1+\af\nb 2+\nb 6+i\nb 7+\nb 8\ra_{\af\ne 0}$\\
$\la \nb 1+\af\nb 3+\bt\nb 4+\nb 6\ra_{\bt\ne 0}$&
$\la\af\nb 1+i\af\nb 5+\nb 6+i\nb 7+\nb 8\ra_{\af\ne 0}$\\
$\la\af\nb 1-i\af\nb 5+\nb 6+i\nb 7+\nb 8\ra_{\af\ne 0}$&
$\la\af\nb 1+\nb 5+\nb 8\ra$\\ 
$\la\af\nb 1+\nb 6+i\nb 7+\nb 8\ra$&
$\la \nb 1+\nb 8\ra$\\
$\la-\frac 12\nb 2+\nb 3+i\nb 4+\nb 5\ra$&
$\la \nb 2+\af\nb 3-i(1+\af)\nb 4+\nb 6+i\nb 7\ra$\\
$\la \nb 2+\af\nb 3+\nb 6\ra$&
$\la\nb 2-2\nb 5+\nb 6+i\nb 7+\nb 8\ra$\\
$\la\af\nb 2+\bt\nb 5+\nb 6+i\nb 7+\nb 8\ra_{\bt\not\in\{0,-2\af\}}$& 
$\la\nb 3+i\nb 4+\nb 5\ra$\\
$\la \af\nb 3+\bt\nb 4+\nb 6\ra_{\bt\ne 0}$&  
$\la i\nb 3+\nb 4+\nb 6+i\nb 7\ra$\\
$\la -i\af\nb 3+\af\nb 4+\nb 6+i\nb 7\ra$& 
$\la \nb 3-2i\nb 4+\nb 6+i\nb 7\ra$\\
$\la -\frac 12\nb 3-\frac i2\nb 4+\nb 6+i\nb 7\ra$&
$\la \af\nb 3+\nb 6\ra$\\
$\la\af\nb 3+\nb 6+i\nb 7\ra_{\af\ne 0}$& 
$\la\nb 4+\nb 5\ra$\\
$\la \nb 8\ra.$ &

\end{longtable}
and those giving Jordan algebras: $\la\nb 1+i\nb 3 +\nb 4\ra$, $\la i\nb 3 +\nb 4\ra$, $\la \nb 3 \ra.$

\begin{longtable}{lllllllllllllllllll}
%{\red $\cd 4{41}(\af,\bt)$} &:&  
%$e_1 e_1 = e_3+e_4$ & $e_1e_2=\af e_4$ & $e_1e_3=\bt e_4$ & $e_2e_1=\af e_4$\\
%$\af\ne\pm\frac i2$&& $ e_2 e_2=e_3$ & $e_2e_3=-i(1+\bt)e_4$ & $e_3e_1=(\bt+1) e_4$ & $e_3e_2=-i\bt e_4$\\ 

%\hline
$\cd 4{39}$ &:&  
$e_1 e_1 = e_3+e_4$ & $e_1e_2=\frac i2 e_4$ & $e_1e_3=e_4$ & $e_2e_1=\frac i2 e_4$\\
&& $ e_2 e_2=e_3$ & $e_2e_3=-2ie_4$ & $e_3e_1=2e_4$ & $e_3e_2=-ie_4$\\

\hline
$\cd 4{40}$ &:&  
$e_1 e_1 = e_3+e_4$ & $e_1e_2=\frac i2 e_4$ & $e_1e_3=-\frac 12e_4$ & $e_2e_1=\frac i2 e_4$\\
&& $ e_2 e_2=e_3$ & $e_2e_3=-\frac i2e_4$ & $e_3e_1=\frac 12e_4$ & $e_3e_2=\frac i2e_4$\\

\hline
$\cd 4{41}$ &:&  
$e_1 e_1 = e_3+e_4$ & $e_1e_2=-\frac i2 e_4$ & $e_1e_3=e_4$ & $e_2e_1=-\frac i2 e_4$\\
&& $ e_2 e_2=e_3$ & $e_2e_3=-2ie_4$ & $e_3e_1=2e_4$ & $e_3e_2=-ie_4$\\

\hline
 $\cd 4{42}$ &:&  
$e_1 e_1 = e_3+e_4$ & $e_1e_2=-\frac i2 e_4$ & $e_1e_3=-\frac 12e_4$ & $e_2e_1=-\frac i2 e_4$\\
&& $ e_2 e_2=e_3$ & $e_2e_3=-\frac i2e_4$ & $e_3e_1=\frac 12e_4$ & $e_3e_2=\frac i2e_4$\\

\hline
$\cd 4{43}(\af)$ &$:$&  
$e_1 e_1 = e_3 + e_4$     & $e_1 e_2 = \af e_4$          & $e_1 e_3 = -\frac 12 e_4$ \\
&&   $e_2 e_1 = \af e_4$      & $e_2 e_2 = e_3$           & $e_3 e_1 = \frac 12 e_4$\\

\hline
$\cd 4{44}(\af,\bt,\gm)$ &: &  
$e_1 e_1 = e_3 + \af e_4$   & $e_1 e_2 = \bt e_4$ & $e_2 e_1 = (\bt+\gm) e_4$ \\
&& $e_2 e_2 = e_3$ & $e_3 e_1 = e_4$             & $e_3 e_3 = e_4$\\

\hline
$\cd 4{45}$ &: &  
$e_1 e_1 = e_3 + 2i e_4$  & $e_1 e_2 = e_4$ & $e_2 e_1 = e_4$& $e_2 e_2 = e_3$\\
&&  $e_3 e_1 = e_4$           & $e_3 e_2 = i e_4$       & $e_3 e_3 = e_4$\\

\hline
$\cd 4{46}(\af)$ &: &
$e_1 e_1 = e_3 - 2i\af e_4$ & $e_1 e_2 = \af e_4$  & $e_2 e_1 = \af e_4$ & $e_2 e_2 = e_3$ \\
$\af\ne$ 0& & $e_3 e_1 = e_4$              & $e_3 e_2 = i e_4$       & $e_3 e_3 = e_4$\\

\hline
$\cd 4{47}(\af,\bt)$ &: & 
$e_1 e_1 = e_3 + e_4$  & $e_1 e_3 = \af e_4$  & $e_2 e_2 = e_3$\\ 
$\bt\ne 0$ && $e_2 e_3 = \bt e_4$  & $e_3 e_1 =(\af+1) e_4$        & $e_3 e_2 = \bt e_4$ \\

\hline		
$\cd 4{48}(\af)$ &: & 
 $e_1 e_1 = e_3 + \af e_4$ &  $e_2 e_1 = i\af e_4$ & $e_2 e_2 = e_3$  \\
$\af\ne$ 0&& $e_3 e_1 = e_4$            & $e_3 e_2 = i e_4$     & $e_3 e_3 = e_4$\\

\hline 
$\cd 4{49}(\af)$ &: & 
$e_1 e_1 = e_3 + \af e_4$ & $e_2 e_1 = -i\af e_4$ & $e_2 e_2 = e_3$\\
$\af\ne$ 0& & $e_3 e_1 = e_4$& $e_3 e_2 = i e_4$ & $e_3 e_3 = e_4$\\

\hline 
$\cd 4{50}(\af)$ &:    & 
$e_1 e_1 = e_3 + \af e_4$ & $e_2 e_1 = e_4$ & $e_2 e_2 = e_3$ & $e_3 e_3 = e_4$\\

\hline
$\cd 4{51}(\af)$ &: &
$e_1 e_1 = e_3 + \af e_4$ &   $e_2 e_2 = e_3$ & $e_3 e_1 = e_4$\\
&& $e_3 e_2 = i e_4$ & $e_3 e_3 = e_4$\\

\hline 
$\cd 4{52}$ &: &  
$e_1 e_1 = e_3 +  e_4$ &   $e_2 e_2 = e_3$ &   $e_3 e_3 = e_4$\\

\hline 
$  \cd 4{53}$ &: &  
$e_1 e_1 = e_3$ &  $e_1e_2=-\frac 12e_4$ & $e_1e_3=e_4$ & $e_2e_1=\frac 12e_4$\\
 && $e_2 e_2 = e_3$ & $e_2e_3=ie_4$ &  $e_3e_1=e_4$ & $e_3e_2=ie_4$ \\

\hline
 $\cd 4{54}(\af)$ &:&  
$e_1 e_1 = e_3$ & $e_1e_2=e_4$ & $e_1e_3=\af e_4$ & $e_2e_1=e_4$\\
&& $ e_2 e_2=e_3$ & $e_2e_3=-i(\af+1)e_4$ & $e_3e_1=(\af+1) e_4$ & $e_3e_2=-i\af e_4$\\

\hline
$\cd 4{55}(\af)$ &: &  
$e_1 e_1 = e_3$ & $e_1 e_2 = e_4$ & $e_1 e_3 = \af e_4$ \\
& & $e_2 e_1 = e_4$ & $e_2 e_2 = e_3$  & $e_3 e_1 = (\af+1) e_4$ \\
		
\hline
$\cd 4{56}$ &: &  
$e_1 e_1 = e_3$ & $e_1 e_2 = e_4$ & $e_2 e_1 = -e_4$ & $e_2 e_2 = e_3$ \\
& & $e_3 e_1 = e_4$ & $e_3 e_2 = i e_4$ & $e_3 e_3 = e_4$\\

\hline
$\cd 4{57}(\af,\bt)$ &: &  
$e_1 e_1 = e_3$ & $e_1 e_2 = \af e_4$ &   $e_2 e_1 = (\af+\bt)e_4$ & $e_2 e_2 = e_3$ \\
$\bt\not\in\{0,-2\af\}$& & $e_3 e_1 = e_4$          & $e_3 e_2 = i e_4$     & $e_3 e_3 = e_4$\\

\hline
$\cd 4{58}$ &: &  
$e_1 e_1 = e_3$ & $e_1 e_3 = e_4$  & $e_2 e_1 = e_4$ & $e_2 e_2 = e_3$\\
&& $e_2 e_3 = i e_4$  & $e_3 e_1 = e_4$ & $e_3 e_2 = i e_4$ \\

\hline 
$\cd 4{59}(\af,\bt)$ &: &   
$e_1 e_1 = e_3$ & $e_1 e_3 = \af e_4$ & $e_2 e_2 = e_3$\\
 $\bt\ne 0$ && $e_2 e_3 = \bt e_4$ & $e_3 e_1 = (\af+1) e_4$  & $e_3 e_2 = \bt e_4$ \\

\hline		
$\cd 4{60}$ &: & 
$e_1 e_1 = e_3$ & $e_1 e_3 = i e_4$  & $e_2 e_2 = e_3$\\
&& $e_2 e_3 = e_4$ & $e_3 e_1 = (i+1) e_4$ & $e_3 e_2 = (i+1) e_4$  \\

\hline
$\cd 4{61}(\af)$ &: & 
 $e_1 e_1 = e_3$ & $e_1 e_3 = -i\af e_4$ & $e_2 e_2 = e_3$\\
 && $e_2 e_3 = \af e_4$ & $e_3 e_1 = (1-i\af) e_4$ & $e_3 e_2 = (\af+i) e_4$  \\

\hline
$\cd 4{62}$ &:&  
$e_1 e_1 = e_3$  & $e_1e_3=e_4$ & $e_2 e_2=e_3$\\
&& $e_2e_3=-2ie_4$ & $e_3e_1=2e_4$ & $e_3e_2=-ie_4$\\

\hline
$\cd 4{63}$ &:&  
$e_1 e_1 = e_3$ & $e_1e_3=-\frac 12e_4$ & $ e_2 e_2=e_3$\\
&& $e_2e_3=-\frac i2e_4$ & $e_3e_1=\frac 12e_4$ & $e_3e_2=\frac i2e_4$\\

\hline 
$\cd 4{64}(\af)$ &: & 
$e_1 e_1 = e_3$ & $e_1 e_3 = \af e_4$ & $e_2 e_2 = e_3$  & $e_3 e_1 = (\af+1) e_4$ 		\\

\hline
$\cd 4{65}(\af)$ &: &  
$e_1 e_1 = e_3$ &   $e_1 e_3 = \af e_4$ & $e_2 e_2 = e_3$ \\
$\af\ne 0$&&  $e_3 e_1 = (\af+1) e_4$  & $e_3 e_2 = i e_4$ \\
\hline
$\cd 4{66}$ &: &  
$e_1 e_1 = e_3$& $e_2 e_1 = e_4$ & $e_2 e_2 = e_3$ \\
&& $e_2 e_3 = e_4$ & $e_3 e_2 = e_4$ \\

\hline 
$\cd 4{67}$ &: & 
$e_1 e_1 = e_3$ &  $e_2 e_2 = e_3$ &  $e_3 e_3 = e_4$\\

\hline
$\cd 4{68}$ &: & 
$e_1 e_1 = e_3+e_4$ & $e_1 e_3 = i e_4$  &$e_2 e_2 = e_3$ &\\
&&$e_2 e_3 = e_4$&$e_3 e_1 =   ie_4$ & $e_3 e_2 = e_4$  \\

\hline
$\cd 4{69}$ &: & 
$e_1 e_1 = e_3$ & $e_1 e_3 =  ie_4$  &$e_2 e_2 = e_3$ &\\
&&$e_2 e_3 = e_4$&$e_3 e_1 =   ie_4$ & $e_3 e_2 = e_4$  \\

\hline
$\cd 4{70}$ &: & 
$e_1 e_1 = e_3$ & $e_1 e_3 =    e_4$  &$e_2 e_2 = e_3$ & $e_3 e_1 = e_4$ &   \\

	\end{longtable}

All these algebras are non-isomorphic, except
\begin{longtable}{lll}
$\cd 4{43}(\af)\cong\cd 4{43}(-\af)$ &
$\cd 4{44}(\af,\bt,\gm)\cong\cd 4{44}(\af,-\bt,-\gm)$& 
$\cd 4{47}(\af,\bt)\cong \cd 4{47}(\af,-\bt)$\\
$\cd 4{50}(\af)=\cd 4{50}(-\af)$& 
$\cd 4{54}(\af)\cong\cd 4{54}(-\af-1)$&
$\cd 4{57}(\af,\bt)\cong \cd 4{57}(\af+\bt,-\bt)$\\
$\cd 4{59}(\af,\bt)\cong\cd 4{59}(\af,-\bt)$.&&
\end{longtable}

		   \subsection{$1$-dimensional central extensions of $\cd {3*}{03}$}
	
	Let us use the following notations 
	$$
	\begin{array}{rclrclrclrcl}
	\nb 1& = &\Dl 13-\Dl 31, &\nb 2& = &\Dl 23-\Dl 32, &\nb 3& = &\Dl 11, &\nb 4& = &\Dl 12,\\
	\nb 5& = &\Dl 13, &\nb 6& = &\Dl 22,        &\nb 7& = &\Dl 23,        &\nb 8& = &\Dl 33.
	\end{array}
	$$
	Take $\0=\sum_{i=1}^8\af_i\nb i\in {\rm H}^2_{\mathfrak{CD}}(\cd {3*}{03})$.
	If 
	$$
	\phi=
	\begin{pmatrix}
	x &    y  &  0\\
    z &    u  &  0\\
    v &    w  &  xu-yz
	\end{pmatrix}\in\aut{\cd {3*}{03}},
	$$
	then
	$$
	\phi^T\begin{pmatrix}
	\af_3        &  \af_4        & \af_1+\af_5\\
	0            &  \af_6        & \af_2+\af_7\\
	-\af_1       &  -\af_2       & \af_8
	\end{pmatrix} \phi=
	\begin{pmatrix}
	\af_3^*      &  \af_4^*-\af^*& \af_1^*+\af_5^*\\
	\af^*        &  \af_6^*      & \af_2^*+\af_7^*\\
	-\af_1^*     &  -\af_2^*     & \af_8^*
	\end{pmatrix},
	$$
	where
	\begin{longtable}{l}
	$\af^*_1 = (xu - yz)(\af_1x + \af_2z -\af_8v)$\\
    $\af^*_2 = (xu - yz)(\af_1y + \af_2u - \af_8w)$\\
    $\af^*_3 = \af_3x^2 + \af_5vx + \af_8v^2 + z(\af_4x + \af_6z + \af_7v)$\\
    $\af^*_4 = x(2\af_3y + \af_4u + \af_5w) + z(\af_4y + 2\af_6u + \af_7w) + v(\af_5y+ \af_7u + 2\af_8w)$\\
    $\af^*_5 = (xu - yz)(\af_5x + \af_7z + 2\af_8v)$\\
    $\af^*_6 = \af_3y^2 + \af_4uy + \af_6u^2 + w(\af_5y + \af_7u + \af_8w)$\\
    $\af^*_7 = (xu - yz)(\af_5y + \af_7u + 2\af_8w)$\\
    $\af^*_8 = \af_8(xu - yz)^2.$
	\end{longtable}
	Hence, $\phi\langle\0\rangle=\langle\0^*\rangle$, where $\0^*=\sum\limits_{i=1}^8 \af_i^*  \nb i$. We are only interested in $\0$ with  $(\af_1,\af_2,\af_5,\af_7,\af_8) \neq (0,0,0,0,0).$
	
	\begin{enumerate}
	    \item $\af_8\ne 0$. Choosing $w=\frac {\af_1y + \af_2u}{\af_8}$ and $v=\frac{\af_1x + \af_2z}{\af_8}$ we get $\af^*_1=\af^*_2=0$ without changing $\af_8^*$. Thus, we may also assume that $\af_1=\af_2=0$ from the very beginning. Taking $v=w=0$, we obtain $\af^*_1=\af^*_2=0$ and 
	    \begin{longtable}{ll}
$\af^*_3 = \af_3x^2 + \af_4xz + \af_6z^2$ & 
$\af^*_4 = 2\af_3xy + \af_4(ux + yz) + 2\af_6uz$\\
$\af^*_5 = (xu - yz)(\af_5x + \af_7z)$ &
$\af^*_6 = \af_3y^2 + \af_4uy + \af_6u^2$\\
$\af^*_7 = (xu - yz)(\af_5y + \af_7u)$ &
$\af^*_8 = \af_8(xu - yz)^2.$
	    \end{longtable}
	    \begin{enumerate}
	        \item $\af_3\ne 0$ and $\af_4^2-4\af_3\af_6\ne 0$. Then the equation $\af_3^*=0$ has two distinct solutions: $x_1=\mu_1z$ and $x_2=\mu_2z$, where $\mu_1,\mu_2\in\Co$. Since $\af_6^*=0$ has the same coefficients, its solutions are $y_1=\mu_1u$ and $y_2=\mu_2u$. Observe that $x_1u-y_2z=(\mu_1-\mu_2)uz$. So, we may take $x=\mu_1z$ and $y=\mu_2z$, where $u,z\ne 0$, and obtain $\af^*_3=\af^*_6=0$. Thus, we obtain the family of representatives 
	        \[ \la\af\nb 4+\bt\nb 5+\gm\nb 7+\nb 8\ra.\]
	        \begin{enumerate}
	            \item $\bt,\gm\ne 0$. Then choosing $v=w=x=u=0$, $y = -\gm$, $z = -\bt$ we obtain the family of representatives of distinct orbits $\la\af\nb 4+\nb 5+\nb 7+\nb 8\ra$.
	            \item $\bt\ne 0$ and $\gm=0$. Then we obtain the representatives $\la\nb 5+\nb 8\ra$ and $\la\nb 4+\nb 5+\nb 8\ra$ depending on whether $\af=0$ or not. Observe that $\la\nb 5+\nb 8\ra$ belongs to the family $\la\af\nb 4+\nb 5+\nb 7+\nb 8\ra$, while $\la\nb 4+\nb 5+\nb 8\ra$ does not.
	            \item $\bt=0$ and $\gm\ne 0$. Then we obtain the same representatives as in the previous item.
	            \item $\bt=\gm=0$. Then we obtain the representatives $\la\nb 8\ra$ and $\la\nb 4+\nb 8\ra$ depending on whether $\af=0$ or not. 
	        \end{enumerate} 

	        %{\red  \sout{$\la\nb 4+\nb 5+\nb 8\ra$,  $\la\nb 5+\nb 8\ra$,}} $\la\nb 4+\nb 8\ra$ and $\la\nb 8\ra$.
	        
	        \item $\af_3\ne 0$ and $\af_4^2-4\af_3\af_6=0$. Observe that $(\af^*_4)^2-4\af^*_3\af^*_6=(\af_4^2 - 4\af_3\af_6)(xu - yz)^2=0$. So, choosing $x=-\frac{\af_4z}{2\af_3}$ we obtain $\af_3^*=0$, and hence $\af_4^*=0$ as well. 
	        \begin{enumerate}
	        	\item $\af_7\ne 0$ and $\af_4\af_5-2\af_3\af_7\ne 0$. Then taking $u=-\frac{\af_5y}{\af_7}$, we obtain $\af^*_7=0$. This gives the representative $\la\nb 5+\nb 6+\nb 8\ra$.

	        	\item $\af_7\ne 0$ and $\af_4\af_5-2\af_3\af_7=0$. Then  $\af^*_5=0$ as well. This leads to the family of representatives of distinct orbits $\la\af\nb 6+\nb 7+\nb 8\ra_{\af\ne 0}$.
	        	    
	        	\item $\af_7=0$, $\af_5 \neq 0$ and $\af_4\neq 0$. 
	        	Then choosing 
	        	$y=0,$
	        	 $z=\sqrt{\frac{\af_3}{\af_8}},$
	        	$u=\frac{\af_5}{\af_8}$ we obtain the representative  $\la\nb 5+ \nb 6+\nb 8\ra$ found above.

	    	\item $\af_7=0$, $\af_5 \neq 0$ and $\af_4= 0$. 
	        	Then choosing  $u=1,$  $z=-\frac{\af_5}{\af_8},$  $y = 1$, we obtain the family  $\la  \af \nb 6+\nb 7   +\nb 8\ra_{\af\neq 0}$ found above.

	    	\item $\af_7=0$ and  $\af_5 = 0$. 
	        	Then choosing  $u=0,$ $z= \sqrt{\frac{\af_3}{\af_8}},$ $y = 1$, we obtain the representative  $\la  \nb 6 +\nb 8\ra.$
	
	        \end{enumerate}
	        
	        \item $\af_3=0.$ 
	   
	        \begin{enumerate}
	            \item $\af_6\ne 0$. Then choosing $x=0$ and $y=z=1$ we get $\af^*_3\ne 0$, so we are in one of the previous cases.
	            
	            \item $\af_6=0$ and $\af_4\ne 0$. Then choosing $x=z=1$ and $y\ne u$, we get $\af^*_3\ne 0$, so we are in one of the previous cases.
	            
	            \item $\af_6=\af_4=0$. Then we obtain the family of representatives $\la\bt\nb 5+\gm\nb 7+\nb 8\ra$ which is a subfamily of a family considered above.
	        \end{enumerate}

\end{enumerate}
	
	    	%\begin{align*}
	        %    \af^*_1 &= (xu - yz)(\af_1x + \af_2z),\\
            %    \af^*_2 &= (xu - yz)(\af_1y + \af_2u),\\
            %    \af^*_3 &= \af_3x^2 + \af_5vx + z(\af_4x + \af_6z + \af_7v),\\
                %\af^*_4 &= x(2\af_3y + \af_4u + \af_5w) + z(\af_4y + 2\af_6u + \af_7w) + %v(\af_5y+ \af_7u),\\
   %             \af^*_5 &= (xu - yz)(\af_5x + \af_7z),\\
  %              \af^*_6 &= \af_3y^2 + \af_4uy + \af_6u^2 + w(\af_5y + \af_7u),\\
 %               \af^*_7 &= (xu - yz)(\af_5y + \af_7u).
%	    \end{align*}

	    \item $\af_8=0$ and $\af_7\ne 0$. Then 
		choosing $u=-\frac{\af_5y}{\af_7}$ we obtain $\af^*_8=0$ and $\af_7^*=0$, so we may assume $\af_7=\af_8=0$ since the beginning. Then taking $y=0$ we get $\af^*_7=\af^*_8=0$ and
		
	\begin{longtable}{ll}
$\af^*_1 = (\af_1x + \af_2z)ux$ &
$\af^*_2 = \af_2u^2x$\\

$\af^*_3 = \af_3x^2 + \af_4xz + \af_5vx + \af_6z^2$ &
$\af^*_4 = \af_4ux + \af_5wx + 2\af_6uz$\\

$\af^*_5 = \af_5ux^2$ &
$\af^*_6 = \af_6u^2.$
	\end{longtable}
	\begin{enumerate}
	    \item $\af_5\ne 0$. Then choosing $x\ne 0$ and $v=-\frac{\af_3x^2 + \af_4xz + \af_6z^2}{\af_5 x}$, $w=-\frac{\af_4ux+ 2\af_6uz}{\af_5 x}$ we obtain $\af^*_3=\af^*_4=0$.
	    \begin{enumerate}
	        \item $\af_2\ne 0$. Then choosing $z=-\frac{\af_1x}{\af_2}$ and $u=\frac{\af_5 x}{\af_2}$ we have two representatives $\la\nb 2+\nb 5\ra$ and $\la\nb 2+\nb 5+\nb 6\ra$ depending on whether $\af_6=0$ or not.
	        \item $\af_2=0$. Then we have two families of representatives of distinct orbits $\la\af\nb 1+\nb 5\ra $ and $\la\af\nb 1+\nb 5+\nb 6\ra $ depending on whether $\af_6=0$ or not.
	    \end{enumerate}
	    \item $\af_5=0$. Then $\af^*_5=0$ and
	    \begin{longtable}{lll}
$\af^*_1 = (\af_1x + \af_2z)ux$ & 
$\af^*_2 = \af_2u^2x$ &
$\af^*_3 = \af_3x^2 + \af_4xz + \af_6z^2$\\
$\af^*_4 = (\af_4x + 2\af_6z)u$ &
$\af^*_6 = \af_6u^2.$
	\end{longtable}
	    \begin{enumerate}
	        \item $\af_6\ne 0$. Then choosing $z=-\frac{\af_4x}{2\af_6}$ we have $\af^*_4=0$, $\af^*_2,\af^*_6$ are as above and 
	        \begin{longtable}{ll}
$\af^*_1 = \left(\af_1-\frac{\af_2\af_4}{2\af_6}\right)ux^2$& 
$\af^*_3 = \left(\af_3-\frac{\af_4^2}{4\af_6}\right)x^2.$
	        \end{longtable}
	       \begin{enumerate}
	           \item If $\af_2\ne 0$, then choosing $x=\frac{\af_6}{\af_2}$, we obtain two families of representatives 
	           \[\la\nb 2+\af\nb 3+\nb 6\ra \mbox{ and }\la\nb 1+\nb 2+\af\nb 3+\nb 6\ra\] depending on whether $\af_1-\frac{\af_2\af_4}{2\af_6}=0$ or not. 
	           
	           $\bullet$ The first family gives two distinct orbits, whose representatives are 
	           \[\la\nb 2+\nb 6\ra \mbox{ and }\la\nb 2+\nb 3+\nb 6\ra.\] 
	           $\bullet$ The second family splits into three distinct orbits, whose representatives are 
\[\la\nb 1+\nb 3+\nb 6\ra_{\af\not\in\{-1,0\}}, \ 
	            \la\nb 1+ \nb 6\ra_{\af=0} \mbox{ and } \la\nb 1+\nb 4\ra_{\af=-1} .\] 
 Note that $\orb\la\nb 1+\nb 3+\nb 6\ra=\orb\la\nb 2+\nb 3+\nb 6\ra$.
	        
	         \item If $\af_2=0$ and $\af_1\neq 1,$ then choosing $u=\frac{\af_1x^2}{\af_6}$, we obtain the family  $\la\nb 1+\af \nb 3+\nb 6\ra$, which gives two distinct orbits with representatives $\la\nb 1+\nb 6\ra$ and $\la\nb 1+\nb 3+\nb 6\ra$. Both of them were found above.

              \item If $\af_2=0$ and $\af_1=0,$ then we obtain two distinct orbits with representatives $\la \nb 6\ra$ and $\la \nb 3+\nb 6\ra$.  They define algebras with $2$-dimensional annihilators 
	       	       \end{enumerate}
	       
	        \item $\af_6=0$. Then $\af^*_6=0$. 
	        
	        \begin{enumerate}
	            \item If $\af_4\ne 0$, then choosing $z=-\frac{\af_3x}{\af_4}$, we get $\af^*_3=0$ and
	        \begin{longtable}{lll}
$\af^*_1 = \left(\af_1-\frac{\af_2\af_3}{\af_4}\right)ux^2$&
$\af^*_2 = \af_2u^2x$&
$\af^*_4 = \af_4ux.$
	        \end{longtable}
	        Hence, we have two families of representatives $\la\af\nb 1+\nb 4\ra$ and $\la\af\nb 1+\nb 2+\nb 4\ra$ depending on whether $\af_2=0$ or not. The first family gives two distinct orbits with representatives $\la\nb 4\ra$ and $\la\nb 1+\nb 4\ra$. Moreover, the representative $\la\nb 4\ra$ determines an algebra with $2$-dimensional annihilator which has been found before. The second family also gives two distinct orbits whose representatives are $\la\nb 2+\nb 4\ra$ and $\la\nb 1+\nb 2+\nb 4\ra$.   Observe that 
	        \[\orb\la\nb 2+\nb 4\ra=\orb\la\nb 1+\nb 4\ra\mbox{ and }\orb\la\nb 1+\nb 2+\nb 4\ra=\orb\la\nb 1+\nb 3+\nb 6\ra.\]
	        
	        \item If $\af_4=0$, then $\af^*_4=0$ and
	         \begin{longtable}{lll}
$\af^*_1 = (\af_1x + \af_2z)ux$ &
$\af^*_2 = \af_2u^2x$ &
$\af^*_3 = \af_3x^2.$
	        \end{longtable}
	        $\bullet$ $\af_3\ne 0$. If $\af_2\ne 0$, then choosing $z=-\frac{\af_1}{\af_3}$, $x=\frac{\af_2}{\af_3}$, $y=0$ and $u=1$, we obtain the representative $\la\nb 2+\nb 3\ra$,   but it has the same orbit as $\orb\la\nb 1+\nb 6\ra$ found above. Otherwise, we have two representatives $\la\nb 3\ra$ and $\la\nb 1+\nb 3\ra$ depending on whether $\af_1=0$ or not.  The representative $\la\nb 1+\nb 3\ra$ will substitute $\la\nb 2+\nb 6\ra$ found above, since they belong to the same orbit.  Observe also that $\la\nb 3\ra$ gives an algebra with $2$-dimensional annihilator which has been found before.
	        
	        $\bullet$ If $\af_3= 0$, then we have only one representative  $\la \nb 2\ra,$ which gives a Lie algebra.
	        
	    \end{enumerate}
	    
	        \end{enumerate}
	\end{enumerate}

	\end{enumerate}

Summarizing, we have the following representatives of different orbits giving non-Lie $\cd{}{}$-algebras:
\begin{longtable}{lll}
    $\la\nb 1+\nb 3\ra$ & 
    $\la\nb 1+\nb 3+\nb 6\ra$ & 
    $\la\nb 1+\nb 4\ra$\\
    $\la\af\nb 1+\nb 5\ra$ & 
    $\la\af\nb 1+\nb 5+\nb 6\ra$ & 
    $\la\nb 1+\nb 6\ra$\\
    $\la\nb 2+\nb 5\ra$ & 
    $\la\nb 2+\nb 5+\nb 6\ra$ & 
    $\la\af\nb 4+\nb 5+\nb 7+\nb 8\ra$\\
    $\la\nb 4+\nb 5+\nb 8\ra$ & 
    $\la\nb 4+\nb 8\ra$ &  
    $\la\nb 5+\nb 6+\nb 8\ra$\\
    $\la\af\nb 6+\nb 7+\nb 8\ra_{ \af\ne 0}$ & 
    $\la\nb 6+\nb 8\ra$ & $\la\nb 8\ra$.
\end{longtable}
and those giving Lie algebras: $\la \nb 2\ra.$
	\begin{longtable}{lllllllllllllllllll}	
		$\cd 4{71}$ &: & $e_1 e_1 = e_4$ & $e_1 e_2 = e_3$  & $e_1 e_3 = e_4$ & $e_2 e_1 = -e_3$  & $e_3 e_1 = -e_4$\\
		\hline
		$\cd 4{72}$ &: & $e_1 e_1 = e_4$ & $e_1 e_2 = e_3$ & $e_1 e_3 = e_4$\\
		&& $e_2 e_1 = -e_3$  & $e_2 e_2 = e_4$ & $e_3 e_1 = -e_4$\\
		\hline
		$\cd 4{73}$ &: & $e_1 e_2 = e_3 + e_4$ & $e_1 e_3 = e_4$ & $e_2 e_1 = -e_3$  & $e_3 e_1 = -e_4$\\
		\hline
		$\cd 4{74}(\af)$ &: & $e_1 e_2 = e_3$ & $e_1 e_3 = (\af+1)e_4$ & $e_2 e_1 = -e_3$ & $e_3 e_1 = -\af e_4$\\
		\hline
		$\cd 4{75}(\af)$ &: & $e_1 e_2 = e_3$ & $e_1 e_3 = (\af+1)e_4$ & $e_2 e_1 = -e_3$ & $e_2 e_2 = e_4$ & $e_3 e_1 = -\af e_4$\\        \hline
		$\cd 4{76}$ &: & $e_1 e_2 = e_3$ & $e_1 e_3 = e_4$ & $e_2 e_1 = -e_3$ & $e_2 e_2 = e_4$ & $e_3 e_1 = -e_4$\\        
		\hline
		$\cd 4{77}$ &: & $e_1 e_2 = e_3$ & $e_1 e_3 = e_4$ & $e_2 e_1 = -e_3$ & $e_2 e_3 = e_4$ & $e_3 e_2 = -e_4$\\
		\hline
		$\cd 4{78}$ &: & $e_1 e_2 = e_3$ & $e_1 e_3 = e_4$ & $e_2 e_1 = -e_3$\\
		&& $e_2 e_2 = e_4$ & $e_2 e_3 = e_4$ & $e_3 e_2 = -e_4$\\
		\hline
		$\cd 4{79}(\af)$ &: & $e_1 e_2 = e_3+\af e_4$ & $e_1 e_3 = e_4$ & $e_2 e_1 = -e_3$ & $e_2 e_3 = e_4$ & $e_3 e_3 = e_4$\\
		\hline
		$\cd 4{80}$ &: & $e_1 e_2 = e_3+e_4$ & $e_1 e_3 = e_4$ & $e_2 e_1 = -e_3$ & $e_3 e_3 = e_4$\\
		\hline
		$\cd 4{81}$ &: & $e_1 e_2 = e_3+e_4$ & $e_2 e_1 = -e_3$ & $e_3 e_3 = e_4$\\
		\hline
		$\cd 4{82}$ &: & $e_1 e_2 = e_3$ & $e_1 e_3 = e_4$ & $e_2 e_1 = -e_3$ & $e_2 e_2 = e_4$ & $e_3 e_3 = e_4$\\
		\hline
		$\cd 4{83}( \af\ne 0)$ &: & $e_1 e_2 = e_3$ & $e_2 e_1 = -e_3$ & $e_2 e_2 = \af e_4$ & $e_2 e_3 = e_4$ & $e_3 e_3 = e_4$\\
		\hline
		$\cd 4{84}$ &: & $e_1 e_2 = e_3$ & $e_2 e_1 = -e_3$ & $e_2 e_2 = e_4$ & $e_3 e_3 = e_4$\\
		\hline
		$\cd 4{85}$ &: & $e_1 e_2 = e_3$ & $e_2 e_1 = -e_3$ & $e_3 e_3 = e_4$\\
		\hline
		$\cd 4{86}$ &: & $e_1 e_2 = e_3$ & $e_2 e_1 = -e_3$ & $e_2 e_3 = e_4$  &  $e_3 e_2 = -e_4$
	\end{longtable}

		  \subsection{$1$-dimensional central extensions of $\cd {3*}{04}$}
	
	Let us use the following notations 
	$$
	\begin{array}{rclrclrclrcl}
	\nb 1& = &\Dl 11, &\nb 2& = &\Dl 12, &\nb 3& = &\Dl 13, &\nb 4& = &\Dl 21,\\
	\nb 5& = &\Dl 23, &\nb 6& = &\Dl 31, &\nb 7& = &\Dl 32, &\nb 8& = &\Dl 33.
	\end{array}
	$$
	Take $\0=\sum_{i=1}^8\af_i\nb i\in {\rm H}^2_{\mathfrak{CD}}(\cd {3*}{04})$.
	If 
	$$
	\phi=
	\begin{pmatrix}
	x      &      y    &  0\\
    -\lb y &    x-y    &  0\\
    z      &      u    &  x^2-xy+\lb y^2
	\end{pmatrix}\in\aut{\cd {3*}{04}},
	$$
	then
	$$
	\phi^T\begin{pmatrix}
	\af_1            &  \af_2        & \af_3\\
	\af_4            &      0        & \af_5\\
	\af_6            &  \af_7        & \af_8
	\end{pmatrix} \phi=
	\begin{pmatrix}
	\af_1^*+\lb\af^*      &  \af_2^*  & \af_3^*\\
	\af_4^*+\af^*         &  \af^*    & \af_5^*\\
	\af_6^*               &  \af_7^*  & \af_8^*
	\end{pmatrix},
	$$
	where
	\begin{longtable}{lcl}
	$\af^*_1$&$ =$&$ \af_1x^2 + \lb(-\af_1 + \af_2 + \af_4) y^2 + \af_8z^2 -\lb\af_8 u^2 -2\lb(\af_2 + \af_4) xy + (\af_3 + \af_6)xz -$ \\
	&&$ \lb(\af_5 + \af_7) yz-\lb(\af_5 + \af_7) ux + \lb(- \af_3 + \af_5 - \af_6 + \af_7) uy$\\
    $\af^*_2$&$ =$&$ \af_2x^2 - \lb\af_4 y^2 + (\af_1 - \af_2)xy + \af_7xz + (\af_6 - \af_7)yz -\af_5\lb uy  + \af_3ux +\af_8uz$\\
    $\af^*_3$&$ =$&$ (x^2 - xy + \lb y^2)(\af_3x - \lb\af_5 y + \af_8z)$\\
    $\af^*_4$&$ =$&$ \af_4x^2 + (- \af_1 + (1-\lb)\af_2 + \af_4)y^2 - \af_8u^2 + (\af_1 - \af_2 - 2\af_4)xy + (\af_3 - \af_5)yz +$\\
    &&$ \af_5xz+ (- \af_5 + \af_6 - \af_7)ux  + (- \af_3 + \af_5 - \af_6 +  (1-\lb)\af_7)uy + \af_8uz$\\
    $\af^*_5$&$ =$&$ (x^2 - xy + \lb y^2)(\af_5x + (\af_3 - \af_5)y + \af_8u)$\\
    $\af^*_6$&$ =$&$ (x^2 - xy + \lb y^2)(\af_6x -\lb\af_7 y + \af_8z)$\\
    $\af^*_7$&$ =$&$ (x^2 - xy + \lb y^2)(\af_7x + (\af_6-\af_7)y + \af_8u)$\\
    $\af^*_8$&$ =$&$ (x^2 - xy + \lb y^2)^2\af_8.$
	\end{longtable}
Hence, $\phi\langle\0\rangle=\langle\0^*\rangle$, where $\0^*=\sum\limits_{i=1}^8 \af_i^*  \nb i$.
All orbits with $\af_8=0,$ which give non-Leibniz algebras are classified in \cite{kkp19geo}.
It is easy to see that only orbit 
$\la \af_1 \nabla_1+\af_2 \nabla_2+\af_4 \nabla_4+\nabla_5 \ra_{\lb=0}$
gives non-split Leibniz algebra.
It is follows that there are only three non-split extensions of $\cd {3*}{04}(0)$ 
which give Leibniz algebras: 
\begin{longtable}{lll}
$\la \nabla_1+\nabla_5 \ra_{\lb=0}$& 
$\la \nabla_5 \ra_{\lb=0}$ &
$\la \nabla_1+\nabla_2+\nabla_5 \ra_{\lb=0}$ 
\end{longtable}
Thanks to \cite{kkp19geo}, the first and the second representatives we can joint 
with families of algebras $\D{4}{38}(\lb \neq 0)$ and  $\D{4}{40}(\lb \neq 0)$ from  \cite{kkp19geo}.
The last representative gives one new algebra:
\begin{longtable}{llllllllllll}
$\D{4}{00} $&$:$&  
$e_1e_1=e_4$ & $e_1e_2=e_4$& $e_2 e_1=e_3$ & $e_2 e_2 = e_3$ & $e_2e_3=e_4$ 

\end{longtable}

From here we are only interested in $\0$ with $ \af_8\ne 0.$ %$(\af_3,\af_5,\af_6,\af_7) \ne (0,0,0,0)$ (if $\lb \ne 0$) and $(\af_3,\af_6,\af_7) \ne (0,0,0)$ (if $\lb = 0$).
		Choosing $u=-\frac{\af_7x + (\af_6-\af_7)y}{\af_8}$ and $z=-\frac{\af_6x -\lb\af_7 y}{\af_8}$, we may make $\af^*_6=\af^*_7=0$, so we shall assume that $\af_6=\af_7=0$. Taking $u=z=0$, we get $\af^*_6=\af^*_7=0$, $\af^*_8$ is as above and
	    \begin{align*}
	        \af^*_1 &= \af_1x^2 -2\lb(\af_2 + \af_4) xy + \lb(-\af_1 + \af_2 + \af_4) y^2, \\
            \af^*_2 &= \af_2x^2 + (\af_1 - \af_2)xy - \lb\af_4 y^2,\\
            \af^*_3 &= (x^2 - xy + \lb y^2)(\af_3x - \lb\af_5 y),\\
            \af^*_4 &= \af_4x^2 + (\af_1 - \af_2 - 2\af_4)xy + (- \af_1 + (1-\lb)\af_2 + \af_4)y^2,\\
            \af^*_5 &= (x^2 - xy + \lb y^2)(\af_5x + (\af_3 - \af_5)y).
	    \end{align*}
	Observe that the equations $\af^*_2=0$ and $\af^*_4=0$ have the same discriminant  $(\af_1 - \af_2)^2+4\lb\af_2\af_4$. So, we have the following subcases.
	\begin{enumerate}
	    \item $\af_1\ne 0$, $\af_2\ne 0$ and $(\af_1 - \af_2)^2+4\lb\af_2\af_4\ne 0$. Note that \[(\af^*_1 - \af^*_2)^2+4\lb\af^*_2\af^*_4=((\af_1 - \af_2)^2+4\lb\af_2\af_4)(x^2 - xy + \lb y^2)^2,\] so the condition $(\af_1 - \af_2)^2+4\lb\af_2\af_4\ne 0$ is invariant under the automorphisms   preserving the condition $\af_6=\af_7=0$.  The equation $\af^*_2=0$ has two distinct solutions $x_1=\mu_1 y$ and $x_2=\mu_2 y$, where $\mu_1,\mu_2\in\Co$. Suppose that $\mu_1^2 - \mu_1 + \lb=\mu_2^2 - \mu_2 + \lb=0$. Then $\mu_1^2-\mu_1=\mu_2^2-\mu_2$, i.e. $(\mu_1-\mu_2)(\mu_1+\mu_2-1)=0$. Since $\mu_1\ne\mu_2$, we conclude that $\mu_1+\mu_2=1$. But $\mu_1+\mu_2=-\frac{\af_1 - \af_2}{\af_2}$, whence $\af_1-\af_2=-\af_2$, which contradicts the assumption that $\af_1\ne 0$. Thus, we may choose $y\ne 0$ and $x=\mu_iy$ for which $\af^*_2=0$ and  $\af_8^*=1$. Then $(\af^*_1 - \af^*_2)^2+4\lb\af^*_2\af^*_4=(\af^*_1)^2$, whence $\af^*_1\ne 0$. This gives the following family of representatives $\la\nb 1+\af\nb 3+\bt\nb 4+\gm\nb 5+\nb 8\ra$. 
	    
	    \begin{enumerate}
	    	\item $\lb\not\in\{0,\frac 14\}$. If $(\af,\bt,\gm)\ne(\af',\bt',\gm')$, such that $\bt,\bt'\ne\frac{1\pm\sqrt{1-4\lb}}{2\lb}$, then 
	    	\[\orb\la\nb 1+\af\nb 3+\bt\nb 4+\gm\nb 5+\nb 8\ra=\orb\la\nb 1+\af'\nb 3+\bt'\nb 4+\gm'\nb 5+\nb 8\ra\] if and only if either $(\af',\bt',\gm')=\left(\pm(\gm-\af\bt)\sqrt{\frac{-\lb}{1-\bt+\lb\bt^2}},\frac 1\lb-\bt,\pm(\frac\gm\lb-\frac\af\lb-\bt\gm)\sqrt{\frac{-\lb}{1-\bt+\lb\bt^2}}\right)$, where $\sqrt{\frac{-\lb}{1-\bt+\lb\bt^2}}$ is one of the two values of the square root, or $(\af',\bt',\gm')=(-\af,\bt,-\gm)$. Moreover,   for $(\af,\gm)\ne(\af',\gm')$ one has 
	    	\[\orb\La\nb 1+\af\nb 3+\frac{1\pm\sqrt{1-4\lb}}{2\lb}\nb 4+\gm\nb 5+\nb 8\Ra=\orb\La\nb 1+\af'\nb 3+\frac{1\pm\sqrt{1-4\lb}}{2\lb}\nb 4+\gm'\nb 5+\nb 8\Ra\] 
	    	if and only if $(\af',\gm')=(-\af,-\gm)$.

	    	\item $\lb=\frac 14$. If $(\af,\bt,\gm)\ne(\af',\bt',\gm')$, such that $\bt,\bt'\ne 2$, then \[\orb\la\nb 1+\af\nb 3+\bt\nb 4+\gm\nb 5+\nb 8\ra=\orb\la\nb 1+\af'\nb 3+\bt'\nb 4+\gm'\nb 5+\nb 8\ra\] 
	    	if and only if either $(\af',\bt',\gm')=\left(\pm\frac{i(\gm-\af\bt)}{\bt-2},4-\bt,\pm\frac{i(4\gm-4\af-\bt\gm)}{\bt-2}\right)$, or $(\af',\bt',\gm')=(-\af,\bt,-\gm)$. Moreover,  for $(\af,\gm)\ne(\af',\gm')$ one has 
	    	\[\orb\la\nb 1+\af\nb 3+2\nb 4+\gm\nb 5+\nb 8\ra=\orb\la\nb 1+\af'\nb 3+2\nb 4+\gm'\nb 5+\nb 8\ra\]
	    	if and only if $(\af',\gm')=(-\af,-\gm)$.

	    	\item $\lb=0$. If $(\af,\bt,\gm)\ne(\af',\bt',\gm')$, then 
	    	\[\orb\la\nb 1+\af\nb 3+\bt\nb 4+\gm\nb 5+\nb 8\ra=\orb\la\nb 1+\af'\nb 3+\bt'\nb 4+\gm'\nb 5+\nb 8\ra\] 
	    	if and only if  $(\af',\bt',\gm')=(-\af,\bt,-\gm)$.

	    \end{enumerate}

	    \item $\af_1\ne 0$, $\af_2\ne 0$ and $(\af_1 - \af_2)^2+4\lb\af_2\af_4=0.$ Observe that 
	    \[(\af^*_1 - \af^*_2)^2+4\lb\af^*_2\af^*_4=((\af_1 - \af_2)^2+4\lb\af_2\af_4)(x^2 - xy + \lb y^2)^2=0,\]  so $(\af_1 - \af_2)^2+4\lb\af_2\af_4=0$ is invariant under the automorphisms preserving the condition $\af_6=\af_7=0$. 
	    \begin{enumerate}
	        \item $\af_4\ne 0$ and $\af_1 - \af_2\ne 0$. It follows that $\lb=-\frac{(\af_1 - \af_2)^2}{4\af_2\af_4}\ne 0$. Observe that if we take $x=-\frac{(\af_1-\af_2)y}{2\af_2}$, then $x^2-xy+\lb y^2=-\frac{(\af_1\af_2 - \af_2^2 - \af_1\af_4 - \af_2\af_4)(\af_1 - \af_2)y^2}{4\af_2^2\af_4}$. So we have the following two subcases.
	        
	        \begin{enumerate}
	        	\item $\af_1\af_2 - \af_2^2 - \af_1\af_4 - \af_2\af_4\ne 0$. Then choosing $x=-\frac{(\af_1-\af_2)y}{2\af_2}$ and an appropriate $y\ne 0$ we have $\af^*_2=0$ and $\af^*_8=1$. But since $(\af^*_1 - \af^*_2)^2+4\lb\af^*_2\af^*_4=0$, we automatically get $\af^*_1=0$. Thus, we obtain the family of representatives $\la\af\nb 3+\bt\nb 4+\gm\nb 5+\nb 8\ra$. 
	        	
	        	$\bullet$ If $\af\ne 0$, then taking $z=u=y=0$ and $x=\af$ we obtain the family $\la \nb 3+\af\nb 4+\bt\nb 5+\nb 8\ra$. 
	        	Since $\lb\ne 0$, then $\la \nb 3+\af\nb 4+\bt\nb 5+\nb 8\ra$ splits into the subfamily $\la\nb 3+\af\nb 4+\bt\nb 5+\nb 8\ra_{\af\ne 0}$ of representatives of distinct orbits and the separate orbits whose representatives are $\la\nb 5+\nb 8\ra$ and $\La\frac{1\pm\sqrt{1-4\lb}}2\nb 3+\nb 5+\nb 8\Ra$.

	        %	$\star$ Let $\lb=0$. {\red Then $\la \nb 3+\af\nb 4+\bt\nb 5+\nb 8\ra$ splits into two separate subfamilies $\la\nb 3+\af\nb 4+\nb 8\ra$ and $\la\nb 3+\af\nb 4+\nb 5+\nb 8\ra$ of representatives of distinct orbits.}
	        	
	        	$\bullet$ If $\af=0$ and $\bt\ne 0$, then taking $z=u=y=0$ and $x=\sqrt{\bt}$ we obtain the family $\la \nb 4+\af\nb 5+\nb 8\ra$.   Observe that $\orb\la\nb 4+\af\nb 5+\nb 8\ra=\orb\la\nb 4+\af'\nb 5+\nb 8\ra$ if and only if $\af=\pm\af'$.
	        	
	        	$\bullet$ If $\af=\bt=0$ and $\gm\ne 0$, then taking $z=u=y=0$ and $x=\gm$ we obtain the representative $\la \nb 5+\nb 8\ra$   found above. 
	        	
	        	$\bullet$ If $\af=\bt=\gm=0$, then we obtain the representative $\la\nb 8\ra$.
	        	
	        	\item $\af_1\af_2 - \af_2^2 - \af_1\af_4 - \af_2\af_4=0$. Observe that $\af_1\ne-\af_2$, since otherwise $\af_1=\af_2=0$. Then $\af_4=\frac{\af_1\af_2 - \af_2^2}{\af_1+\af_2}$, and substituting this into $(\af_1 - \af_2)^2+4\lb\af_2\af_4$ we obtain $\frac{((4\lb-1)\af_2^2 + \af_1^2)(\af_1 - \af_2)}{\af_1 + \af_2}$. Hence we must have $(4\lb-1)\af_2^2 + \af_1^2=0$. Observe that $\lb\ne\frac 14$, since otherwise $\af_1=0$. Then $\af_1=\pm\sqrt{1-4\lb}\af_2$ and $\af_4=\frac{\pm\sqrt{1-4\lb} - 1}{\pm\sqrt{1-4\lb} + 1}\af_2 =-\frac 1\lb\left(\frac{1 \mp\sqrt{1-4\lb}}2\right)^2\af_2$. Hence, we obtain the following family of representatives 
	        	\[\La \pm\sqrt{1-4\lb}\nb 1 + \nb 2 + \af\nb 3-\frac 1\lb\left(\frac{1 \mp \sqrt{1-4\lb}}2\right)^2\nb 4+\bt\nb 5+\nb 8\Ra.\] 
	        	
 	        	$\bullet$ If $\bt\not\in\left\{\frac{1+\sqrt{1-4\lb}}{2\lb}\af,\frac{1-\sqrt{1-4\lb}}{2\lb}\af\right\}$, then choosing $x=-\frac{2(\af - \bt)}{(\pm\sqrt{1-4\lb}-1)\bt+2\af}$ and $y=\frac{2\bt}{(\pm\sqrt{1-4\lb}-1)\bt+2\af}$ we obtain 2 families representatives \[\La \pm\sqrt{1-4\lb}\nb 1 + \nb 2+\af\nb 3-\frac 1\lb\left(\frac{1 \mp \sqrt{1-4\lb}}2\right)^2\nb 4+\nb 8\Ra_{\af\ne 0}.\] Observe that $\orb\La \pm\sqrt{1-4\lb}\nb 1 + \nb 2+\af\nb 3-\frac 1\lb\left(\frac{1 \mp \sqrt{1-4\lb}}2\right)^2\nb 4+\nb 8\Ra=\orb\La \pm\sqrt{1-4\lb}\nb 1 + \nb 2+\af'\nb 3-\frac 1\lb\left(\frac{1 \mp \sqrt{1-4\lb}}2\right)^2\nb 4+\nb 8\Ra$ if and only if $\af=\pm\af'$.

	        	$\bullet$ If $\bt\in\left\{\frac{1+\sqrt{1-4\lb}}{2\lb}\af,\frac{1-\sqrt{1-4\lb}}{2\lb}\af\right\}$, then we obtain 6 representatives depending on whether $\af=0$ or not. These are $\La \sqrt{1-4\lb}\nb 1 + \nb 2 + \nb 3-\frac 1\lb\left(\frac{1 - \sqrt{1-4\lb}}2\right)^2\nb 4+\frac{1\pm\sqrt{1-4\lb}}{2\lb}\nb 5+\nb 8\Ra$, $\La -\sqrt{1-4\lb}\nb 1 + \nb 2 + \nb 3-\frac 1\lb\left(\frac{1 + \sqrt{1-4\lb}}2\right)^2\nb 4+\frac{1\pm\sqrt{1-4\lb}}{2\lb}\nb 5+\nb 8\Ra$ and $\La \pm\sqrt{1-4\lb}\nb 1 + \nb 2-\frac 1\lb\left(\frac{1 \mp \sqrt{1-4\lb}}2\right)^2\nb 4+\nb 8\Ra$. The last two representatives will be joined with the family above. 
	        \end{enumerate}

	        \item $\af_4\ne 0$ and $\af_1 - \af_2=0$. Then $(\af_1 - \af_2)^2+4\lb\af_2\af_4=0$ implies $\lb=0$. We have
\begin{longtable}{ll}
$\af^*_1 = \af_2x^2$& 
$\af^*_2 = \af_2x^2$\\
$\af^*_3 = \af_3(x - y)x^2$&
$\af^*_4 = \af_4(x-y)^2$\\
$\af^*_5 = (\af_5x + (\af_3 - \af_5)y)(x - y)x$&
$\af^*_8 = \af_8x^2(x-y)^2.$
	        \end{longtable}
	        \begin{enumerate}
	        	\item $\af_3\ne 0$ and $\af_3 - \af_5\ne 0$. Then choosing $y=-\frac{\af_5x}{\af_3 - \af_5}$, we obtain $x^2-xy=\frac{\af_3x^2}{\af_3 - \af_5}\ne 0$ and $\af^*_5=0$. 
 If we take
$x=\frac{\sqrt{\alpha_2} (\alpha_3-\alpha_5)}{\alpha_3 \sqrt{\alpha_8}}$, then we obtain the family of representatives $\la\nb 1+\nb 2+\af\nb 3+\bt\nb 4 + \nb 8\ra_{\af,\bt\ne 0}$. Note that 
\[\orb\la\nb 1+\nb 2+\af\nb 3+\bt\nb 4 + \nb 8\ra=\orb\la\nb 1+\nb 2+\af'\nb 3+\bt'\nb 4 + \nb 8\ra\] 
if and only if $(\af',\bt')=(\pm\af,\bt)$.
	        	
	        	\item $\af_3\ne 0$ and $\af_3 - \af_5=0$. Then $\af^*_3=\af^*_5$, so taking 
	     $x=\sqrt{\frac{\alpha_4}{\alpha_8}}$ and 
	     $y=\sqrt{\frac{\alpha_4}{\alpha_8}}-\sqrt{\frac{\alpha_2}{\alpha_8}}$,  we obtain the family $\la\nb 1+\nb 2+\af\nb 3+\nb 4 + \af\nb 5 + \nb 8\ra_{\af\ne 0}$, where 
	     \[\orb\la\nb 1+\nb 2+\af\nb 3+\nb 4 + \af\nb 5 + \nb 8\ra=\orb\la\nb 1+\nb 2+\af'\nb 3+\nb 4 + \af'\nb 5 + \nb 8\ra\] 
	     if and only if $\af=\pm\af'$.
	        	
	        	\item $\af_3=0$. Then $\af^*_3=0$, so we obtain the family of representatives \[\la\nb 1+\nb 2+\af\nb 4 + \bt\nb 5 + \nb 8\ra_{\af\ne 0}.\] 
	        	
	        	$\bullet$ If $\bt\ne 0$, then taking $z=u=0$, $x=\bt$ and $y=\bt-1$ we obtain the family $\la\nb 1+\nb 2+\af\nb 4 + \nb 5 + \nb 8\ra_{\af\ne 0}$ of representatives of distinct orbits. 
	        	
	        	$\bullet$ If $\bt=0$, then we obtain the representative $\la\nb 1+\nb 2+\nb 4  + \nb 8\ra$ which will be joined with the family $\la\nb 1+\nb 2+\af\nb 3+\nb 4 + \af\nb 5 + \nb 8\ra_{\af\ne 0}$.
	        \end{enumerate}

	        \item $\af_4=0$. Then $(\af_1 - \af_2)^2+4\lb\af_2\af_4=0$ implies $\af_1 - \af_2=0$. Hence
\begin{longtable}{ll}
	        $\af^*_1 = \af_2x(x - 2\lb y)$& 
            $\af^*_2 = \af_2x^2$\\
            $\af^*_3 = (x^2 - xy + \lb y^2)(\af_3x - \lb\af_5 y)$&
            $\af^*_4 = -\lb\af_2y^2$\\
            $\af^*_5 = (x^2 - xy + \lb y^2)(\af_5x + (\af_3 - \af_5)y)$&
            $\af^*_8 = (x^2 - xy + \lb y^2)^2\af_8.$
	        \end{longtable}
	        
	        We have the following subcases. 
	        
\begin{enumerate}
    \item 	        If  $\lb\ne 0$, then choosing   $x=0$, we obtain a representative of the family $\la\af\nb 3+\bt\nb 4+\gm\nb 5+\nb 8\ra$ considered above. 
    
    %$x=2\lb y$ we get $\af^*_1=0$ and
%	         \begin{longtable}{lll}
%$\af^*_2 = 4\lb^2 y^2 \af_2$&
%$\af^*_3 = (2\af_3 - \af_5)(4\lb - 1)\lb^2y^3$& 
%$\af^*_4 = -\lb y^2 \af_2$\\
%$\af^*_5 = (\af_3 + (2\lb-1)\af_5)(4\lb - 1)\lb y^3$ &
%$\af^*_8 = (4\lb - 1)^2\lb^2y^4\af_8.$
%	        \end{longtable}
%	        Applying one more automorphism with $y=\frac 1{4\lb-1}\sqrt{\frac{\af_2}{\lb\af_8}}$, we obtain the family $\la 4\lb\nb 2 +\af\nb 3 - \nb 4 + \bt\nb 5 + \nb 8\ra$. Observe that 
%{\gr 	    pode apagar  $4\lambda,$ e receber a familia que está em cima   \[\orb\la 4\lb\nb 2 +\af\nb 3 - \nb 4 + \bt\nb 5 + \nb 8\ra=\orb\la 4\lb\nb 2 +\af'\nb 3 - \nb 4 + \bt'\nb 5 + \nb 8\ra\] } if and only if $(\af,\bt)=(-\af',-\bt')$.
	        
%\item	        If $\lb=\frac 14$, then
%\begin{longtable}{ll}
%$\af^*_1 = \af_2x\left(x - \frac y2\right)$&
%$\af^*_2 = \af_2x^2$\\
%$\af^*_3 = \left(\af_3x - \frac {\af_5}4y\right)\left(x - \frac y2\right)^2$&
%$\af^*_4 = -\frac{\af_2}4y^2$\\
%$\af^*_5 = (\af_5x + (\af_3 - \af_5)y)\left(x - \frac y2\right)^2$&
%	        \end{longtable}
%	        Choosing $x=0$ and $y=\sqrt{-\frac{4\af_2}{\af_8}}$, we obtain the family $\la\af\nb 3+\nb 4+\bt\nb 5+\nb 8\ra$. 
%	        
%	        $\bullet$ However, if $\af\ne 0$, this family coincides with $\la\nb 3+\af\nb 4+\bt\nb 5+\nb 8\ra_{\af\ne 0}.$ 
%	        
%	        $\bullet$ If $\af=0$ it coincides with $\la\nb 4+\bt\nb 5+\nb 8\ra$, both families being found above.
	        
\item	  If $\lb=0$, then $\af^*_4=0$ and
	        \begin{longtable}{lcl}
	        $\af^*_1 = \af_2x^2$&
            $\af^*_2 = \af_2x^2$& 
            $\af^*_3 = \af_3 x^2(x - y)$\\
\multicolumn{2}{l}{            $\af^*_5 = (\af_5x + (\af_3 - \af_5)y)x(x-y)$}&
            $\af^*_8 = \af_8 x^2(x - y)^2.$
	        \end{longtable}
	        So, we have the family $\la\nb 1+\nb 2+\af\nb 3+\bt\nb 5+\nb 8\ra$. 
	        
	        $\bullet$ If $\af\ne 0$ and $\af\ne\bt$, then taking $z=u=0$, $x=\frac{\bt}{\af}-1$ and $y=\frac{\bt}{\af}$ we obtain the family $\la\nb 1+\nb 2+\af\nb 3+\nb 8\ra_{\af\ne 0}$ which will be joined with the family $\la\nb 1+\nb 2+\af\nb 3+\bt\nb 4 + \nb 8\ra_{\af,\bt\ne 0}$ found above.  
	        
	        $\bullet$ If  $\af=\bt$, then we obtain the family 
	        $\la\nb 1+\nb 2+\af\nb 3+\af\nb 5+\nb 8\ra$,  where \[\orb\la\nb 1+\nb 2+\af\nb 3+\af\nb 5+\nb 8\ra=\orb\la\nb 1+\nb 2+\af'\nb 3+\af'\nb 5+\nb 8\ra\] 
	        if and only if $\af=\pm\af'$.
	        
	        $\bullet$ If $\af=0$ and $\bt\ne 0$, then choosing $x=\bt$ and $y=\bt+1$ we obtain the representative $\la\nb 1+\nb 2+\nb 5 + \nb 8\ra$ which will be joined with the family $\la\nb 1+\nb 2+\af\nb 4 + \nb 5 + \nb 8\ra_{\af\ne 0}$ found above.
	    
\end{enumerate}   \end{enumerate}
	    \item $\af_1\ne 0$ and $\af_2=0$. Then
\begin{longtable}{lr}
$\af^*_1 = \af_1x^2 -2\lb\af_4 xy + \lb(-\af_1 + \af_4) y^2$&
$\af^*_2 = y(\af_1x - \lb\af_4 y).$
	    \end{longtable}
	    Clearly, there are $x,y$ such that $\af^*_1\ne 0$, $\af^*_2\ne 0$ and $x^2-xy+\lb y^2\ne 0$, so we are in the case considered above.
	    
	    \item $\af_1=0$ and $\af_2\ne 0$. Then
\begin{longtable}{lr}
$\af^*_1 = \lb(\af_2 + \af_4)y(y - 2x)$ &
$\af^*_2 = \af_2x^2 -\af_2xy - \lb\af_4 y^2$\\
\multicolumn{2}{c}{$\af^*_4 = \af_4x^2 - (\af_2 + 2\af_4)xy + ((1-\lb)\af_2 + \af_4)y^2.$}
	    \end{longtable}
	    
	    \begin{enumerate}
	    	\item $\af_2+\af_4\ne 0$ and $\lb\ne 0$. Then we can find $x,y$ such that $\af^*_1\ne 0$, $\af^*_2\ne 0$ and $x^2-xy+\lb y^2\ne 0$, and this case has been considered above.
	    	
	    	\item  $\af_2+\af_4\ne 0$ and $\lb=0$. Then $\af^*_1=0$ and
	    	\begin{longtable}{ll}
	    	$\af^*_2 = \af_2x(x - y)$ &
	    	$\af^*_4 =   (\af_4x - (\af_2+\af_4)y)(x - y)$\\
	    	$\af^*_8 = \af_8(x - y)^2x^2.$
	    	\end{longtable}
	    %	The discriminant of $\af^*_4=0$ is $\af_2^2\ne 0$. Hence, there are distinct $\mu_1,\mu_2\in\Co$ such that $y_1=\mu_1 x$ and $y_1=\mu_2 x$ are the solutions of $\af^*_4=0$. 
	    	
	%    	$\bullet$ If $\{\mu_1,\mu_2\}\ne\{0,1\}$, then we can find $y_i=\mu_i x$ such that $\af^*_4=0$ and $(x - y)^2x^2\ne 0$. 
	
	    	Then  choosing $x=\sqrt{\frac{\af_2 + \af_4}{\af_8}}$ and $y=\frac{\af_4}{\af_2 + \af_4}\sqrt{\frac{\af_2 + \af_4}{\af_8}}$  we have the family $\la\nb 2+\af\nb 3+\bt\nb 5+\nb 8\ra$. 
	    	Moreover,  for $(\af,\bt)\ne(\af',\bt')$ one has     \[\orb\la\nb 2+\af\nb 3+\bt\nb 5+\nb 8\ra=\orb\la\nb 2+\af'\nb 3+\bt'\nb 5+\nb 8\ra\] 
	    	if and only if $(\af,\bt)=(-\af',-\bt')$. 
	    	
	  %  	$\bullet$ If $\{\mu_1,\mu_2\}=\{0,1\}$, then {\gr ???? $\frac{\af_2 + 2\af_4}{\af_2 + \af_4}\in\{0,1\}$.
	   % 	---
	    %	} This implies $\af_4\in\{0,-\frac 12\af_2\}$. If $\af_4=0$, then we are in the previous case. {\gr apagar - If $\af_4=-\frac 12\af_2$, then
	   % 	\begin{align*}
	    %	\af^*_2 &= \af_2x(x - y),\\
	    %	\af^*_4 &= -\frac{\af_2}2(x - y)(x + y),\\
	    %	\af^*_8 &= \af_8(x - y)^2x^2.
	    %	\end{align*}
	    %	Choosing $y=-x$ and $x=\sqrt{\frac{\af_2}{2\af_8}}$ we obtain the family $\la\nb 2+\af\nb 3+\bt\nb 5+\nb 8\ra$ found above.}
	    	
	    	\item $\af_2+\af_4=0$. 
	    	%Then $(\af_1 - \af_2)^2+4\lb\af_2\af_4=\af_2^2(1-4\lb)$, so $\lb=\frac 14$. 
	    	In this case $\af^*_1=0$ and 
	    	\begin{longtable}{ll}
$\af^*_2 = \af_2(x^2 -xy + \lb y^2)$ &
$\af^*_4 = -\af_2(x^2 - xy +\lb y^2).$
	    	\end{longtable}
	    	Hence, we have the family of representatives $\la\nb 2+\af\nb 3-\nb 4+\bt\nb 5+\nb 8\ra$. 
	    	
	    	$\bullet$ If   $\lb\ne 0$ and $\af\ne \frac{1\pm\sqrt{1-4\lb}}{2}\bt$, then choosing $x=  \bt \sqrt{\frac \lb{\af^2 -\af\bt+\lb\bt^2}} $ and $y=  \frac\af\lb \sqrt{\frac \lb{\af^2 -\af\bt+\lb\bt^2}}$ we obtain the family $\la\nb 2-\nb 4+\af\nb 5+\nb 8\ra_{\af\ne 0}.$ 
	    		Moreover, one has
	    		\[\orb\la\nb 2-\nb 4+\af\nb 5+\nb 8\ra=\orb\la\nb 2-\nb 4+\af'\nb 5+\nb 8\ra\] 
	    	if and only if $\af=\pm\af'$.
	    	
	     	$\bullet$ If  $\lb\not\in\{0,\frac 14\}$, $\af=\frac{1\pm\sqrt{1-4\lb}}{2}\bt$ and $\bt\ne 0$, then choosing $x=\frac{\af\bt^2 - \bt^3 + \af}{\bt(2\af - \bt)}$ and $y=\frac{1-\bt^2}{2\af - \bt}$ we obtain two representatives $\La\nb 2+\frac{1\pm\sqrt{1-4\lb}}2\nb 3-\nb 4+\nb 5+\nb 8\Ra$.
	    	
	    	$\bullet$ If $\lb\not\in\{0,\frac 14\}$ and $\af=\bt=0$, then we obtain the representative $\la\nb 2-\nb 4+\nb 8\ra$ which will be joined with the family $\la\nb 2-\nb 4+\af\nb 5+\nb 8\ra_{\af\ne 0}$.
	    	
	    	$\bullet$ If $\lb=\frac 14$ and $\af=\frac 12\bt$, then we obtain the family of representatives $\la\nb 2+\af\nb 3-\nb 4+2\af\nb 5+\nb 8\ra$. Moreover, one has
	    		\[\orb\la\nb 2+\af\nb 3-\nb 4+2\af\nb 5+\nb 8\ra=\orb\la\nb 2+\af'\nb 3-\nb 4+2\af'\nb 5+\nb 8\ra\] 
	    	if and only if $\af=\pm\af'$.
	    	
	    	$\bullet$ If $\lb=0$ and $\af\ne 0$, then choosing $x=\frac 1\af$ and $y=\frac 1\af-\af$ we obtain the family $\la\nb 2+\nb 3-\nb 4+\af\nb 5+\nb 8\ra$ of representatives of distinct orbits.
	    	
	    	$\bullet$ If $\lb=0$ and $\af=0$, then we obtain two representatives $\la\nb 2-\nb 4+\nb 8\ra$ and $\la\nb 2-\nb 4+\nb 5+\nb 8\ra$ depending on whether $\bt=0$ or not.
 
	    \end{enumerate}

	    \item $\af_1=\af_2=0$. Then
	    \begin{longtable}{lll}
	        $\af^*_1 = \lb\af_4y(y - 2x)$ &
            $\af^*_2 =  -\lb\af_4 y^2 $&
            $\af^*_4 = \af_4(x-y)^2.$
	    \end{longtable}
	    
	    \begin{enumerate}
	         \item If $\af_4\ne 0$ and $\lb\ne 0$, then we can find $x,y$ such that $\af^*_1\ne 0$, $\af^*_2\ne 0$ and $\af^*_8\ne 0$, so we are in one of the previous situations. 
	    	\item $\af_4=0$ and $\af_3^2-\af_3\af_5+\lb\af_5^2\ne 0$. Then choosing $x= \af_3 - \af_5$ and $y=-\af_5$, we obtain the representative $\la\nb 3+\nb 8\ra$. 
	    	
	    	$\bullet$ If $\lb=0$, then it  will be joined with the family $\la\nb 3+\af\nb 4+\nb 8\ra_{\af\ne 0}$ found  below   . 
	    	
	    	$\bullet$ If $\lb\ne 0$, then   $\la\nb 3+\nb 8\ra$ has the same orbit as $\la\nb 5+\nb 8\ra$ found above.
	    	
	    	\item $\af_4=0$ and $\af_3^2-\af_3\af_5+\lb\af_5^2=0$. 
	    	
	    	$\bullet$ If $\lb\ne 0$, then we have the representatives $\La\frac{1\pm\sqrt{1-4\lb}}2\nb 3+\nb 5+\nb 8\Ra$ and $\la\nb 8\ra$ found above. 
	    	
	    	%$\bullet$ If $\lb=\frac 14$, then we have 2 representatives $\la\nb 3+2\nb 5+\nb 8\ra$ and $\la\nb 8\ra$ found above. 
	    	
	    	$\bullet$ If $\lb=0$, then either $\af_3=0$, or $\af_3=\af_5$. So we have 3 representatives $\la\nb 3+\nb 5+\nb 8\ra$, $\la\nb 5+\nb 8\ra$ and $\la\nb 8\ra$.
	    	
	    	\item $\af_4\ne 0$ and $\lb=0$. 
	    	
	    	$\bullet$ If   $\af_3\ne 0$ and $\af_3\ne\af_5$, then choosing $x=  \af_3 - \af_5$ and $y=  -\af_5 $, we obtain the family $\la\nb 3+\af\nb 4+\nb 8\ra_{ \af\ne 0}$ of representatives of distinct orbits. 
	    	%\begin{enumerate}
	    %	\item If $\af\ne 0$, then $\orb\la\af\nb 3+\nb 4+\nb 8\ra=\orb\la\nb 3+\frac 1{\af^2}\nb 4+\nb 8\ra$, so this family has been found above. 
	    	
	    %	    \item And if $\af=0$, then we obtain the representative $\la\nb 4+\nb 8\ra$ belonging to the family $\la\nb 4+\af\nb 5+\nb 8\ra$ found above. 
	    	
	    %	\end{enumerate}
	    	$\bullet$ If $\af_3=\af_5\ne 0$, then we have the representative $\la\nb 3+\nb 4+\nb 5+\nb 8\ra$. 
	    	
	    	$\bullet$ If $\af_3=0$, then we obtain the  family $\la\nb 4+\af\nb 5+\nb 8\ra$, where $\orb\la\nb 4+\af\nb 5+\nb 8\ra=\orb\la\nb 4+\af'\nb 5+\nb 8\ra$ if and only if $\af=\pm\af'$.
	    \end{enumerate}
	\end{enumerate}

All the found orbits are pairwise distinct, except those whose equality is indicated explicitly in the table below.

{\tiny
\begin{longtable}{|ll|}
				
				\hline
\multicolumn{2}{|c|}{	 $\lb\not\in\{0,\frac 14\}$}\\
				\hline		

$\La \pm\sqrt{1-4\lb}\nb 1 + \nb 2 + \nb 3  -\frac 1\lb\left(\frac{1 \mp \sqrt{1-4\lb}}2\right)^2\nb 4+\frac{1+\sqrt{1-4\lb}}{2\lb}\nb 5+\nb 8\Ra$& 
$\la\nb 3+\af\nb 4+\bt\nb 5+\nb 8\ra_{\af\ne 0}$\\

$\La \pm\sqrt{1-4\lb}\nb 1 + \nb 2 + \nb 3 -\frac 1\lb\left(\frac{1 \mp \sqrt{1-4\lb}}2\right)^2\nb 4+\frac{1-\sqrt{1-4\lb}}{2\lb}\nb 5+\nb 8\Ra$& 
$\La\frac{1\pm\sqrt{1-4\lb}}2\nb 3+\nb 5+\nb 8\Ra$\\

$\La \pm\sqrt{1-4\lb}\nb 1 + \nb 2+\af\nb 3-\frac 1\lb\left(\frac{1 \mp \sqrt{1-4\lb}}2\right)^2\nb 4+\nb 8\Ra^{ O(\af)=O(-\af)}$&
$\la\nb 4+\af\nb 5+\nb 8\ra^{O(\af)=O(-\af)}$\\

				%\La\frac{1\pm\sqrt{1-4\lb}}2\nb 1+\nb 2-\frac{1\pm\sqrt{1-4\lb}}2\nb 3-\nb 5+\frac{1\pm\sqrt{1-4\lb}}2\nb 6+\nb 7\Ra,\\
				%\La\frac{1\mp\sqrt{1-4\lb}}2\nb 1+\nb 2-\frac{1\pm\sqrt{1-4\lb}}2\nb 3-\nb 5+\frac{1\pm\sqrt{1-4\lb}}2\nb 6+\nb 7\Ra,\\
				%\la\nb 1+\af\nb 2-\nb 3+\nb 6\ra,\\
				%\La \pm\sqrt{1-4\lb}\nb 1 + \nb 2 +\frac{\sqrt{1-4\lb} \mp 1}{\sqrt{1-4\lb} \pm 1}\nb 4+\nb 8\Ra,\\  	   	   
				%\La 2\lb\nb 1+\frac{-1\pm\sqrt{1-4\lb}}2\nb 3+\nb 4+\frac{\sqrt{1-4\lb}\mp 1}{\sqrt{1-4\lb}\pm 1}\nb 5+\frac{1\pm\sqrt{1-4\lb}}2\nb 6+\nb 7 \Ra,\\
				%\La \frac{1\pm\sqrt{1-4\lb}}2\nb 1+\frac{-1\pm\sqrt{1-4\lb}}2\nb 3+\nb 4+\frac{\sqrt{1-4\lb}\mp 1}{\sqrt{1-4\lb}\pm 1}\nb 5+\frac{1\pm\sqrt{1-4\lb}}2\nb 6+\nb 7 \Ra,\\
      
				%\La \nb 1 +  \frac{1\pm\sqrt{1-4\lb}}{2}\nb 3 + \nb 5 \Ra,\\
				%\la\nb 1+\af\nb 3+\bt\nb 5+\nb 6\ra_{\af,\bt\ne 0},\\
				%\La\nb 1+\frac{1\pm\sqrt{1-4\lb}}2\af\nb 3+\af\nb 5+\frac{1\pm\sqrt{1-4\lb}}2\nb 6+\nb 7\Ra_{\af\ne 0},\\
				%\La\nb 1+\frac{1\mp\sqrt{1-4\lb}}2\af\nb 3+\af\nb 5+\frac{1\pm\sqrt{1-4\lb}}2\nb 6+\nb 7\Ra_{\af^2+\frac{1\pm\sqrt{1-4\lb}}{2\lb}\af+1\ne 0},\\
				%\La \nb 1+\af\nb 3+ \frac{1\pm\sqrt{1-4\lb}}{2}\nb 6+\nb 7\Ra_{\af\ne 0},\\

				%{\red \sout{\la 4\lb\nb 2 +\af\nb 3 - \nb 4 + \bt\nb 5 + \nb 8\ra, O(\af,\bt)=O(-\af,-\bt),}}\\

$\la\nb 2-\nb 4+\af\nb 5+\nb 8\ra^{O(\af)=O(-\af)}$ & $\la\nb 5+\nb 8\ra$  \\

$\La\nb 2+\frac{1\pm\sqrt{1-4\lb}}2\nb 3-\nb 4+\nb 5+\nb 8\Ra$ &
$\la\nb 8\ra$\\

				%\la \nb 2+\af\nb 3+\nb 6\ra,\\
				%\La\nb 2+\af\nb 4+\frac{1\pm\sqrt{1-4\lb}}{2\lb}\nb 5+\nb 6\Ra_{\af\ne 0},\\
				%\la\nb 2+\af\nb 5+\nb 6\ra,\\

				%\la \nb 3 \ra,\\
				%\la \nb 3 + \nb 4 \ra,\\
				%\La \frac{1\pm\sqrt{1-4\lb}}{2}\nb 3 + \nb 5 \Ra,\\
				%\la\af\nb 3+\bt\nb 5+\nb 6\ra_{\af,\bt\ne 0},\\
				%\La\frac{1\pm\sqrt{1-4\lb}}2\af\nb 3+\af\nb 5+\frac{1\pm\sqrt{1-4\lb}}2\nb 6+\nb 7\Ra_{\af\ne 0},\\
				%\La\frac{1\mp\sqrt{1-4\lb}}2\af\nb 3+\af\nb 5+\frac{1\pm\sqrt{1-4\lb}}2\nb 6+\nb 7\Ra,\\
				%\la \af\nb 3+\nb 6\ra_{\af\ne 0},\\
				%\La \af\nb 3+ \frac{1\pm\sqrt{1-4\lb}}{2}\nb 6+\nb 7\Ra_{\af\ne 0},\\ 

				%\La\nb 4+\frac{1\pm\sqrt{1-4\lb}}{2\lb}\nb 5+\nb 6\Ra,\\

				%\la\af\nb 5+\nb 6\ra,\\

\multicolumn{2}{|l|}{
$\la\nb 1+\af\nb 3+\bt\nb 4+\gm\nb 5+\nb 8\ra, \mbox{ where }$}\\
\multicolumn{2}{|l|}{
${O(\af,\bt,\gm)=O(-\af,\bt,-\gm)},
O(\af,\bt,\gm)=O\left((\gm-\af\bt)\sqrt{\frac{-\lb}{1-\bt+\lb\bt^2}},\frac 1\lb-\bt,(\frac\gm\lb-\frac\af\lb-\bt\gm)\sqrt{\frac{-\lb}{1-\bt+\lb\bt^2}}\right) \mbox{ if }\bt\ne\frac{1\pm\sqrt{1-4\lb}}{2\lb}$}\\	 
				\hline

						\hline
\multicolumn{2}{|c|}{	 $\lb =0$}\\
				\hline		
				
$\la\nb 1+\nb 2+\af\nb 3+\nb 4 + \af\nb 5 + \nb 8\ra^{O(\af)=O(-\af)}$ &
$\la\nb 2-\nb 4+\nb 8\ra$\\

$\la\nb 1+\nb 2+\af\nb 3+\bt\nb 4 + \nb 8\ra_{\af\ne 0}^{ O(\af,\bt)=O(-\af,\bt)}$& 
$ \la\nb 3+\nb 4+\nb 5+\nb 8\ra $
			\\ 
				%\la\nb 1+\nb 2-\nb 3-\nb 5+\nb 6+\nb 7\ra,\\

$\la\nb 1+\nb 2+\af\nb 3+\af\nb 5+\nb 8\ra^{O(\af)=O(-\af)}$&
$\la\nb 3+\af\nb 4+\nb 8\ra$\\
				%\la\nb 1+\af\nb 2-\nb 3+\nb 6\ra,\\

$\la\nb 1+\nb 2+\af\nb 4 + \nb 5 + \nb 8\ra$ &
$  \la\nb 3+\nb 5+\nb 8\ra $
\\
				%\la\nb 1+\nb 2-\nb 5+\nb 7\ra,\\

$\la\nb 1+\af\nb 3+\bt\nb 4+\gm\nb 5+\nb 8\ra^{O(\af,\bt,\gm)=O(-\af,\bt,-\gm)}$& 
$\la\nb 4+\af\nb 5+\nb 8\ra^{ O(\af)=O(-\af)}$\\
				%\la\nb 1-\nb 3+\af\nb 4+\nb 6+\nb 7\ra,\\
				%\la \nb 1 + \nb 3 + \nb 5 \ra,\\
				%\la\nb 1+\af\nb 3+\bt\nb 5+\nb 6\ra_{\af,\bt\ne 0},\\
				%\la\nb 1+\af\nb 3+\af\nb 5+\nb 6+\nb 7\ra_{\af\ne 0},\\
				%\la\nb 1+\af\nb 3+\af\nb 5+\nb 7\ra_{\af^2+\af+1\ne 0},\\
				%\la\nb 1+\af\nb 3+\nb 7\ra_{\af\ne 0},\\
				
				%\la\nb 1+\nb 4+\nb 6+\nb 7\ra,\\
				%\la\nb 1+\af\nb 5+\nb 6+\nb 7\ra,\\
				%\la\nb 1+\af\nb 5+\nb 7\ra_{\af\ne 0},\\
				
$  \la\nb 2+\nb 3-\nb 4+\af\nb 5+\nb 8\ra $& 
$	\la\nb 5+\nb 8\ra$\\
				%\la\nb 2-\nb 3-\nb 5+\nb 6+\nb 7\ra,\\

$\la\nb 2+\af\nb 3+\bt\nb 5+\nb 8\ra^{ O(\af,\bt)=O(-\af,-\bt)}$& 
$\la\nb 8\ra$\\
				%\la \nb 2+\af\nb 3+\nb 6\ra,\\
				%\la\nb 2+\af\nb 4+\nb 5+\nb 6\ra_{\af\ne 0},\\

$  \la\nb 2-\nb 4+\nb 5+\nb 8\ra $&\\

				%\la\nb 2+\af\nb 5+\nb 6\ra,\\
				%\la\nb 2-\nb 5+\nb 7\ra,\\
				
				%\la \nb 3 \ra,\\
				%\la \nb 3 + \nb 4 \ra,\\
				%\la\af\nb 3+\nb 4+\nb 6+\nb 7\ra,\\
				%\la \nb 3 + \nb 5 \ra,\\
				%\la\af\nb 3+\bt\nb 5+\nb 6\ra_{\af,\bt\ne 0},\\
				%\la\af\nb 3+\af\nb 5+\nb 6+\nb 7\ra_{\af\ne 0},\\
				%\la\af\nb 3+\af\nb 5+\nb 7\ra,\\
				%\la \af\nb 3+\nb 6\ra_{\af\ne 0},\\
				%\la\af\nb 3+\nb 6+\nb 7\ra_{\af\ne 0},\\
				%\la\af\nb 3+\nb 7\ra_{\af\ne 0},\\
				
				%\la\nb 4+\nb 5+\nb 6\ra,\\
				
				%\la\af\nb 5+\nb 6\ra,\\
				%\la\af\nb 5+\nb 6+\nb 7\ra,\\
				%\la\af\nb 5+\nb 7\ra_{\af\ne 0},\\
 
						\hline
\multicolumn{2}{|c|}{	 $\lb =\frac{1}{4}$}\\
				\hline

\multicolumn{2}{|l|}{				
$\la\nb 1+\af\nb 3+\bt\nb 4+\gm\nb 5+\nb 8\ra, \mbox{ where }
O(\af,\bt,\gm)=O\left(\pm\frac{i(\gm-\af\bt)}{\bt-2},4-\bt,\pm\frac{i(4\gm-4\af-\bt\gm)}{\bt-2}\right),
				\mbox{ if }\bt\ne 2, O(\af,\bt,\gm)=O(-\af,\bt,-\gm)$}\\

$\la\nb 2+\af\nb 3-\nb 4+2\af\nb 5+\nb 8\ra_{  \af\ne 0}^{O(\af)=O(-\af)}$& 
$\la\nb 4+\af\nb 5+\nb 8\ra^{ O(\af)=O(-\af)}$\\

$\la\nb 2-\nb 4+\af\nb 5+\nb 8\ra^{   O(\af)=O(-\af)}$ & 
$	\la\nb 5+\nb 8\ra$\\

$\la\nb 3+\af\nb 4+\bt\nb 5+\nb 8\ra_{\af\ne 0}$ & $\la\nb 8\ra$\\

$ \la\frac 12\nb 3+\nb 5+\nb 8\ra $&\\
\hline				
		\end{longtable}
				}

	%%\newpage
	
	%\begin{landscape}

Denote $\Theta=\frac{1+\sqrt{1-4\lambda}}{2}.$

{\tiny
\begin{longtable}{ll}
$\La (2 \Theta-1)\nb 1 + \nb 2 + \nb 3- (1-\Theta)^2\lb^{-1}\nb 4+ \Theta\lb^{-1}\nb 5+\nb 8\Ra_{\lb\not\in\{0,\frac 14\}}$&
$\la\nb 2+\nb 3-\nb 4+\af\nb 5+\nb 8\ra_{\lb=0,\af\ne 1}$

\\
$\La (1-2 \Theta)\nb 1 + \nb 2 + \nb 3- \Theta^2\lb^{-1} \nb 4+ \Theta\lb^{-1} \nb 5+\nb 8\Ra_{\lb\not\in\{0,\frac 14\}}$&
$\la\nb 2+\af\nb 3-\nb 4+2\af\nb 5+\nb 8\ra_{\lb=\frac 14, \af\not\in\{0,\frac 12\}}^{ O(\af)=O(-\af)}$

\\	    	
$\La (2 \Theta-1)\nb 1 + \nb 2 + \nb 3 -  (1-\Theta)^2\lb^{-1} \nb 4+ (1-\Theta)\lb^{-1}\nb 5+\nb 8\Ra_{\lb\not\in\{0,\frac 14\}}$& 
   $\la\nb 2+\af\nb 3+\bt\nb 5+\nb 8\ra_{\lb=0}^{O(\af,\bt)=O(-\af,-\bt)}$

\\
$\La (1-2 \Theta) \nb 1 + \nb 2 + \nb 3-  \Theta^2\lb^{-1}\nb 4+ (1-\Theta)\lb^{-1} \nb 5+\nb 8\Ra_{\lb\not\in\{0,\frac 14\}}$& 
$\la\nb 2-\nb 4+\af\nb 5+\nb 8\ra_{\lb\ne 0}$

\\
$\La (2 \Theta-1)\nb 1 + \nb 2+\af\nb 3- (1-\Theta)^2\lb^{-1}\nb 4+\nb 8\Ra_{\lb\not\in\{0,\frac 14\}}^{O(\af)=O(-\af)}$&
$\la\nb 2-\nb 4+\nb 8\ra_{\lb=0}$

\\
$\La (1-2 \Theta) \nb 1 + \nb 2+\af\nb 3- \Theta^2\lb^{-1} \nb 4+\nb 8\Ra_{\lb\not\in\{0,\frac 14\}}^{O(\af)=O(-\af)}$& $\la\nb 3+\nb 4+\nb 5+\nb 8\ra_{\lb=0}$

\\
$\la\nb 1+\nb 2+\af\nb 3+\nb 4 + \af\nb 5 + \nb 8\ra_{\lb=0}^{O(\af)=O(-\af)}$&
$\la\nb 3+\af\nb 4+\bt\nb 5+\nb 8\ra_{\lb\ne 0,\af\ne 0}$

\\
$\la\nb 1+\nb 2+\af\nb 3+\bt\nb 4 + \nb 8\ra_{\lb=0,\af\ne 0}^{O(\af,\bt)=O(-\af,\bt)}$& 
$\la\nb 3+\af\nb 4+\nb 8\ra_{\lb=0}$

\\
$\la\nb 1+\nb 2+\af\nb 3+\af\nb 5+\nb 8\ra_{\lb=0}^{O(\af)=O(-\af)}$&
$\La \Theta \nb 3+\nb 5+\nb 8\Ra $

\\
$\la\nb 1+\nb 2+\af\nb 4 + \nb 5 + \nb 8\ra_{\lb=0}$&
$\La (1-\Theta) \nb 3+\nb 5+\nb 8\Ra_{\lb\not\in\{0,\frac 14\}}$
 	  
\\
$\La\nb 2+\Theta \nb 3-\nb 4+\nb 5+\nb 8\Ra $&
$\la\nb 4+\af\nb 5+\nb 8\ra^{ O(\af)=O(-\af)}$

\\
$\La\nb 2+(1-\Theta) \nb 3-\nb 4+\nb 5+\nb 8\Ra_{\lb\ne\frac 14}$& 
$\la\nb 5+\nb 8\ra $ \quad\quad\quad\quad$\la\nb 8\ra $\\

\multicolumn{2}{l}{$\la\nb 1+\af\nb 3+\bt\nb 4+\gm\nb 5+\nb 8\ra,  \mbox{ where } \ O(\af,\bt,\gm)=O(-\af,\bt,-\gm) \mbox{ and } \mbox{ if }\lb\ne 0\mbox{ and }\bt\ne\frac{1\pm\sqrt{1-4\lb}}{2\lb}$}\\
\multicolumn{2}{r}{$O(\af,\bt,\gm)=O\left(\pm(\gm-\af\bt)\sqrt{\frac{-\lb}{1-\bt+\lb\bt^2}},\frac 1\lb-\bt,\pm(\frac\gm\lb-\frac\af\lb-\bt\gm)\sqrt{\frac{-\lb}{1-\bt+\lb\bt^2}}\right)$
}\\
	    	
	            \end{longtable}}

{\tiny	 
\begin{longtable}{llllllll}

%\hline
$\cd {4}{87}(\lambda)$&$:$& 
$e_1 e_1 = \lambda e_3+(2 \Theta-1)e_4$& $e_1 e_2=e_4$ & $e_1e_3=e_4$& 
$e_2 e_1=e_3-(1- \Theta)^2 \lb^{-1}e_4$\\
$\lb \neq 0, \frac{1}{4}$&& $e_2 e_2=e_3$& $e_2e_3=\Theta\lb^{-1}e_4$ 
&$e_3e_3=e_4$\\

\hline$\cd {4}{88}(\lambda)$&$:$& 
$e_1 e_1 = \lambda e_3+(1-2 \Theta)e_4$& $e_1 e_2=e_4$ & $e_1e_3=e_4$& 
$e_2 e_1=e_3- \Theta^2 \lb^{-1}e_4$\\
$\lb \neq 0, \frac{1}{4}$&& $e_2 e_2=e_3$& $e_2e_3=\Theta\lb^{-1}e_4$ 
&$e_3e_3=e_4$\\

\hline$\cd {4}{89}(\lambda)$&$:$& 
$e_1 e_1 = \lambda e_3+(2 \Theta-1)e_4$& $e_1 e_2=e_4$ & $e_1e_3=e_4$& 
$e_2 e_1=e_3-(1- \Theta)^2 \lb^{-1}e_4$\\
$\lb \neq 0, \frac{1}{4}$&& $e_2 e_2=e_3$& $e_2e_3=(1-\Theta)\lb^{-1}e_4$ 
&$e_3e_3=e_4$\\

\hline$\cd {4}{90}(\lambda)$&$:$& 
$e_1 e_1 = \lambda e_3+(1-2 \Theta)e_4$& $e_1 e_2=e_4$ & $e_1e_3=e_4$& 
$e_2 e_1=e_3- \Theta^2 \lb^{-1}e_4$\\
$\lb \neq 0, \frac{1}{4}$&& $e_2 e_2=e_3$& $e_2e_3=(1-\Theta)\lb^{-1}e_4$ 
&$e_3e_3=e_4$\\

\hline$\cd {4}{91}(\lambda, \af)$&$:$& 
$e_1 e_1 = \lambda e_3+(2 \Theta-1)e_4$& $e_1 e_2=e_4$ & $e_1e_3=\af e_4$& \\
$\lb \neq 0, \frac{1}{4}$&& $e_2 e_1=e_3- (1-\Theta)^2 \lb^{-1}e_4$&
$e_2 e_2=e_3$&  $e_3e_3=e_4$\\

\hline$\cd {4}{92}(\lambda, \af)$&$:$& 
$e_1 e_1 = \lambda e_3+(1-2 \Theta)e_4$& $e_1 e_2=e_4$ & $e_1e_3=\af e_4$& \\
$\lb \neq 0, \frac{1}{4}$&& $e_2 e_1=e_3-  \Theta^2 \lb^{-1}e_4$&
$e_2 e_2=e_3$&  $e_3e_3=e_4$\\

\hline$\cd {4}{93}( \af)$&$:$& 
$e_1 e_1 = e_4$ & $e_1 e_2=e_4$ &$e_1e_3=\af e_4$ & $e_2e_1=e_3+e_4$\\
&&$ e_2 e_2=e_3$ &$e_2e_3=\af e_4$ & $e_3e_3=e_4$\\

\hline$\cd {4}{94}(\af, \bt)$&$:$& 
$e_1 e_1 = e_4$ & $e_1e_2=e_4$& $e_1e_3=\af e_4$&\\
$\af\neq0$&&$e_2 e_1=e_3+\beta e_4$   & $e_2 e_2=e_3$   &$e_3e_3=e_4$\\

\hline$\cd {4}{95}(\af)$&$:$& 
$e_1 e_1 =  e_4$ & $e_1e_2=e_4$ &$e_1e_3=\af e_4$ & $e_2 e_1=e_3$\\
&& $e_2 e_2=e_3$ & $e_2e_3=\af e_4$&$e_3e_3=e_4$\\

\hline$\cd {4}{96}(\af)$&$:$& 
$e_1 e_1 =  e_4$ & $e_1e_2=e_4$ & $e_2 e_1=e_3+\af e_4$\\
&& $e_2 e_2=e_3$ & $e_2e_3= e_4$&$e_3e_3=e_4$\\

\hline$\cd {4}{97}(\lambda)$&$:$& 
$e_1 e_1 = \lambda e_3$ & $e_1e_2=e_4$& $e_1e_3=\Theta e_4 $&$e_2 e_1=e_3-e_4$ \\
&& $e_2 e_2=e_3$ &$e_2e_3=e_4$ &$e_3e_3=e_4$\\

\hline$\cd {4}{98}(\lambda)$&$:$& 
$e_1 e_1 = \lambda e_3$ & $e_1e_2=e_4$& $e_1e_3=(1-\Theta) e_4 $&$e_2 e_1=e_3-e_4$ \\
$\lb \neq \frac{1}{4}$ && $e_2 e_2=e_3$ &$e_2e_3=e_4$ &$e_3e_3=e_4$\\

\hline$\cd {4}{99}(\af)$&$:$& 
$e_1e_2=e_4$& $e_1e_3=e_4$& $e_2 e_1=e_3-e_4$\\ 
$\af\neq1$ && $e_2 e_2=e_3$ &$e_2e_3=\af e_4$&$e_3e_3=e_4$\\

\hline$\cd {4}{100}(\af)$&$:$& 
$e_1 e_1 = \frac{1}{4} e_3$ & $e_1e_2=e_4$& $e_1e_3=\af e_4$& $e_2 e_1=e_3-e_4$ \\
$\af\notin\{0, \frac{1}{2}\}$&&  $e_2 e_2=e_3$ &$e_2e_3=2 \af e_4$ &$e_3e_3=e_4$\\

 \hline$\cd {4}{101}(\af, \bt)$&$:$& 
 $e_1e_2=e_4$& $e_1e_3=\af e_4$& $e_2 e_1=e_3$  \\
 && $e_2 e_2=e_3$&$e_2e_3=\bt e_4$ &$e_3e_3=e_4$\\
 
 \hline$\cd {4}{102}(\lb, \af)$&$:$& 
 $e_1 e_1 = \lambda e_3$ & $e_1e_2=e_4$& $e_2 e_1=e_3-e_4$ \\
 $\lb\neq0$&& $e_2 e_2=e_3$&$ e_2e_3=\af e_4$  &$e_3e_3=e_4$\\

\hline$\cd {4}{103}$&$:$& 
$e_1 e_2 = e_4$ & $e_2 e_1=e_3-e_4$  & $e_2 e_2=e_3$ &$e_3e_3=e_4$\\
 
 \hline$\cd {4}{104}$&$:$& 
 $e_1 e_3 =  e_4$ & $e_2 e_1=e_3+e_4$  & $e_2 e_2=e_3$ &$e_2e_3=e_4$&$e_3e_3=e_4$\\
 
 \hline$\cd {4}{105}(\lambda, \af,\bt)$&$:$& 
 $e_1 e_1 = \lambda e_3$ & $e_1e_3=e_4$&$e_2 e_1=e_3+\af e_4$  \\
 $  \lb\ne 0, \af\ne 0$&& $e_2 e_2=e_3$ & $e_2e_3=\bt e_4$&$e_3e_3=e_4$\\
 
 \hline$\cd {4}{106}(\af)$&$:$& $e_1 e_3 = e_4$ & $e_2 e_1=e_3+\af e_4$  & $e_2 e_2=e_3$ &$e_3e_3=e_4$\\
 
 \hline$\cd {4}{107}(\lambda)$&$:$& 
 $e_1 e_1 = \lambda e_3$ & $e_1e_3=\Theta e_4$& $e_2 e_1=e_3$  \\
 && $e_2 e_2=e_3$ & $e_2e_3=e_4$&$e_3e_3=e_4$\\
 
 \hline$\cd {4}{108}(\lambda)$&$:$& 
 $e_1 e_1 = \lambda e_3$ & $e_1e_3=(1-\Theta) e_4$& $e_2 e_1=e_3$  \\
 $\lb \not\in \{0, \frac{1}{4}\}$&& $e_2 e_2=e_3$ & $e_2e_3=e_4$&$e_3e_3=e_4$\\
 
 \hline$\cd {4}{109}(\lambda,\af)$&$:$& $e_1 e_1 = \lambda e_3$ & $e_2 e_1=e_3+e_4$  & $e_2 e_2=e_3$ &$e_2e_3=\alpha e_4$&$e_3e_3=e_4$\\

  \hline$\cd {4}{110}(\lambda)$&$:$& $e_1 e_1 = \lambda e_3$ & $e_2 e_1=e_3$  & $e_2 e_2=e_3$ &$e_2e_3= e_4$&$e_3e_3=e_4$\\
 
   \hline$\cd {4}{111}(\lambda)$&$:$& $e_1 e_1 = \lambda e_3$ & $e_2 e_1=e_3$  & $e_2 e_2=e_3$ &$e_3e_3=e_4$\\
 
\hline$\cd {4}{112}(\lambda, \af, \bt, \gamma)$&$:$& 
$e_1 e_1 = \lambda e_3+e_4$ & $e_1e_3=\af e_4$ & $e_2 e_1=e_3+\bt e_4$  \\
&&$e_2 e_2=e_3$& $e_2e_3=\gamma e_4$&$e_3e_3=e_4$\\
 
\end{longtable}
}
All these algebras are non-isomorphic, except

{\tiny
\begin{longtable}{lll}
$\cd {4}{91}(\lb, \af)\cong\cd {4}{91}(\lb, -\af)$ &
$\cd {4}{92}(\lb, \af)\cong\cd {4}{92}(\lb, -\af)$ & 
$\cd {4}{93}(\af)\cong\cd {4}{93}(-\af)$ \\
$\cd {4}{94}(\af,\bt)\cong\cd {4}{94}(-\af,\bt)$ &
$\cd {4}{95}(\af)\cong\cd {4}{95}(-\af)$ &
$\cd {4}{100}(\af)\cong\cd {4}{100}(-\af)$ \\
$\cd {4}{101}(\af,\bt)\cong\cd {4}{101}(-\af,-\bt)$ &
$\cd {4}{109}(\lb,\af)\cong\cd {4}{109}(\lb,-\af)$ &
$\cd {4}{112}(\lb,\af,\bt,\gm)\cong\cd {4}{112}(\lb,-\af,\bt,-\gm)$ \\
\multicolumn{3}{l}{$\cd {4}{112}(\lb,\af,\bt,\gm)\cong \cd {4}{112}\left(\lb,(\gm-\af\bt)\sqrt{\frac{-\lb}{1-\bt+\lb\bt^2}},\frac 1\lb-\bt,(\frac\gm\lb-\frac\af\lb-\bt\gm)\sqrt{\frac{-\lb}{1-\bt+\lb\bt^2}}\right)$, if $\lb\ne 0$, $\bt\ne\frac{1\pm\sqrt{1-4\lb}}{2\lb}$}
\end{longtable} 
}

\section{Classification theorem}
 
   \begin{theorem}\label{teor} 
Let $\mathfrak{CD}$ be a complex $4$-dimensional nilpotent $\mathfrak{CD}$-algebra.
Then $\mathfrak{CD}$ is isomorphic to an algebra from the following list:
 
{\tiny
\begin{longtable}{lllll llll}
\hline$\cd {4*}{01}$&$:$& $e_1 e_1 = e_2$\\
\hline$\cd {4*}{02}$&$:$& $e_1 e_1 = e_3$ &$ e_2 e_2=e_3$ \\
\hline$\cd {4*}{03}$&$:$& $e_1 e_2=e_3$ & $e_2 e_1=-e_3$ \\
\hline$\cd {4*}{04}(\lambda)$&$:$& $e_1 e_1 = \lambda e_3$ & $e_2 e_1=e_3$  & $e_2 e_2=e_3$\\

\hline$\cd {4*}{05}$ & $:$ & $e_1e_1 = e_3$&$ e_2e_2=e_4$ \\ 
\hline$\cd {4*}{06}$ & $:$ & $e_1e_2 = e_4$&$ e_3e_1 = e_4$   \\ 
\hline$\cd {4*}{07}$ & $:$ & $e_1e_2 = e_3$&$ e_2e_1 = e_4$&$  e_2e_2 = -e_3$\\ 
\hline$\cd {4*}{08}(\af)$ & $:$&$ e_1e_1 = e_3$& $e_1e_2 = e_4$&  $e_2e_1 = -\af e_3$&$e_2e_2 = -e_4$ \\
\hline$\cd {4*}{09}(\alpha)$&$:$&$e_1e_1 = e_4$&$ e_1e_2 = \alpha e_4$&$  e_2e_1 = -\alpha e_4$&$ e_2e_2 = e_4$&$  e_3e_3 = e_4$\\
\hline$\cd {4*}{10}$&$:$&$ e_1e_2 = e_4$&$ e_1e_3 = e_4$&$ e_2e_1 = -e_4$&$ e_2e_2 = e_4$&$ e_3e_1 = e_4$\\ 
\hline$\cd {4*}{11}$&$:$&$ e_1e_1 = e_4$&$ e_1e_2 = e_4$&$ e_2e_1 = -e_4$&$ e_3e_3 = e_4$\\
\hline$\cd {4*}{12}$&$:$&$ e_1e_2 = e_3$&$ e_2e_1 = e_4$  \\ 
\hline$\cd {4*}{13}$&$:$&$ e_1e_1 = e_4$&$ e_1e_2 = e_3$&$ e_2e_1 = -e_3$&$ e_2e_2=2e_3+e_4$\\
\hline$\cd {4*}{14}(\af)$&$:$&$ e_1e_2 = e_4$&$ e_2e_1 =\alpha e_4$&$ e_2e_2 = e_3$\\
\hline$\cd {4*}{15}$&$:$&$e_1e_2 = e_4$&$ e_2e_1 = -e_4$&$ e_3e_3 = e_4$\\
\hline

$\D{4}{00} $&$:$&  
$e_1e_1=e_4$ & $e_1e_2=e_4$& $e_2 e_1=e_3$ & $e_2 e_2 = e_3$ & $e_2e_3=e_4$ 
\\\hline

$\D{4}{01}(\lambda,\alpha,\beta)$&$:$& 
$e_1 e_1 = \lambda e_3 + e_4$ & $e_1 e_3 = \alpha e_4$ & $e_2 e_1=e_3$ \\
&&$ e_2 e_2 = e_3$ & $e_2 e_3 = \beta e_4$ & $e_3e_1 = e_4$ \\
\hline
$\D{4}{02}(\lambda,\alpha,\beta)$&$:$& 
$e_1 e_1 = \lambda e_3$ & $e_1 e_3 = \alpha e_4$ & $e_2 e_1=e_3$ \\
&& $e_2 e_2 = e_3$ & $e_2 e_3 = \beta e_4$ & $e_3e_1 = e_4$ \\
\hline
$\D{4}{03}(\lambda,\alpha)$&$:$& 
$e_1 e_1 = \lambda e_3$ & $e_1 e_2 = e_4$ & $e_2 e_1=e_3$ \\
&& $e_2 e_2 = e_3$ & $e_2 e_3 = \alpha e_4$ & $e_3e_1 = e_4$ \\
\hline
$\D{4}{04}(\lambda,\alpha)$&$:$& 
$e_1 e_1 = \lambda e_3$ & $e_2 e_1=e_3$ & $e_2 e_2 = e_3$ \\
&& $e_2 e_3 = \alpha e_4$ & $e_3e_1 = e_4$ \\
\hline
$\D{4}{05}(\lambda,\alpha)$&$:$& 
$e_1 e_1 = \lambda e_3$ & $e_1 e_2 = \lambda e_4$ & $e_2 e_1=e_3 + \lambda\alpha e_4$ \\
&& $e_2 e_2 = e_3$ & $e_2 e_3 = 	\Theta e_4$ & $e_3e_1 = \lambda e_4$ \\
	    %\D{4}{05}(\lambda\alpha)&:& e_1 e_1 = \lambda e_3 & e_1 e_2 = e_4 & e_2 e_1=e_3 + \alpha e_4 & e_2 e_2 = e_3 & e_2 e_3 = 	\frac{1+\sqrt{1-4\lb}}{2\lb} e_4 & e_3e_1 = e_4 \\
\hline
$\D{4}{06}(\lambda,\alpha)$&$:$& 
$e_1 e_1 = \lambda e_3$ & $e_1 e_2 = e_4$ & $e_2 e_1=e_3 + \alpha e_4$ \\
&& $e_2 e_2 = e_3$ & $e_2 e_3 = 	\Theta^{-1} e_4$ & $e_3e_1 = e_4$ \\
\hline
	    %\D{4}{06}(\lambda\alpha)&:& e_1 e_1 = \lambda e_3 & e_1 e_2 = e_4 & e_2 e_1=e_3 + \alpha e_4 & e_2 e_2 = e_3 & e_2 e_3 = 	\frac{1-\sqrt{1-4\lb}}{2\lb} e_4 & e_3e_1 = e_4 \\
$\D{4}{07}(\lambda)$&$:$& $e_1 e_1 = \lambda e_3$ & $e_2 e_1=e_3 + \lambda e_4$ & $e_2 e_2 = e_3$ & $e_2 e_3 = 	\Theta e_4$ & $e_3e_1 = \lambda e_4$ \\
\hline
	    %\D{4}{07}(\lambda)&:& e_1 e_1 = \lambda e_3 & e_2 e_1=e_3 + e_4 & e_2 e_2 = e_3 & e_2 e_3 = 	\frac{1+\sqrt{1-4\lb}}{2\lb} e_4 & e_3e_1 = e_4 \\
$\D{4}{08}(\lambda)$&$:$& 
$e_1 e_1 = \lambda e_3$ &  $e_2 e_1=e_3 + e_4$ & $e_2 e_2 = e_3$ & $e_2 e_3 = 	\Theta^{-1} e_4$ & $e_3e_1 = e_4$ \\
\hline
	    %\D{4}{08}(\lambda)&:& e_1 e_1 = \lambda e_3 &  e_2 e_1=e_3 + e_4 & e_2 e_2 = e_3 & e_2 e_3 = 	\frac{1-\sqrt{1-4\lb}}{2\lb} e_4 & e_3e_1 = e_4 \\
$\D{4}{09}(\lambda,\alpha)$&$:$& 
$e_1 e_1 = \lambda e_3$ & $e_1 e_2 = e_4$& $e_1 e_3 = \alpha e_4$ \\
&& $e_2 e_1=e_3$ & $e_2 e_2 = e_3 $& $e_3e_1 = e_4$ \\
\hline
$\D{4}{10}(\lambda,\alpha)$&$:$& 
$e_1 e_1 = \lambda e_3$& $e_1 e_3 = \alpha e_4$ &  $e_2 e_1=e_3$ & $e_2 e_2 = e_3$ & $e_3e_1 = e_4$ \\
\hline        
$\D{4}{11}(\lambda,\alpha)$&$:$& 
$e_1 e_1 = \lambda e_3 + e_4$& $e_1e_2 = \alpha e_4$ & $e_1 e_3 = -e_4$ \\
&&  $e_2 e_1=e_3$ & $e_2 e_2 = e_3$ & $e_3e_1 = e_4$ \\
\hline
$\D{4}{12}(\lambda,\alpha)$&$:$& 
$e_1 e_1 = \lambda e_3 + e_4$& $e_1e_3 = \alpha e_4$  &  $e_2 e_1=e_3$ \\
&& $e_2 e_2 = e_3$ & $e_3e_1 = \Theta e_4$ & $e_3e_2 = e_4$ \\
\hline
$\D{4}{13}(\lambda,\alpha)$&$:$& 
$e_1 e_1 = \lambda e_3 + e_4$& $e_1e_3 = \alpha e_4$  &  $e_2 e_1=e_3$ \\
&& $e_2 e_2 = e_3$ & $e_3e_1 = (1-\Theta)e_4$ & $e_3e_2 = e_4$ \\
\hline
$\D{4}{14}(\lambda,\alpha)$&$:$& 
$e_1 e_1 = \lambda e_3$& $e_1e_3 = \alpha e_4 $ &  $e_2 e_1=e_3$ \\
&& $e_2 e_2 = e_3$ & $e_3e_1 = \Theta e_4$ & $e_3e_2 = e_4$ \\
\hline
$\D{4}{15}(\lambda,\alpha)$&$:$& 
$e_1 e_1 = \lambda e_3$& $e_1e_3 = \alpha e_4$  &  $e_2 e_1=e_3$ \\
&& $e_2 e_2 = e_3$ & $e_3e_1 = (1-\Theta)e_4$ & $e_3e_2 = e_4$ \\
\hline
$\D{4}{16}(\alpha)$&$:$& $e_1e_3 = \alpha e_4$  &  $e_2 e_1=e_3 + e_4$ & $e_2 e_2 = e_3$ & $e_3e_1 = e_4$ & $e_3e_2 = e_4$ \\
\hline
$\D{4}{17}(\alpha)$&$:$& 
$e_1 e_1 = e_4$ & $e_1e_3 = -e_4$  &  $e_2 e_1=e_3 + \alpha e_4$ \\
&& $e_2 e_2 = e_3$ & $e_3e_1 = e_4$ & $e_3e_2 = e_4$ \\
\hline
$\D{4}{18}(\lambda,\alpha)$&$:$& $e_1 e_1 = \lambda e_3 + e_4$& $e_1e_3 = \Theta\alpha e_4$  &  $e_2 e_1=e_3$ & $e_2 e_2 = e_3$ \\ && $e_2 e_3 = \alpha e_4$ & $e_3e_1 = \Theta e_4$ & $e_3e_2 = e_4$ \\

\hline
$\D{4}{19}(\lambda,\alpha)$&$:$& $e_1 e_1 = \lambda e_3 + e_4$& $e_1e_3 = (1-\Theta)\alpha e_4$  &  $e_2 e_1=e_3$ & $e_2 e_2 = e_3$ \\&& $e_2 e_3 = \alpha e_4$ & $e_3e_1 = (1-\Theta)e_4$& $e_3e_2 = e_4$ \\
\hline   
$\D{4}{20}(\lambda,\alpha)$&$:$& $e_1 e_1 = \lambda e_3$& $e_1e_3 = \Theta\alpha e_4$  &  $e_2 e_1=e_3$ & $e_2 e_2 = e_3$ \\
&& $e_2 e_3 = \alpha e_4$ & $e_3e_1 = \Theta e_4$& $e_3e_2 = e_4$ \\
\hline
$\D{4}{21}(\lambda,\alpha)$&$:$& $e_1 e_1 = \lambda e_3$& $e_1e_3 = (1-\Theta)\alpha e_4$  &  $e_2 e_1=e_3$ & $e_2 e_2 = e_3$ \\&& $e_2 e_3 = \alpha e_4$ & $e_3e_1 = (1-\Theta)e_4$ & $e_3e_2 = e_4$ \\
\hline
$\D{4}{22}(\lambda)$&$:$& $e_1 e_1 = \lambda e_3 + (1-2\lambda)e_4$& $e_1 e_2 = e_4$& $e_1e_3 = (\Theta - 1) e_4$  &  $e_2 e_1=e_3$ \\&& $e_2 e_2 = e_3$ & $e_2 e_3 = (1-\Theta^{-1}) e_4$& $e_3e_1 = \Theta e_4$ & $e_3e_2 = e_4$ \\
\hline
$\D{4}{23}(\lambda)$&$:$& $e_1 e_1 = \lambda e_3 + \lambda(1-2\lambda) e_4$& $e_1 e_2 = \lambda e_4$& $e_1e_3 = -\lambda\Theta e_4 $&  $e_2 e_1=e_3$ \\&& 
$e_2 e_2 = e_3$ & $e_2 e_3 = -\Theta^2 e_4$ & $e_3e_1 = \lambda(1-\Theta)e_4$ & $e_3e_2 = \lambda e_4$ \\
\hline        
$\D{4}{24}(\lambda)$&$:$& $e_1 e_1 = \lambda e_3 + \Theta e_4$& $e_1 e_2 = e_4$& $e_1e_3 = (\Theta - 1) e_4$  &  $e_2 e_1=e_3$ \\&& $e_2 e_2 = e_3$ & $e_2 e_3 = (1-\Theta^{-1}) e_4$ & $e_3e_1 = \Theta e_4$ & $e_3e_2 = e_4$ \\
\hline
$\D{4}{25}(\lambda)$&$:$& $e_1 e_1 = \lambda e_3 + \lambda(1-\Theta)e_4$& $e_1 e_2 = \lambda e_4$& $e_1e_3 = -\lambda\Theta e_4$  &  $e_2 e_1=e_3$ \\&& $e_2 e_2 = e_3$ & $e_2 e_3 = -\Theta^2 e_4$ 
        & $e_3e_1 = \lambda (1-\Theta)e_4$ & $e_3e_2 = \lambda e_4$ \\
\hline
$\D{4}{26}(\lambda)$&$:$& $e_1 e_1 = \lambda e_3 + \Theta e_4$& $e_1 e_2 = e_4$& $e_1e_3 = -\Theta e_4$  &  $e_2 e_1=e_3$\\
& & $e_2 e_2 = e_3$ & $e_2 e_3 = -e_4$
         & $e_3e_1 = \Theta e_4$ & $e_3e_2 = e_4 $\\
\hline
$\D{4}{27}(\lambda)$&$:$& $e_1 e_1 = \lambda e_3 + (1-\Theta)e_4$& $e_1 e_2 = e_4$& $e_1e_3 = (\Theta-1) e_4$  &  $e_2 e_1=e_3$ \\&& $e_2 e_2 = e_3$ & $e_2 e_3 = -e_4$  & $e_3e_1 = (1-\Theta)e_4$ & $e_3e_2 = e_4$ \\
\hline
$\D{4}{28}(\lambda)$&$:$& $e_1 e_1 = \lambda e_3 + (1-\Theta)e_4$& $e_1 e_2 = e_4$& $e_1e_3 = -\Theta e_4$  &  $e_2 e_1=e_3$ \\&& $e_2 e_2 = e_3$ & $e_2 e_3 = -e_4$ & $e_3e_1 = \Theta e_4$ & $e_3e_2 = e_4$ \\
\hline
$\D{4}{29}(\lambda)$&$:$& $e_1 e_1 = \lambda e_3 + \Theta e_4$& $e_1 e_2 = e_4$& $e_1e_3 = (\Theta-1) e_4$  &  $e_2 e_1=e_3$ \\
&& $e_2 e_2 = e_3$ & $e_2 e_3 = -e_4$ & $e_3e_1 = (1-\Theta)e_4$ & $e_3e_2 = e_4$ \\
$\D{4}{30}(\lambda)$&$:$& $e_1 e_1 = \lambda e_3 + e_4$ & $e_1e_3 = (\Theta-1) e_4$  &  $e_2 e_1=e_3$ & $e_2 e_2 = e_3$ \\
& & $e_2 e_3 = -e_4$ & $e_3e_1 = \Theta e_4$ & $e_3e_2 = e_4$ \\
\hline
$\D{4}{31}(\lambda)$&$:$& $e_1 e_1 = \lambda e_3 + e_4$ & $e_1e_3 = -\Theta e_4$  &  $e_2 e_1=e_3$ & $e_2 e_2 = e_3$ \\
&& $e_2 e_3 = -e_4$ & $e_3e_1 = (1-\Theta)e_4$ & $e_3e_2 = e_4$ \\
\hline        
$\D{4}{32}(\lambda)$&$:$& $e_1 e_1 = \lambda e_3$ & $e_1e_3 = (\Theta-1) e_4$  &  $e_2 e_1=e_3$ & $e_2 e_2 = e_3$ \\
&& $e_2 e_3 = -e_4$ & $e_3e_1 = \Theta e_4$ & $e_3e_2 = e_4$ \\
\hline
$\D{4}{33}(\lambda)$&$:$& $e_1 e_1 = \lambda e_3$ & $e_1e_3 = -\Theta e_4$  & $ e_2 e_1=e_3$ & $e_2 e_2 = e_3$ \\
&& $e_2 e_3 = -e_4$ & $e_3e_1 = (1-\Theta)e_4$ & $e_3e_2 = e_4$ \\
\hline 
$\D{4}{34}$&$:$& $e_1 e_1 = e_4$ &  $e_2 e_1=e_3 + e_4$ & $e_2 e_2 = e_3$ & $e_3 e_1 = e_4$& $e_3 e_2 = e_4$ \\
\hline
$\D{4}{35}(\lambda)$&$:$& $e_1 e_1 = \lambda e_3$ & $e_1e_3 = e_4$  &  $e_2 e_1=e_3 + e_4$ & $e_2 e_2 = e_3$ \\
\hline
$\D{4}{36}(\lambda)$&$:$& $e_1 e_1 = \lambda e_3$ & $e_1e_3 = e_4$  &  $e_2 e_1=e_3$ & $e_2 e_2 = e_3$ \\
\hline 
$\D{4}{37}(\lambda)$&$:$& $e_1 e_1 = \lambda e_3 + e_4$ & $e_1e_3 = \Theta e_4$ &  $e_2 e_1=e_3$ & $e_2 e_2 = e_3$&  $e_2 e_3 = e_4$ \\
\hline
$\D{4}{38}(\lambda)$&$:$& $e_1 e_1 = \lambda e_3 + e_4$ & $e_1e_3 = (1-\Theta)e_4$ &  $e_2 e_1=e_3$ & $e_2 e_2 = e_3$&  $e_2 e_3 = e_4$ \\
\hline
$\D{4}{39}(\lambda)$&$:$& $e_1 e_1 = \lambda e_3$ & $e_1e_3 = \Theta e_4$ &  $e_2 e_1=e_3$ & $e_2 e_2 = e_3$&  $e_2 e_3 = e_4$ \\
\hline
$\D{4}{40}(\lambda)$&$:$& $e_1 e_1 = \lambda e_3$ & $e_1e_3 = (1-\Theta)e_4$ &  $e_2 e_1=e_3$ & $e_2 e_2 = e_3$&  $e_2 e_3 = e_4$\\
\hline
$\cd 4{01}$&$:$& $e_1 e_1 = e_2$ & $e_2 e_2=e_3$ \\
\hline
$\cd 4{02}$&$:$& $e_1 e_1 = e_2$ & $e_2 e_1= e_3$ & $e_2 e_2=e_3$ \\
\hline$\cd 4{03}$&$:$& $e_1 e_1 = e_2$ & $e_2 e_1=e_3$ \\
\hline$\cd 4{04}(\lambda)$&$:$& $ e_1 e_1 = e_2$ & $e_1 e_2=e_3$ & $e_2 e_1=\lambda e_3$ \\ \hline

$\cd 4{05}$ & $:$ & $e_1 e_1 = e_2$ & $e_2 e_1=e_4$ & $e_2 e_2=e_3$ \\
\hline$\cd 4{06}$ & $:$ & $ e_1 e_1 = e_2$ & $e_1 e_2=e_4$ & $e_2 e_1=e_3$  \\

\hline$\cd 4{07}(\lambda)$& $:$ & $e_1 e_1 = e_2$ & $e_1 e_2=e_4$ & $e_2 e_1=\lambda e_4$ & $e_2 e_2=e_3$ \\\hline
$\cd 4{08}(\alpha)$ &$:$&
$e_1 e_1 = e_2$ & $e_1e_3=e_4$ & $e_2 e_1=e_3$ & $e_2e_2=\alpha e_4$ & $e_3e_1=-2e_4$ \\

\hline$\cd 4{09}$ &$:$& 
$e_1 e_1 = e_2$ & $e_1e_2=e_4$ & $e_1e_3=e_4$ \\
&& $e_2 e_1=e_3$ & $e_2e_2=- e_4$ & $e_3e_1=-2e_4$\\
\hline

$\cd 4{10}(\alpha)$ &$:$&
$e_1 e_1 = e_2$ & $e_1 e_2=e_3$ & $e_1e_3=-e_4$ \\
&& $e_2 e_1=  e_3 +e_4$ & $e_2e_2=\alpha e_4$ & $e_3e_1=- e_4$ \\

\hline$\cd 4{11}(\lambda)$ &$:$& 
$e_1 e_1 = e_2$ & $e_1 e_2=e_3$ & $e_1e_3=(\lambda-2)e_4$ \\
$\lambda \neq1$&& $e_2 e_1=  \lb e_3+e_4$  &
$e_2e_2=-(\lambda+1)^2 e_4$ &$e_3e_1= (1-2\lb) e_4$ \\

\hline$\cd 4{12}(\alpha, \lambda)$ &$:$& 
$e_1 e_1 = e_2$ & $e_1 e_2=e_3$ & $e_1e_3=(\lambda-2)e_4$ \\&& $e_2 e_1=\lambda e_3$ & $e_2e_2=\alpha e_4$ & $e_3e_1=(1-2\lambda) e_4$ \\
\hline

$\cd {4}{13}(\af)$&$:$& 
$e_1 e_1 = e_2$ & $e_1e_2=e_4$ & $e_1e_3=e_4$& $e_2e_1=e_4$ \\
$\af\neq \frac{1}{2}$&&  $e_2e_3=\alpha e_4$& $e_3e_1=e_4$& \multicolumn{2}{l}{$e_3e_2=(\alpha +1)e_4$}\\

\hline
$\cd {4}{14}(\af, \beta)$&$:$& 
$e_1 e_1 = e_2$ & $e_1 e_2 = e_4$&  $e_1 e_3 = \af e_4$& $e_2 e_1 = e_4$\\ && 
$e_2 e_2 = e_4$& $e_3 e_1 = \af e_4$& $e_3 e_2 = e_4$& $e_3 e_3 =\beta  e_4$& \\

\hline
$\cd {4}{15}(\af)$&$:$& 
$e_1 e_1 = e_2$ &  $e_1 e_2 = \af e_4$&   $e_1 e_3 = e_4$&
 $e_2 e_1 =(\af+1)  e_4$ &  $e_3 e_1 = e_4$\\

\hline
$\cd {4}{16}$&$:$& 
$e_1 e_1 = e_2$ & $e_1e_2=e_4$ & $e_2 e_1 = e_4$\\ && $e_2 e_3 = -\frac{1}{2} e_4$& $e_3 e_2 =\frac{1}{2} e_4$& $e_3 e_3 = e_4$&\\

\hline
$\cd {4}{17}(\af)$&$:$& 
$e_1 e_1 = e_2$ &  $e_1 e_2 = e_4$&  $e_2 e_1 = e_4$&  $e_2 e_3 =\af e_4$&  $e_3 e_2 = (\af+1) e_4$&\\

\hline
$\cd {4}{18}(\af)$&$:$& 
$e_1 e_1 = e_2$ &  $e_1 e_2 =\af  e_4$&   \multicolumn{2}{l}{$e_2 e_1 =(\af+1) e_4$}&   $e_3 e_3 = e_4$&\\

\hline
$\cd {4}{19}$&$:$& 
$e_1 e_1 = e_2$ &   $e_1 e_2 = e_4$&   $e_2 e_1 = e_4$&    $e_3 e_1 = e_4$&
 $e_3 e_3 = e_4$&\\
 
\hline
$\cd {4}{20}$&$:$& 
$e_1 e_1 = e_2$ &    $e_1 e_2 = e_4$&   $e_2 e_1 = e_4$&   $e_3 e_1 = e_4$&\\

\hline
$\cd {4}{21}(\af)$&$:$& 
$e_1 e_1 = e_2$ & $e_1 e_3 =\af e_4$&   $e_2 e_1 = e_4$&  $e_2 e_2 = e_4$&\\
&&  $e_2 e_3 = e_4$&  $e_3 e_1 =\af e_4$&  $e_3 e_2 = e_4$&  $e_3 e_3 = e_4$\\

\hline
$\cd {4}{22}$&$:$& 
$e_1 e_1 = e_2$ &   $e_1 e_3 = e_4$&    $e_2 e_3 = -\frac{1}{2} e_4$\\&&
 $e_3 e_1= e_4$&   $e_3 e_2 =\frac{1}{2} e_4$&   $e_3 e_3 = e_4$&\\

\hline
$\cd {4}{23}(\af)$&$:$& 
  $e_1 e_1 = e_2$&   $e_1 e_3 = e_4$&   $e_2 e_3 = \af e_4$&
  $e_3 e_1 = e_4$&   $e_3 e_2 =(\af +1)e_4$&\\

\hline
$\cd {4}{24}(\af)$&$:$& 
  $e_1 e_1 = e_2$&   $e_1 e_3 = e_4$&    $e_2 e_2 = e_4$\\
  && $e_3 e_1 = e_4$&   $e_3 e_2 = e_4$& $e_3 e_3 =\af  e_4$&\\

\hline
$\cd {4}{25} $&$:$& 
  $e_1 e_1 = e_2$&  $e_1 e_3 = e_4$& $e_2 e_1 = e_4$& $e_3 e_1 = e_4$&    $e_2 e_2 = e_4$&\\

\hline
$\cd {4}{26}(\af)$&$:$& 
  $e_1 e_1 = e_2$&   $e_1 e_3 = \af e_4$&    $e_2 e_2 = e_4$& \multicolumn{2}{l}{$e_3 e_1 =(\af+1) e_4$}&\\

%\hline
%$\cd {4}{27}(\af)$&$:$& 
%  $e_1 e_1 = e_2$& $e_1 e_3 = \af e_4$& \multicolumn{2}{l}{$e_3 e_1 =(\af +1) e_4$}&\\
  
  \hline
$\cd {4}{27}$&$:$& 
  $e_1 e_1 = e_2$&   $e_2 e_1 = e_4$& $e_2 e_3 = e_4$&  $e_3 e_2 = e_4$&   $e_3 e_3 = e_4$&\\

  \hline
$\cd {4}{28}(\af)$&$:$& 
  $e_1 e_1 = e_2$&   $e_2 e_1 = e_4$&    $e_2 e_2 = e_4$\\
  $\af\ne 1$&& $e_2 e_3 =   e_4$&   $e_3 e_2 =  e_4$&   $e_3 e_3 =\af  e_4$&\\

  \hline
$\cd {4}{29}$&$:$& 
  $e_1 e_1 = e_2$&    $e_2 e_3 =-\frac{1}{2} e_4$&   $e_3 e_2 =\frac{1}{2} e_4$&   $e_3 e_3 = e_4$&\\

  \hline
$\cd {4}{30}$&$:$& 
  $e_1 e_1 = e_2$&   $e_2 e_1 = e_4$& $e_2 e_3 = e_4$&   $e_3 e_2 = e_4$&\\

  \hline
$\cd {4}{31}$&$:$& 
  $e_1 e_1 = e_2$&  $e_2 e_3 = e_4$&     $e_3 e_1 = e_4$&
$e_3 e_2 = e_4$&\\

  \hline
$\cd {4}{32}(\af)$&$:$& 
  $e_1 e_1 = e_2$&   $e_2 e_3 = \af e_4$&   \multicolumn{2}{l}{$e_3 e_2 = (\af+1) e_4$}& \\

  \hline
$\cd {4}{33}$&$:$& 
  $e_1 e_1 = e_2$&   $e_2 e_1 = e_4$&   $e_2 e_2 = e_4$& $e_3 e_3 = e_4$& \\

 \hline
$\cd {4}{34}$&$:$& 
  $e_1 e_1 = e_2$&    $e_2 e_2 = e_4$& $e_3 e_3 = e_4$& \\

 \hline
$\cd {4}{35}$&$:$& 
  $e_1 e_1 = e_2$&    $e_2 e_2 = e_4$&  $e_3 e_1 = e_4$& $e_3 e_3 = e_4$& \\

 \hline
$\cd {4}{36}(\af)$&$:$& 
  $e_1 e_1 = e_2$&     $e_2 e_2 = e_4$&   $e_3 e_2= e_4$& $e_3 e_3 =\af  e_4$&\\

 %\hline
%$\cd {4}{38}$&$:$& 
%  $e_1 e_1 = e_2$&    $e_3 e_1 = e_4$&    $e_3 e_3 = e_4$&\\

 \hline
$\cd {4}{37}$&$:$& 
  $e_1 e_1= e_2$&    $e_1 e_2= e_4$&   $e_2 e_1= e_4$&    $e_3 e_3= e_4$\\

%  \hline
%$\cd {4}{40}$&$:$& 
%  $e_1 e_1= e_2$&   $e_1 e_3= e_4$&   $e_3 e_1= e_4$&\\

  \hline
$\cd {4}{38}$&$:$& 
  $e_1 e_1= e_2$&   $e_2 e_3= e_4$&   $e_3 e_2= e_4$&\\

%\hline
%$\cd {4}{42}$&$:$& 
%  $e_1 e_1= e_2$&      $e_3 e_3= e_4$\\
  
\hline
$\cd 4{39}$ & $:$ &  
$e_1 e_1 = e_3+e_4$ & $e_1e_2=\frac i2 e_4$ & $e_1e_3=e_4$ & $e_2e_1=\frac i2 e_4$\\
&& $ e_2 e_2=e_3$ & $e_2e_3=-2ie_4$ & $e_3e_1=2e_4$ & $e_3e_2=-ie_4$\\

\hline
$\cd 4{40}$ & $:$ &  
$e_1 e_1 = e_3+e_4$ & $e_1e_2=\frac i2 e_4$ & $e_1e_3=-\frac 12e_4$ & $e_2e_1=\frac i2 e_4$\\
&& $ e_2 e_2=e_3$ & $e_2e_3=-\frac i2e_4$ & $e_3e_1=\frac 12e_4$ & $e_3e_2=\frac i2e_4$\\

\hline
$\cd 4{41}$ & $:$ &  
$e_1 e_1 = e_3+e_4$ & $e_1e_2=-\frac i2 e_4$ & $e_1e_3=e_4$ & $e_2e_1=-\frac i2 e_4$\\
&& $ e_2 e_2=e_3$ & $e_2e_3=-2ie_4$ & $e_3e_1=2e_4$ & $e_3e_2=-ie_4$\\

\hline
 $\cd 4{42}$ & $:$ &  
$e_1 e_1 = e_3+e_4$ & $e_1e_2=-\frac i2 e_4$ & $e_1e_3=-\frac 12e_4$ & $e_2e_1=-\frac i2 e_4$\\
&& $ e_2 e_2=e_3$ & $e_2e_3=-\frac i2e_4$ & $e_3e_1=\frac 12e_4$ & $e_3e_2=\frac i2e_4$\\

\hline
$\cd 4{43}(\af)$ &$:$&  
$e_1 e_1 = e_3 + e_4$     & $e_1 e_2 = \af e_4$          & $e_1 e_3 = -\frac 12 e_4$ \\
&&   $e_2 e_1 = \af e_4$      & $e_2 e_2 = e_3$           & $e_3 e_1 = \frac 12 e_4$\\

\hline
$\cd 4{44}(\af,\bt,\gm)$ & $:$ &  
$e_1 e_1 = e_3 + \af e_4$   & $e_1 e_2 = \bt e_4$ & $e_2 e_1 = (\bt+\gm) e_4$ \\
&& $e_2 e_2 = e_3$ & $e_3 e_1 = e_4$             & $e_3 e_3 = e_4$\\

\hline
$\cd 4{45}$ & $:$ &  
$e_1 e_1 = e_3 + 2i e_4$  & $e_1 e_2 = e_4$ & $e_2 e_1 = e_4$& $e_2 e_2 = e_3$\\
&&  $e_3 e_1 = e_4$           & $e_3 e_2 = i e_4$       & $e_3 e_3 = e_4$\\

\hline
$\cd 4{46}(\af)$ & $:$ &
$e_1 e_1 = e_3 - 2i\af e_4$ & $e_1 e_2 = \af e_4$  & $e_2 e_1 = \af e_4$ & $e_2 e_2 = e_3$ \\
$\af\ne 0$& & $e_3 e_1 = e_4$              & $e_3 e_2 = i e_4$       & $e_3 e_3 = e_4$\\

\hline
$\cd 4{47}(\af,\bt)$ & $:$ & 
$e_1 e_1 = e_3 + e_4$  & $e_1 e_3 = \af e_4$  & $e_2 e_2 = e_3$\\ 
$\bt\ne 0$ && $e_2 e_3 = \bt e_4$  & $e_3 e_1 =(\af+1) e_4$        & $e_3 e_2 = \bt e_4$ \\

\hline		
$\cd 4{48}(\af)$ & $:$ & 
 $e_1 e_1 = e_3 + \af e_4$ &  $e_2 e_1 = i\af e_4$ & $e_2 e_2 = e_3$  \\
$\af\ne$ 0&& $e_3 e_1 = e_4$            & $e_3 e_2 = i e_4$     & $e_3 e_3 = e_4$\\

\hline 
$\cd 4{49}(\af)$ & $:$ & 
$e_1 e_1 = e_3 + \af e_4$ & $e_2 e_1 = -i\af e_4$ & $e_2 e_2 = e_3$\\
$\af\ne 0$& & $e_3 e_1 = e_4$& $e_3 e_2 = i e_4$ & $e_3 e_3 = e_4$\\

\hline 
$\cd 4{50}(\af)$ & $:$ & 
$e_1 e_1 = e_3 + \af e_4$ & $e_2 e_1 = e_4$ & $e_2 e_2 = e_3$ & $e_3 e_3 = e_4$\\

\hline
$\cd 4{51}(\af)$ & $:$ &
$e_1 e_1 = e_3 + \af e_4$ &   $e_2 e_2 = e_3$ & $e_3 e_1 = e_4$\\
&& $e_3 e_2 = i e_4$ & $e_3 e_3 = e_4$\\

\hline 
$\cd 4{52}$ & $:$ &  
$e_1 e_1 = e_3 +  e_4$ &   $e_2 e_2 = e_3$ &   $e_3 e_3 = e_4$\\

\hline 
$  \cd 4{53}$ & $:$ &  
$e_1 e_1 = e_3$ &  $e_1e_2=-\frac 12e_4$ & $e_1e_3=e_4$ & $e_2e_1=\frac 12e_4$\\
 && $e_2 e_2 = e_3$ & $e_2e_3=ie_4$ &  $e_3e_1=e_4$ & $e_3e_2=ie_4$ \\

\hline
 $\cd 4{54}(\af)$ & $:$ &  
$e_1 e_1 = e_3$ & $e_1e_2=e_4$ & $e_1e_3=\af e_4$ & $e_2e_1=e_4$\\
&& $ e_2 e_2=e_3$ & $e_2e_3=-i(\af+1)e_4$ & $e_3e_1=(\af+1) e_4$ & $e_3e_2=-i\af e_4$\\

\hline
$\cd 4{55}(\af)$ & $:$ &  
$e_1 e_1 = e_3$ & $e_1 e_2 = e_4$ & $e_1 e_3 = \af e_4$ \\
& & $e_2 e_1 = e_4$ & $e_2 e_2 = e_3$  & $e_3 e_1 = (\af+1) e_4$ \\
		
\hline
$\cd 4{56}$ & $:$ &  
$e_1 e_1 = e_3$ & $e_1 e_2 = e_4$ & $e_2 e_1 = -e_4$ & $e_2 e_2 = e_3$ \\
& & $e_3 e_1 = e_4$ & $e_3 e_2 = i e_4$ & $e_3 e_3 = e_4$\\

\hline
$\cd 4{57}(\af,\bt)$ & $:$ &  
$e_1 e_1 = e_3$ & $e_1 e_2 = \af e_4$ &   $e_2 e_1 = (\af+\bt)e_4$ & $e_2 e_2 = e_3$ \\
$\bt\not\in\{0,-2\af\}$& & $e_3 e_1 = e_4$          & $e_3 e_2 = i e_4$     & $e_3 e_3 = e_4$\\

\hline
$\cd 4{58}$ & $:$ &  
$e_1 e_1 = e_3$ & $e_1 e_3 = e_4$  & $e_2 e_1 = e_4$ & $e_2 e_2 = e_3$\\
&& $e_2 e_3 = i e_4$  & $e_3 e_1 = e_4$ & $e_3 e_2 = i e_4$ \\

\hline 
$\cd 4{59}(\af,\bt)$ & $:$ &   
$e_1 e_1 = e_3$ & $e_1 e_3 = \af e_4$ & $e_2 e_2 = e_3$\\
 $\bt\ne 0$ && $e_2 e_3 = \bt e_4$ & $e_3 e_1 = (\af+1) e_4$  & $e_3 e_2 = \bt e_4$ \\

\hline		
$\cd 4{60}$ & $:$ & 
$e_1 e_1 = e_3$ & $e_1 e_3 = i e_4$  & $e_2 e_2 = e_3$\\
&& $e_2 e_3 = e_4$ & $e_3 e_1 = (i+1) e_4$ & $e_3 e_2 = (i+1) e_4$  \\

\hline
$\cd 4{61}(\af)$ & $:$ & 
 $e_1 e_1 = e_3$ & $e_1 e_3 = -i\af e_4$ & $e_2 e_2 = e_3$\\
 && $e_2 e_3 = \af e_4$ & $e_3 e_1 = (1-i\af) e_4$ & $e_3 e_2 = (\af+i) e_4$  \\

\hline
$\cd 4{62}$ & $:$ &  
$e_1 e_1 = e_3$  & $e_1e_3=e_4$ & $e_2 e_2=e_3$\\
&& $e_2e_3=-2ie_4$ & $e_3e_1=2e_4$ & $e_3e_2=-ie_4$\\

\hline
$\cd 4{63}$ & $:$ &  
$e_1 e_1 = e_3$ & $e_1e_3=-\frac 12e_4$ & $ e_2 e_2=e_3$\\
&& $e_2e_3=-\frac i2e_4$ & $e_3e_1=\frac 12e_4$ & $e_3e_2=\frac i2e_4$\\

\hline 
$\cd 4{64}(\af)$ & $:$ & 
$e_1 e_1 = e_3$ & $e_1 e_3 = \af e_4$ & $e_2 e_2 = e_3$  & $e_3 e_1 = (\af+1) e_4$ 		\\

\hline
$\cd 4{65}(\af)$ & $:$ &  
$e_1 e_1 = e_3$ &   $e_1 e_3 = \af e_4$ & $e_2 e_2 = e_3$ \\
$\af\ne 0$&&  $e_3 e_1 = (\af+1) e_4$  & $e_3 e_2 = i e_4$ \\
\hline
$\cd 4{66}$ & $:$ &  
$e_1 e_1 = e_3$& $e_2 e_1 = e_4$ & $e_2 e_2 = e_3$ \\
&& $e_2 e_3 = e_4$ & $e_3 e_2 = e_4$ \\

\hline 
$\cd 4{67}$ & $:$ & 
$e_1 e_1 = e_3$ &  $e_2 e_2 = e_3$ &  $e_3 e_3 = e_4$\\

\hline
$\cd 4{68}$ & $:$ & 
$e_1 e_1 = e_3+e_4$ & $e_1 e_3 = i e_4$  &$e_2 e_2 = e_3$ &\\
&&$e_2 e_3 = e_4$&$e_3 e_1 =   ie_4$ & $e_3 e_2 = e_4$  \\

\hline
$\cd 4{69}$ & $:$ & 
$e_1 e_1 = e_3$ & $e_1 e_3 =  ie_4$  &$e_2 e_2 = e_3$ &\\
&&$e_2 e_3 = e_4$&$e_3 e_1 =   ie_4$ & $e_3 e_2 = e_4$  \\

\hline
$\cd 4{70}$ & $:$ & 
$e_1 e_1 = e_3$ & $e_1 e_3 =    e_4$  &$e_2 e_2 = e_3$ & $e_3 e_1 = e_4$ &   \\
\hline
		$\cd 4{71}$ & $:$ & $e_1 e_1 = e_4$ & $e_1 e_2 = e_3$  & $e_1 e_3 = e_4$ & $e_2 e_1 = -e_3$  & $e_3 e_1 = -e_4$\\
		\hline
		$\cd 4{72}$ & $:$ & $e_1 e_1 = e_4$ & $e_1 e_2 = e_3$ & $e_1 e_3 = e_4$\\
		&& $e_2 e_1 = -e_3$  & $e_2 e_2 = e_4$ & $e_3 e_1 = -e_4$\\
		\hline
		$\cd 4{73}$ & $:$ & $e_1 e_2 = e_3 + e_4$ & $e_1 e_3 = e_4$ & $e_2 e_1 = -e_3$  & $e_3 e_1 = -e_4$\\
		\hline
		$\cd 4{74}(\af)$ & $:$ & $e_1 e_2 = e_3$ & $e_1 e_3 = (\af+1)e_4$ & $e_2 e_1 = -e_3$ & $e_3 e_1 = -\af e_4$\\
		\hline
		$\cd 4{75}(\af)$ & $:$ & $e_1 e_2 = e_3$ & $e_1 e_3 = (\af+1)e_4$ & $e_2 e_1 = -e_3$ & $e_2 e_2 = e_4$ & $e_3 e_1 = -\af e_4$\\        \hline
		$\cd 4{76}$ & $:$ & $e_1 e_2 = e_3$ & $e_1 e_3 = e_4$ & $e_2 e_1 = -e_3$ & $e_2 e_2 = e_4$ & $e_3 e_1 = -e_4$\\        
		\hline
		$\cd 4{77}$ & $:$ & $e_1 e_2 = e_3$ & $e_1 e_3 = e_4$ & $e_2 e_1 = -e_3$ & $e_2 e_3 = e_4$ & $e_3 e_2 = -e_4$\\
		\hline
		$\cd 4{78}$ & $:$ & $e_1 e_2 = e_3$ & $e_1 e_3 = e_4$ & $e_2 e_1 = -e_3$\\
		&& $e_2 e_2 = e_4$ & $e_2 e_3 = e_4$ & $e_3 e_2 = -e_4$\\
		\hline
		$\cd 4{79}(\af)$ & $:$ & $e_1 e_2 = e_3+\af e_4$ & $e_1 e_3 = e_4$ & $e_2 e_1 = -e_3$ & $e_2 e_3 = e_4$ & $e_3 e_3 = e_4$\\
		\hline
		$\cd 4{80}$ & $:$ & $e_1 e_2 = e_3+e_4$ & $e_1 e_3 = e_4$ & $e_2 e_1 = -e_3$ & $e_3 e_3 = e_4$\\
		\hline
		$\cd 4{81}$ & $:$ & $e_1 e_2 = e_3+e_4$ & $e_2 e_1 = -e_3$ & $e_3 e_3 = e_4$\\
		\hline
		$\cd 4{82}$ & $:$ & $e_1 e_2 = e_3$ & $e_1 e_3 = e_4$ & $e_2 e_1 = -e_3$ & $e_2 e_2 = e_4$ & $e_3 e_3 = e_4$\\
		\hline
		$\cd 4{83}( \af\ne 0)$ & $:$ & $e_1 e_2 = e_3$ & $e_2 e_1 = -e_3$ & $e_2 e_2 = \af e_4$ & $e_2 e_3 = e_4$ & $e_3 e_3 = e_4$\\
		\hline
		$\cd 4{84}$ & $:$ & $e_1 e_2 = e_3$ & $e_2 e_1 = -e_3$ & $e_2 e_2 = e_4$ & $e_3 e_3 = e_4$\\
		\hline
		$\cd 4{85}$ & $:$ & $e_1 e_2 = e_3$ & $e_2 e_1 = -e_3$ & $e_3 e_3 = e_4$\\
		\hline
		$\cd 4{86}$ & $:$ & $e_1 e_2 = e_3$ & $e_2 e_1 = -e_3$ & $e_2 e_3 = e_4$  &  $e_3 e_2 = -e_4$\\
\hline		
		
$\cd {4}{87}(\lambda)$&$:$& 
$e_1 e_1 = \lambda e_3+(2 \Theta-1)e_4$& $e_1 e_2=e_4$ & $e_1e_3=e_4$& 
\multicolumn{2}{l}{$e_2 e_1=e_3-(1- \Theta)^2 \lb^{-1}e_4$}\\
$\lb \neq 0, \frac{1}{4}$&& $e_2 e_2=e_3$& $e_2e_3=\Theta\lb^{-1}e_4$ 
&$e_3e_3=e_4$\\

\hline$\cd {4}{88}(\lambda)$&$:$& 
$e_1 e_1 = \lambda e_3+(1-2 \Theta)e_4$& $e_1 e_2=e_4$ & $e_1e_3=e_4$& 
\multicolumn{2}{l}{$e_2 e_1=e_3- \Theta^2 \lb^{-1}e_4$}\\
$\lb \neq 0, \frac{1}{4}$&& $e_2 e_2=e_3$& $e_2e_3=\Theta\lb^{-1}e_4$ 
&$e_3e_3=e_4$\\

\hline$\cd {4}{89}(\lambda)$&$:$& 
$e_1 e_1 = \lambda e_3+(2 \Theta-1)e_4$& $e_1 e_2=e_4$ & $e_1e_3=e_4$& 
\multicolumn{2}{l}{$e_2 e_1=e_3-(1- \Theta)^2 \lb^{-1}e_4$}\\
$\lb \neq 0, \frac{1}{4}$&& $e_2 e_2=e_3$& $e_2e_3=(1-\Theta)\lb^{-1}e_4$ 
&$e_3e_3=e_4$\\

\hline$\cd {4}{90}(\lambda)$&$:$& 
$e_1 e_1 = \lambda e_3+(1-2 \Theta)e_4$& $e_1 e_2=e_4$ & $e_1e_3=e_4$& 
\multicolumn{2}{l}{$e_2 e_1=e_3- \Theta^2 \lb^{-1}e_4$}\\
$\lb \neq 0, \frac{1}{4}$&& $e_2 e_2=e_3$& $e_2e_3=(1-\Theta)\lb^{-1}e_4$ 
&$e_3e_3=e_4$\\

\hline$\cd {4}{91}(\lambda, \af)$&$:$& 
$e_1 e_1 = \lambda e_3+(2 \Theta-1)e_4$& $e_1 e_2=e_4$ & $e_1e_3=\af e_4$& \\
$\lb \neq 0, \frac{1}{4}$&& $e_2 e_1=e_3- (1-\Theta)^2 \lb^{-1}e_4$&
$e_2 e_2=e_3$&  $e_3e_3=e_4$\\

\hline$\cd {4}{92}(\lambda, \af)$&$:$& 
$e_1 e_1 = \lambda e_3+(1-2 \Theta)e_4$& $e_1 e_2=e_4$ & $e_1e_3=\af e_4$& \\
$\lb \neq 0, \frac{1}{4}$&& $e_2 e_1=e_3-  \Theta^2 \lb^{-1}e_4$&
$e_2 e_2=e_3$&  $e_3e_3=e_4$\\

\hline$\cd {4}{93}( \af)$&$:$& 
$e_1 e_1 = e_4$ & $e_1 e_2=e_4$ &$e_1e_3=\af e_4$ & $e_2e_1=e_3+e_4$\\
&&$ e_2 e_2=e_3$ &$e_2e_3=\af e_4$ & $e_3e_3=e_4$\\

\hline$\cd {4}{94}(\af, \bt)$&$:$& 
$e_1 e_1 = e_4$ & $e_1e_2=e_4$& $e_1e_3=\af e_4$&\\
$\af\neq0$&&$e_2 e_1=e_3+\beta e_4$   & $e_2 e_2=e_3$   &$e_3e_3=e_4$\\

\hline$\cd {4}{95}(\af)$&$:$& 
$e_1 e_1 =  e_4$ & $e_1e_2=e_4$ &$e_1e_3=\af e_4$ & $e_2 e_1=e_3$\\
&& $e_2 e_2=e_3$ & $e_2e_3=\af e_4$&$e_3e_3=e_4$\\

\hline$\cd {4}{96}(\af)$&$:$& 
$e_1 e_1 =  e_4$ & $e_1e_2=e_4$ & $e_2 e_1=e_3+\af e_4$\\
&& $e_2 e_2=e_3$ & $e_2e_3= e_4$&$e_3e_3=e_4$\\

\hline$\cd {4}{97}(\lambda)$&$:$& 
$e_1 e_1 = \lambda e_3$ & $e_1e_2=e_4$& $e_1e_3=\Theta e_4 $&$e_2 e_1=e_3-e_4$ \\
&& $e_2 e_2=e_3$ &$e_2e_3=e_4$ &$e_3e_3=e_4$\\

\hline$\cd {4}{98}(\lambda)$&$:$& 
$e_1 e_1 = \lambda e_3$ & $e_1e_2=e_4$& $e_1e_3=(1-\Theta) e_4 $&$e_2 e_1=e_3-e_4$ \\
$\lb \neq \frac{1}{4}$ && $e_2 e_2=e_3$ &$e_2e_3=e_4$ &$e_3e_3=e_4$\\

\hline$\cd {4}{99}(\af)$&$:$& 
$e_1e_2=e_4$& $e_1e_3=e_4$& $e_2 e_1=e_3-e_4$\\ 
$\af\neq1$ && $e_2 e_2=e_3$ &$e_2e_3=\af e_4$&$e_3e_3=e_4$\\

\hline$\cd {4}{100}(\af)$&$:$& 
$e_1 e_1 = \frac{1}{4} e_3$ & $e_1e_2=e_4$& $e_1e_3=\af e_4$& $e_2 e_1=e_3-e_4$ \\
$\af\notin\{0, \frac{1}{2}\}$&&  $e_2 e_2=e_3$ &$e_2e_3=2 \af e_4$ &$e_3e_3=e_4$\\

 \hline$\cd {4}{101}(\af, \bt)$&$:$& 
 $e_1e_2=e_4$& $e_1e_3=\af e_4$& $e_2 e_1=e_3$  \\
 && $e_2 e_2=e_3$&$e_2e_3=\bt e_4$ &$e_3e_3=e_4$\\
 
 \hline$\cd {4}{102}(\lb, \af)$&$:$& 
 $e_1 e_1 = \lambda e_3$ & $e_1e_2=e_4$& $e_2 e_1=e_3-e_4$ \\
 $\lb\neq0$&& $e_2 e_2=e_3$&$ e_2e_3=\af e_4$  &$e_3e_3=e_4$\\

\hline$\cd {4}{103}$&$:$& 
$e_1 e_2 = e_4$ & $e_2 e_1=e_3-e_4$  & $e_2 e_2=e_3$ &$e_3e_3=e_4$\\
 
 \hline$\cd {4}{104}$&$:$& 
 $e_1 e_3 =  e_4$ & $e_2 e_1=e_3+e_4$  & $e_2 e_2=e_3$ &$e_2e_3=e_4$&$e_3e_3=e_4$\\
 
 \hline$\cd {4}{105}(\lambda, \af,\bt)$&$:$& 
 $e_1 e_1 = \lambda e_3$ & $e_1e_3=e_4$&$e_2 e_1=e_3+\af e_4$  \\
 $  \lb\ne 0, \af\ne 0$&& $e_2 e_2=e_3$ & $e_2e_3=\bt e_4$&$e_3e_3=e_4$\\
 
 \hline$\cd {4}{106}(\af)$&$:$& $e_1 e_3 = e_4$ & $e_2 e_1=e_3+\af e_4$  & $e_2 e_2=e_3$ &$e_3e_3=e_4$\\
 
 \hline$\cd {4}{107}(\lambda)$&$:$& 
 $e_1 e_1 = \lambda e_3$ & $e_1e_3=\Theta e_4$& $e_2 e_1=e_3$  \\
 && $e_2 e_2=e_3$ & $e_2e_3=e_4$&$e_3e_3=e_4$\\
 
 \hline$\cd {4}{108}(\lambda)$&$:$& 
 $e_1 e_1 = \lambda e_3$ & $e_1e_3=(1-\Theta) e_4$& $e_2 e_1=e_3$  \\
 $\lb \not\in \{0, \frac{1}{4}\}$&& $e_2 e_2=e_3$ & $e_2e_3=e_4$&$e_3e_3=e_4$\\
 
 \hline$\cd {4}{109}(\lambda,\af)$&$:$& $e_1 e_1 = \lambda e_3$ & $e_2 e_1=e_3+e_4$  & $e_2 e_2=e_3$ &$e_2e_3=\alpha e_4$&$e_3e_3=e_4$\\

  \hline$\cd {4}{110}(\lambda)$&$:$& $e_1 e_1 = \lambda e_3$ & $e_2 e_1=e_3$  & $e_2 e_2=e_3$ &$e_2e_3= e_4$&$e_3e_3=e_4$\\
 
   \hline$\cd {4}{111}(\lambda)$&$:$& $e_1 e_1 = \lambda e_3$ & $e_2 e_1=e_3$  & $e_2 e_2=e_3$ &$e_3e_3=e_4$\\
 
\hline$\cd {4}{112}(\lambda, \af, \bt, \gamma)$&$:$& 
$e_1 e_1 = \lambda e_3+e_4$ & $e_1e_3=\af e_4$ & $e_2 e_1=e_3+\bt e_4$  \\
&&$e_2 e_2=e_3$& $e_2e_3=\gamma e_4$&$e_3e_3=e_4$\\
\end{longtable}}

All these algebras are non-isomorphic, except

{\tiny \begin{longtable}{c}
$\D{4}{01}(\lambda,0,\beta) \cong \D{4}{02}(\lambda,0,\beta) \cong \D{4}{04}(\lambda,\beta),\quad \D{4}{01}(\lambda,\alpha,0)_{\alpha \neq -1} \cong \D{4}{02}(\lambda,\alpha,0) \cong \D{4}{10}(\lambda,\alpha),\quad \D{4}{01}(\lambda,-1,0) \cong \D{4}{11}(\lambda,0),$\\

$\D{4}{03}(\lambda,0) \cong \D{4}{09}(\lambda,0),\quad \D{4}{03}\left(\lambda,(1-\Theta)^{-1}\right)_{\lambda \neq 0} \cong \D{4}{05}(\lambda,0)_{\lambda \neq 0}, \D{4}{03}\left(\lambda,\Theta^{-1}\right)\cong \D{4}{06}(\lambda,0),\quad \D{4}{04}(\lambda,0) \cong \D{4}{10}(\lambda,0),$\\

$\D{4}{05}(1/4,\alpha) \cong \D{4}{06}(1/4,\alpha),\quad \D{4}{07}(1/4) \cong \D{4}{08}(1/4),$\\

$\D{4}{05}(0,\alpha) \cong \D{4}{07}(0) \cong \D{4}{23}(0) \cong  \D{4}{25}(0) \cong  \D{4}{40}(0),$\\

$\D{4}{12}(\lambda,0) \cong \D{4}{18}(\lambda,0),\quad  \D{4}{12}(1/4,\alpha) \cong \D{4}{13}(1/4,\alpha),\quad \D{4}{12}(0,\alpha)_{\alpha \neq -1} \cong \D{4}{14}(0,\alpha),\quad \D{4}{12}(0,-1) \cong \D{4}{17}(0),$\\

$\D{4}{13}(\lambda,0) \cong \D{4}{19}(\lambda,0),\quad \D{4}{14}(\lambda,0) \cong \D{4}{20}(\lambda,0),\quad \D{4}{14}(1/4,\alpha) \cong \D{4}{15}(1/4,\alpha),\quad \D{4}{15}(\lambda,0) \cong \D{4}{21}(\lambda,0),$\\

$\D{4}{18}(1/4,\alpha) \cong \D{4}{19}(1/4,\alpha),\quad \D{4}{18}(0,0) \cong \D{4}{22}(0) \cong  \D{4}{24}(0),\quad \D{4}{18}(1/4,-1) \cong \D{4}{19}(1/4,-1) \cong \D{4}{30}(1/4) \cong \D{4}{31}(1/4),$\\

$\D{4}{20}(1/4,\alpha) \cong \D{4}{21}(1/4,\alpha),\quad \D{4}{20}(1/4,-1) \cong \D{4}{21}(1/4,-1) \cong \D{4}{32}(1/4) \cong \D{4}{33}(1/4),$\\

$\D{4}{22}(1/4) \cong \D{4}{23}(1/4) \cong \D{4}{24}(1/4) \cong \D{4}{25}(1/4) \cong \D{4}{26}(1/4) \cong \D{4}{27}(1/4) \cong \D{4}{28}(1/4) \cong \D{4}{29}(1/4),$\\

$ \D{4}{37}(1/4) \cong \D{4}{38}(1/4),\quad \D{4}{39}(1/4) \cong \D{4}{40}(1/4).$
\end{longtable}

\begin{longtable}{lll}
$\cd 4{43}(\af)\cong\cd 4{43}(-\af)$ &
$\cd 4{44}(\af,\bt,\gm)\cong\cd 4{44}(\af,-\bt,-\gm)$& 
$\cd 4{47}(\af,\bt)\cong \cd 4{47}(\af,-\bt)$\\
$\cd 4{50}(\af)=\cd 4{50}(-\af)$& 
$\cd 4{54}(\af)\cong\cd 4{54}(-\af-1)$&
$\cd 4{57}(\af,\bt)\cong \cd 4{57}(\af+\bt,-\bt)$\\
$\cd 4{59}(\af,\bt)\cong\cd 4{59}(\af,-\bt)$&&\\
$\cd {4}{91}(\lb, \af)\cong\cd {4}{91}(\lb, -\af)$ &
$\cd {4}{92}(\lb, \af)\cong\cd {4}{92}(\lb, -\af)$ & 
$\cd {4}{93}(\af)\cong\cd {4}{93}(-\af)$ \\
$\cd {4}{94}(\af,\bt)\cong\cd {4}{94}(-\af,\bt)$ &
$\cd {4}{95}(\af)\cong\cd {4}{95}(-\af)$ &
$\cd {4}{100}(\af)\cong\cd {4}{100}(-\af)$ \\
$\cd {4}{101}(\af,\bt)\cong\cd {4}{101}(-\af,-\bt)$ &
$\cd {4}{109}(\lb,\af)\cong\cd {4}{109}(\lb,-\af)$ &
$\cd {4}{112}(\lb,\af,\bt,\gm)\cong\cd {4}{112}(\lb,-\af,\bt,-\gm)$ \\
\multicolumn{3}{l}{$\cd {4}{112}(\lb,\af,\bt,\gm)\cong \cd {4}{112}\left(\lb,(\gm-\af\bt)\sqrt{\frac{-\lb}{1-\bt+\lb\bt^2}},\frac 1\lb-\bt,(\frac\gm\lb-\frac\af\lb-\bt\gm)\sqrt{\frac{-\lb}{1-\bt+\lb\bt^2}}\right)$, if $\lb\ne 0$, $\bt\ne\frac{1\pm\sqrt{1-4\lb}}{2\lb}$}
\end{longtable} 
} 
\end{theorem}

\end{document}